\documentclass[11pt]{article}

\usepackage{latexsym}

\usepackage{amsmath,amsfonts}
\usepackage{enumerate}
\usepackage{dsfont}

\font\tenbb=msbm10
\font\sevenbb=msbm7
\font\fivebb=msbm5
\newfam\bbfam
\textfont\bbfam=\tenbb
\scriptfont\bbfam=\sevenbb
\scriptscriptfont\bbfam=\fivebb
\def\Bbb{\fam\bbfam\tenbb}

\newcommand{\begeq}[1]{\begin{equation} \label{#1}}
\newcommand{\fineq}{\end{equation}}

\newcommand{\tend}{\rightarrow}
\newcommand{\g}{\gamma}

\newcommand{\al}{\alpha}

\newcommand{\be}{\beta}

\def\bR{{\Bbb R}}                  
\def\bE{{\Bbb E}}                  

\def\tinf{{\rightarrow\infty}}     

\begin{document}

\newtheorem{theor}{Theorem}
\newtheorem{prop}{Proposition}
\newtheorem{lemma}{Lemma}
\newtheorem{coro}{Corollary}
\newtheorem{defi}{Definition}
\newtheorem{rque}{Remark}
\begin{titlepage}

\title{Revisiting R\'ev\'esz's stochastic approximation method for the estimation of a regression function}
\author{Abdelkader Mokkadem \and Mariane Pelletier \and Yousri Slaoui}
\date{\small{
Universit\'e de Versailles-Saint-Quentin\\
D\'epartement de Math\'ematiques\\
45, Avenue des Etats-Unis\\
78035 Versailles Cedex\\
France\\}
(mokkadem,pelletier,slaoui)@math.uvsq.fr}
\maketitle

\begin{abstract}
In a pioneer work, R\'ev\'esz (1973) introduces the stochastic 
approximation 
method to build up a recursive kernel estimator of the regression function 
$x\mapsto \bE(Y|X=x)$. However, according to R\'ev\'esz (1977), 
his estimator has two main drawbacks: 
on the one hand, 
its convergence rate is smaller than that of the nonrecursive 
Nadaraya-Watson's kernel regression estimator, and, on the other hand, 
the required assumptions on the density of the random variable $X$ 
are stronger than those usually needed in the framework of regression 
estimation. 
We first come back on the study of the convergence rate of R\'ev\'esz's 
estimator. An approach in the proofs completely different from that used in 
R\'ev\'esz (1977) allows us to show that R\'ev\'esz's recursive estimator 
may reach the same optimal convergence rate as Nadaraya-Watson's 
estimator, but the required assumptions on the density of 
$X$ remain stronger than the usual ones, and this is inherent to the 
definition of R\'ev\'esz's estimator. 
To overcome this drawback, we introduce the averaging principle of stochastic 
approximation algorithms 
to construct the averaged R\'ev\'esz's regression estimator, and 
give its asymptotic behaviour. Our assumptions on the density of $X$ 
are then usual in the framework of regression estimation. We prove that 
the averaged R\'ev\'esz's regression estimator may reach the same 
optimal convergence rate as Nadaraya-Watson's estimator. Moreover, 
we show that, according to the estimation by confidence intervals 
point of view, it is better to use the averaged R\'ev\'esz's estimator 
rather than Nadaraya-Watson's estimator.
\end{abstract}
\vspace{1cm}

\noindent
{\bf Key words and phrases~:}\\
Regression estimation; Stochastic approximation algorithm; Averaging  principle.
\vspace{0.5cm}

\end{titlepage}

\section{Introduction}

The use of stochastic approximation algorithms in the framework of
regression estimation was introduced by Kiefer and Wolfowitz (1952). 
It allows the construction of online estimators. 
The great advantage of recursive estimators on nonrecursive ones is that their 
update, from a sample of size $n$ to one of size $n+1$, requires considerably 
less computations. This property is particularly important in the framework of 
regression estimation, since the number of points at which the 
function is estimated is usually very large. 
The famous Kiefer and Wolfowitz algorithm allows the approximation of
the point at which a regression function reaches its maximum. This
pioneer work was widely discussed and extended in many directions
(see, among many others, Blum (1954), Fabian (1967), Kushner and Clark
(1978), Hall and Heyde (1980), Ruppert (1982), Chen (1988), Spall
(1988), Polyak and Tsybakov (1990), Dippon and Renz (1997), Spall
(1997), Chen, Duncan and Pasik-Duncan (1999), Dippon (2003), and
Mokkadem and Pelletier (2007)).  The question of applying
Robbins-Monro's procedure to construct a stochastic approximation
algorithm, which allows the estimation of a regression function at a
given point (instead of approximating its mode) was introduced by
R\'ev\'esz (1973). \\

Let us recall that Robbins-Monro's procedure consists in building up
stochastic approximation algorithms, which allow the search of the
zero $z^*$ of an unknown function $h~:\bR\tend\bR$. These algorithms
are constructed in the following way~: (i) $Z_0\in\bR$ is arbitrarily
chosen; (ii) the sequence $(Z_n)$ is recursively defined by setting
\begin{equation}
\label{aRM def}
Z_n  =  Z_{n-1}+\g_n{\cal W}_n
\end{equation}
where ${\cal W}_n$ is an observation of the function $h$ at the point
$Z_{n-1}$, and where the stepsize $(\g_n)$ is a sequence of positive
real numbers that goes to zero.\\
 
Let $(X,Y), (X_1,Y_1),\ldots ,(X_n,Y_n)$ be independent, identically
distributed pairs of random variables, and let $f$ denote the
probability density of $X$. In order to construct a stochastic
algorithm for the estimation of the regression function $r:x\mapsto
\bE(Y|X=x)$ at a point $x$ such that $f(x)\neq 0$, R\'ev\'esz (1973)
defines an algorithm, which approximates the zero of the function
$h~:y\mapsto f(x)r(x)-f(x)y$. Following Robbins-Monro's procedure,
this algorithm is defined by setting $r_0(x)\in\bR$ and, for $n\geq
1$,
\begin{equation}
\label{arevesz}
r_n(x)  =  r_{n-1}(x)+\frac{1}{n}{\cal W}_n(x)
\end{equation}
where ${\cal W}_n(x)$ is an ``observation'' of the function $h$ at the
point $r_{n-1}(x)$. To define ${\cal W}_n(x)$, R\'ev\'esz (1973)
introduces a kernel $K$ (that is, a function satisfying $\int_\bR
K(x)dx=1$) and a bandwidth $(h_n)$ (that is, a sequence of positive
real numbers that goes to zero), and sets
\begin{equation}
\label{a1intro}
{\cal W}_n(x)=h_n^{-1}Y_nK(h_n^{-1}[x-X_n])-h_n^{-1}K(h_n^{-1}[x-X_n])r_{n-1}(x).
\end{equation}
R\'ev\'esz (1977) chooses the bandwidth $(h_n)$ equal to $(n^{-a})$
with $a\in]1/2,1[$, and establishes a central limit theorem for
$r_n(x)-r(x)$ under the assumption $f(x)> (1-a)/2$, as well as an
upper bound of the uniform strong convergence rate of $r_n$ on any
bounded interval $I$ on which $\inf_{x\in I}f(x)> (1-a)/2$. The two
drawbacks of his approach are the following.  First, the condition
$a>1/2$ on the bandwidth leads to a pointwise weak convergence rate of
$r_n$ slower than $n^{1/4}$, whereas the optimal pointwise weak
convergence rate of the kernel estimator of a regression function
introduced by Nadaraya (1964) and Watson (1964) is $n^{2/5}$ (and
obtained by choosing $a=1/5$). Then, the condition $f(x)> (1-a)/2$ (or
$\inf_{x\in I}f(x)> (1-a)/2$) is stronger than the condition $f(x)>0$
(or $\inf_{x\in I}f(x)>0$) usually required to establish the
convergence rate of regression's estimators.\\

Our first aim in this paper is to come back on the study of the
asymptotic behaviour of R\'ev\'esz's estimator. The technic we use is
totally different from that employed by R\'ev\'esz (1977).  Noting
that the estimator $r_n$ can be rewritten as
\begin{eqnarray*}
r_n(x) & = & 
\left(1-\frac{1}{n}h_n^{-1}K\left(\frac{x-X_n}{h_n}\right)\right)r_{n-1}(x)
+\frac{1}{n}h_n^{-1}Y_nK\left(\frac{x-X_n}{h_n}\right) \\
& = & 
\left(1-\frac{1}{n}f(x)\right)r_{n-1}(x)
+\frac{1}{n}\left[f(x)-h_n^{-1}K\left(\frac{x-X_n}{h_n}\right)\right]r_{n-1}(x)
+\frac{1}{n}h_n^{-1}Y_nK\left(\frac{x-X_n}{h_n}\right),
\end{eqnarray*}
we approximate the sequence $(r_n)$ by the unobservable sequence
$(\rho_n)$ recursively defined by
\begin{equation}
\label{a4intro}
\rho_n(x) = 
\left(1-\frac{1}{n}f(x)\right)\rho_{n-1}(x)
+\frac{1}{n}\left[f(x)-h_n^{-1}K\left(\frac{x-X_n}{h_n}\right)\right]r(x)
+\frac{h_n^{-1}}{n}Y_nK\left(\frac{x-X_n}{h_n}\right).
\end{equation}
The asymptotic properties (pointwise weak and strong convergence rate,
upper bound of the uniform strong convergence rate) of the
approximating algorithm (\ref{a4intro}) are established by applying
different results on the sums of independent variables and on the
martingales. To show that the asymptotic properties of the
approximating algorithm (\ref{a4intro}) are also satisfied by
R\'ev\'esz's estimator, we use a technic of successive upper bounds.
It turns out that our technique of demonstration allows the choice of
the bandwidth $(h_n)=(n^{-1/5})$, which makes R\'ev\'esz's estimator
converge with the optimal pointwise weak convergence rate
$n^{2/5}$. However, to establish the asymptotic convergence rate of
R\'ev\'esz's estimator, we need the same kind of conditions on the
marginal density of $X$ as R\'ev\'esz (1977) does.  To understand why
this second drawback is inherent in the definition of R\'ev\'esz's
estimator, let us come back on the algorithm (\ref{aRM def}). 
The convergence rate of stochastic approximation algorithms
constructed following Robbins-Monro's scheme, and used for the
search of the zero $z^*$ of an unknown function $h$, was widely
studied. It is now well known (see, among many others, Nevels'on and Has'minskii (1976),
Kushner and Clark (1978),
Ljung, Pflug, and Walk (1992), and Duflo (1996)) 
that the convergence rate of algorithms defined as (\ref{aRM def}) is
obtained under the condition that the limit of the sequence $(n\g_n)$
as $n$ goes to infinity is larger than a quantity, which involves the
differential of $h$ at the point $z^*$.  
Now, let us recall that the  R\'ev\'esz's estimator (\ref{arevesz})
is an algorithm approximating the zero $y^*=r(x)$ of the function 
$y\mapsto f(x)r(x)-f(x)y$ (whose differential at the point $y^*$  
equals $-f(x)$), and let us enlighten that the stepsize used to define his 
algorithm is $(\g_n)=(n^{-1})$ (so that $\lim_{n\rightarrow\infty}n\g_n=1$); 
the condition on the limit of $(n\g_n)$, which is usual in the framework of 
stochastic approximation algorithms, comes down, in the case of R\'ev\'esz's 
estimator, to a condition on the probability density $f$. \\

Our second aim in this paper is to introduce the averaging principle
of stochastic approximation algorithms in the framework of regression
estimation. As a matter of fact, in the framework of stochastic approximation, 
this principle is now well known to allow to get rid of tedious conditions 
on the stepsize. It was independently introduced
by Ruppert (1988) and Polyak (1990), and then widely discussed and
extended (see, among many others, Yin (1991), Delyon and Juditsky
(1992), Polyak and Juditsky (1992), Kushner and Yang (1993), Le Breton
(1993), Le Breton and Novikov (1995), Dippon and Renz (1996, 1997),
and Pelletier (2000)).  This procedure
consists in: (i) running the approximation algorithm by using 
slower stepsizes; (ii) computing an average of the approximations 
obtained in (i). We thus need to generalize the definition of R\'ev\'esz's 
estimator before defining the averaged R\'ev\'esz's estimator. 
\\

Let $(\g_n)$ be a sequence of positive numbers going to zero.
The generalized R\'ev\'esz's estimator is defined by setting
$r_0(x)\in\bR$, and, for $n\geq 1$,
\begin{equation}
\label{arun}
r_n(x)  =  r_{n-1}(x)+\g_n{\cal W}_n(x),
\end{equation}
where ${\cal W}_n(x)$ is defined in (\ref{a1intro}).  (R\'ev\'esz's
estimator clearly corresponds to the choice of the stepsize
$(\g_n)=(n^{-1})$).  The estimator (\ref{arun}) with $(\g_n)$ not necessary 
equal to $(n^{-1})$ was introduced in Gy\"orfi et al. (2002), where the strong universal 
convergence rate of $r_n(x)$ is also proved.
Now, let the stepsize in (\ref{arun}) satisfy
$\lim_{n\tinf}n\g_n=\infty$, and let $(q_n)$ be a positive sequence
such that $\sum q_n=\infty$.  The averaged R\'ev\'esz's estimator is
defined by setting
\begin{equation}
\label{a3intro}
\overline r_n(x) = \frac{1}{\sum_{k=1}^nq_k}\sum_{k=1}^nq_kr_k(x)
\end{equation}
(where the $r_k(x)$ are given by the algorithm (\ref{arun})). 
Let us enlighten that the estimator $\overline r_n$ is still recursive. 
\\

We establish the asymptotic behaviour (pointwise weak and strong
convergence rate, upper bound of the uniform strong convergence rate)
of $\overline r_n$.  The condition we require on the density $f$ to
prove the pointwise (respectively uniform) convergence rate of
the averaged R\'ev\'esz's estimator is the usual condition $f(x)>0$
(respectively $\inf_{x\in I}f(x)>0$). Concerning the bandwidth, the
choice $(h_n)=(n^{-1/5})$, which leads to the optimal convergence rate
$n^{2/5}$, is allowed.  Finally, we show that to construct confidence
intervals by slightly undersmoothing, it is preferable to use the
averaged R\'ev\'esz's estimator $\overline r_n$ (with an adequate choice
of weights $(q_n)$) rather than Nadaraya-Watson's estimator, since
the asymptotic variance of this lattest estimator is larger than that of $\overline r_n$.\\

Our paper is organized as follows. Our assumptions and main results
are stated in Section \ref{asection 2}, simulations study is performed
in Section \ref{asection simul}, the outlines of the proofs given in
Section \ref{asection 3}, whereas Section \ref{alem} is devoted to the
proof of several lemmas.

\section{Assumptions and main results} \label{asection 2}
Let us first define the class of positive sequences that will be used
in the statement of our assumptions.

\begin{defi}
Let $\gamma \in \mathbb{R} $ and $\left(v_n\right)_{n\geq 1}$ be a nonrandom positive sequence. We say that $\left(v_n\right) \in \mathcal{GS}\left(\gamma \right)$ if
\begin{eqnarray}\label{aGS}
\lim_{n \to \infty} n\left[1-\frac{v_{n-1}}{v_{n}}\right]=\gamma.
\end{eqnarray}
\end{defi}
Condition~\eqref{aGS} was introduced by Galambos and Seneta (1973) to define regularly varying sequences (see also Bojanic and Seneta (1973)); it was used in 
Mokkadem and Pelletier (2007) in the context of stochastic approximation algorithms. Typical sequences in $\mathcal{GS}\left(\gamma\right)$ are, for $b\in \mathbb{R}$, $n^{\gamma}\left(\log n\right)^{b}$, $n^{\gamma}\left(\log \log n\right)^{b}$, and so on.\\

Let $g\left(s,t\right)$ denote the density of the couple
$\left(X,Y\right)$ (in particular
$f\left(x\right)=\int_{\mathbb{R}}g\left(x,t\right)dt$), and set
$a\left(x\right)=r\left(x\right)f\left(x\right)$. The assumptions to
which we shall refer for our pointwise results are the following.
\begin{description}
\item(A1) $K:{\mathbb R}\rightarrow {\mathbb R}$ is a nonnegative, continuous, bounded function satisfying $\int_{\mathbb{R}}K\left( z\right) dz=1$, $\int_{\mathbb{R}}zK\left( z\right) dz=0$ and $\int_{\mathbb{R}}z^2K\left( z\right) dz<\infty$.
 \item(A2) $i)$ $\left(\gamma_n\right) \in \mathcal{GS}\left(-\alpha \right)$ with $\alpha\in \left]\frac{3}{4},1\right]$; moreover the limit of $\left(n\gamma_n\right)^{-1}$ as $n$ goes to infinity exists.\\
$ii)$ $\left( h_{n}\right)\in \mathcal{GS} \left(-a\right)$ with $ a \in \left]\frac{1-\alpha}{4},\frac{\alpha}{3}\right[$.
\item(A3) $i)$ $g\left(s,t\right)$ is two times continuously differentiable with respect to $s$.\\
 $ii)$ For $q\in\left\{0,1,2\right\}$, $s \mapsto \int_{\mathbb{R}}t^qg\left(s,t\right)dt$ is a bounded function continuous at $s=x$.\\ 
 For $q\in \left[2,3\right]$, $s \mapsto \int_{\mathbb{R}}\left|t\right|^qg\left(s,t\right)dt$ is a bounded function.\\
$iii)$ For $q\in \left\{0,1\right\}$, $\int_{\mathbb{R}}\left|t\right|^q\left|\frac{\partial g}{\partial x}\left(x,t\right)\right|dt<\infty$, and $s\mapsto \int_{\mathbb{R}}t^q\frac{\partial^2 g}{\partial s^2}\left(s,t\right)dt$ is a bounded function continuous at $s=x$.
\end{description}
\begin{rque} (A3) $ii)$ says in particular that $f$ is continuous at $x$.
\end{rque}

For our uniform results, we shall also need the following additional assumption.
\begin{description}
\item $\left(A4\right)$ $i)$$K$ is Lipschitz-continuous.\\
$ii)$ There exists $t^*>0$ such that $\mathbb{E}\left(\exp\left(t^*\left|Y\right|\right)\right) <\infty$.\\
$iii)$ $a\in\left]1-\alpha,\alpha-2/3\right[$.\\
$iv)$ For $q\in \left\{0,1\right\}$, $x \mapsto \int_{\mathbb{R}}\left|t\right|^q\left|\frac{\partial g}{\partial x}\left(x,t\right)\right|dt$ is bounded on the set $\left\{x,f\left(x\right)>0\right\}$.
\end{description}
Throughout this paper we shall use the following notation :
\begin{eqnarray}
\xi &= &\lim_{n\to \infty}\left(n\gamma_n\right)^{-1},\label{axi1p}
\end{eqnarray}
and, for $f\left(x\right)\not=0$,
\begin{eqnarray*}
m^{\left(2\right)}\left(x\right)&=&\frac{1}{2f\left(x\right)}\left[\int_{\mathbb{R}}t\frac{\partial^2 g}{\partial x^2}\left(x,t\right)dt-r\left(x\right)\int_{\mathbb{R}}\frac{\partial^2 g}{\partial x^2}\left(x,t\right)dt\right]\int_{\mathbb{R}}z^2K\left(z\right)dz\nonumber .
\end{eqnarray*}
The asymptotic properties of the generalized R\'ev\'esz's estimator defined in~\eqref{arun} are stated in Section \ref{asection 2.1}, those of the averaged estimator defined in~\eqref{a3intro} in Section \ref{asection 2.2}.
\subsection{Asymptotic behaviour of the generalized R\'ev\'esz's estimator} \label{asection 2.1}

For stepsizes satisfying $\lim_{n\to\infty}n\gamma_n=\infty$, the strong 
universal consistency of the generalized R\'ev\'esz's estimator was 
established by Walk (2001) and Gy\"orfi et al. (2002). 
The aim of this section is to state the convergence rate of the
estimator defined by (\ref{arun}). Theorems
\ref{aP:ra}, \ref{aT:rb}, and \ref{aT:c} below give its weak pointwise
convergence rate, its strong pointwise convergence rate, and an upper
bound of its strong uniform convergence rate, respectively. Let us 
enlighten that the particular 
choice of stepsize $(\g_n)=(n^{-1})$ gives the asymptotic behaviour of
R\'ev\'esz's estimator defined in (\ref{arevesz}).

\begin{theor}[Weak pointwise convergence rate of $r_n$]\label{aP:ra} $ $\\
Let Assumptions $\left( A1\right)-\left( A3\right) $ hold for $x \in
\mathbb{R}$ such that $f\left(x\right)\not=0$.
\begin{enumerate}
\item If there exists $c \geq 0 $ such that
$\gamma_n^{-1}h_n^5\rightarrow c$, and if $\lim_{n\to
\infty}\left(n\gamma_n\right)>\left(1-a\right)/\left(2f\left(x\right)\right)$,
then
\begin{eqnarray*}
\lefteqn{\sqrt{\gamma_n^{-1}h_n}\left( r_{n}\left( x\right)-r\left( x\right) \right)}\nonumber\\
&&\stackrel{\mathcal{D}}{\rightarrow}\mathcal{N}\left(\frac{\sqrt{c}f\left(x\right)m^{\left(2\right)}\left(x\right)}{f\left(x\right)-2a\xi},\frac{Var\left[Y\vert X=x\right]f\left(x\right)}{\left(2f\left(x\right)-\left(1-a\right)\xi \right)}\int_{\mathbb{R}} K^2\left(z\right)dz\right).
\end{eqnarray*}
\item If  $\gamma_n^{-1}h_n^5\rightarrow \infty $, and if $\lim_{n\to \infty}\left(n\gamma_n\right)>2a/f\left(x\right)$, then 
\begin{eqnarray*}
\frac{1}{h_{n}^{2}}\left( r_{n}\left( x\right)-r\left( x\right) \right) \stackrel{\mathbb{P}}{\rightarrow }\frac{f\left(x\right)m^{\left(2\right)}\left(x\right)}{\left(f\left(x\right)-2a \xi\right)} ,
\end{eqnarray*}
\end{enumerate}
where $\stackrel{\mathcal{D}}{\rightarrow}$ denotes the convergence in distribution, $\mathcal{N}$ the Gaussian-distribution and $\stackrel{\mathbb{P}}{\rightarrow}$ the convergence in probability.
\end{theor}

The combination of Parts 1 and 2 of Theorem~\ref{aP:ra} ensures that the optimal weak pointwise 
convergence rate of $r_n$ equals $n^{2/5}$, and is obtained by 
choosing $a=1/5$ and $(\g_n)$ such that $\lim_{n\to \infty}\left(n\gamma_n\right)\in]2/(5f\left(x\right)),\infty[$.

\begin{theor}[Strong pointwise convergence rate of $r_n$]\label{aT:rb} $ $\\
Let Assumptions $\left( A1\right)-\left( A3\right) $ hold for $x \in \mathbb{R}$ such that $f\left(x\right)\not=0$.
\begin{enumerate}
\item If there exists $c \geq 0 $ such that $\gamma_n^{-1}h_n^5/\ln\left(\sum_{k=1}^n \gamma_k\right) \rightarrow c$, and if $\lim_{n\to \infty}\left(n\gamma_n\right)>\left(1-a\right)/\left(2f\left(x\right)\right)$, then, with probability one, the sequence
\[
\left(\sqrt{\frac{\gamma_n^{-1}h_n}{2 \ln\left(\sum_{k=1}^n \gamma_k\right)}}\left( r_{n}\left( x\right)-r\left( x\right) \right) \right) 
\]
is relatively compact and its limit set is the interval
\begin{eqnarray*}
\left[\sqrt{\frac{c}{2}}\frac{f\left(x\right)m^{\left(2\right)}\left( x\right)}{f\left(x\right)-2a\xi}-\sqrt{\frac{Var\left[Y\vert X=x\right]f\left(x\right)\int_{\mathbb{R}} K^2\left(z\right)dz}{\left(2f\left(x\right)-\left(1-a\right)\xi\right)}},\right.\\
\left.\sqrt{\frac{c}{2}}\frac{f\left(x\right)m^{\left(2\right)}\left( x\right)}{f\left(x\right)-2a\xi}+\sqrt{\frac{Var\left[Y\vert X=x\right]f\left(x\right)\int_{\mathbb{R}} K^2\left(z\right)dz}{\left(2f\left(x\right)-\left(1-a\right)\xi\right)}}\right].
\end{eqnarray*}
\item If $\gamma_n^{-1}h_n^5/\ln\left(\sum_{k=1}^n \gamma_k\right)\rightarrow \infty $, and if $\lim_{n\to \infty}\left(n\gamma_n\right)>2a/f\left(x\right)$, then, with probability one,
\begin{eqnarray*}
\lim_{n\rightarrow \infty}\frac{1}{h_{n}^{2}}\left( r_{n}\left( x\right)-r\left( x\right) \right) = \frac{f\left(x\right)m^{\left(2\right)}\left( x\right)}{f\left(x\right)-2a\xi} . 
\end{eqnarray*}
\end{enumerate}
\end{theor}

\begin{theor}[Strong uniform convergence rate of $r_n$]\label{aT:c} $ $\\
Let $I$ be a bounded open interval of $\mathbb{R}$ on which $\varphi=\inf_{x\in I}f\left(x\right)>0$, and let Assumptions $\left( A1\right)-\left( A4\right) $ hold for all $x \in I$. 
\begin{enumerate}
\item If the sequence $\left(\gamma_n^{-1}h_n^5/[\ln n]^2\right)$ is bounded and if $\lim_{n\to\infty}\left(n\gamma_n\right)>\left(1-a\right)/\left(2\varphi\right)$, then
\begin{eqnarray*}
\sup_{x \in I}\left|r_n\left(x\right)-r\left(x\right)\right|=O\left(\sqrt{\gamma_nh_n^{-1}}\, \ln n \right) \quad a.s.
\end{eqnarray*}
\item If $\lim_{n\to\infty}\left(\gamma_n^{-1}h_n^5/[\ln n]^2\right)=\infty$ and if $\lim_{n\to\infty}\left(n\gamma_n\right)>2a/\varphi$, then
\begin{eqnarray*}
\sup_{x \in I}\left|r_n\left(x\right)-r\left(x\right)\right|=O\left(h_n^2 \right) \quad a.s.
\end{eqnarray*}
\end{enumerate}
\end{theor}

Parts 1 of Theorems \ref{aP:ra} and \ref{aT:c} were obtained by
R\'ev\'esz (1977) for the choices $(\g_n)=(n^{-1})$ and
$(h_n)=(n^{-a})$ with $a\in]1/2,1[$.  Let us underline that, for this
choice of stepsize, the conditions $\lim_{n\to
\infty}\left(n\gamma_n\right)>\left(1-a\right)/\left(2f\left(x\right)\right)$
and $\lim_{n\to \infty}\left(n\gamma_n\right)>\left(1-a\right)/
\left(2\inf_{x\in I}f\left(x\right)\right)$ come down to the following
conditions on the unknown density $f$: $f(x)>(1-a)/2$ and $\inf_{x\in
I}f(x)>(1-a)/2$. Let us also mention that our assumption (A2) implies
$a\in]0,1/3[$, so that our results on the generalized R\'ev\'esz's
estimator do not include the results given in R\'ev\'esz (1977).
However, our assumptions include the choice $(\g_n)=(n^{-1})$ and
$a=1/5$, which leads to the optimal weak convergence rate $n^{2/5}$,
whereas the condition on the bandwidth required by R\'ev\'esz leads to
a convergence rate of $r_n$ slower than $n^{1/4}$. \\

Although the optimal convergence rate we obtain for the generalized
R\'ev\'esz's estimator $r_n$ has the same order as that of 
Nadaraya-Watson's estimator, this estimator has a main drawback: to make
$r_n$ converge with its optimal rate, one must set $a=1/5$ and choose
$(\g_n)$ such that $\lim_{n\tend\infty}n\g_n=\g^*\in]0,\infty[$ with
$\g^*> 2/[5f(x)]$ whereas the density $f$ is unknown. This tedious
condition disappears as soon as the stepsize is chosen such that
$\lim_{n\tend\infty}n\g_n=\infty$ (for instance when $(\g_n)=((\ln
n)^bn^{-1})$ with $b>0$), but the optimal convergence rate $n^{2/5}$
is then not reached any more.

\subsection{Asymptotic behaviour of the averaged R\'ev\'esz's estimator} \label{asection 2.2}
To state the asymptotic properties of the averaged R\'ev\'esz's
estimator defined in (\ref{a3intro}), we need the following additional
assumptions.
\begin{description}
 \item(A5) $\lim_{n\to \infty}n\gamma_n\left(\ln \left(\sum_{k=1}^n \gamma_k\right)\right)^{-1}=\infty$ and $a\in\left]1-\alpha,\left(4\alpha-3\right)/2\right[$.
 \item(A6) $\left(q_n\right) \in \mathcal{GS}\left(-q \right)$ with $q<\min\left\{1-2a,\left(1+a\right)/2\right\}$.
\end{description}
Theorems~\ref{aAv:ra},~\ref{aT:rba} and~\ref{aT:Ac} below give the weak
pointwise convergence rate, the strong pointwise convergence rate, and
an upper bound of the strong uniform convergence rate of the averaged
R\'ev\'esz's estimator.
\begin{theor}[Weak pointwise convergence rate of $\bar{r}_n$]\label{aAv:ra} $ $\\
Let Assumptions $\left( A1\right)-\left( A3\right)$, $\left(A5\right)$
and $\left(A6\right)$ hold for $x \in \mathbb{R}$ such that
$f\left(x\right)\not=0$.
\begin{enumerate}

\item If there exists $c \geq 0 $ such that  $nh_n^5\rightarrow  c$, then
\begin{eqnarray*}
\lefteqn{\sqrt{nh_n}\left( \bar{r}_{n}\left( x\right)-r\left( x\right) \right)}\\
&&\stackrel{\mathcal{D}}{\rightarrow}\mathcal{N}\left(c^{\frac{1}{2}}\frac{1-q}{1-q-2a}m^{\left(2\right)}\left(x\right),\frac{\left(1-q\right)^2}{1+a-2q}\frac{Var\left[Y\vert X=x\right]}{f\left(x\right)}\int_{\mathbb{R}} K^2\left(z\right)dz\right).
\end{eqnarray*}
\item If  $nh_n^5\rightarrow \infty $, then  
\begin{eqnarray*}
\frac{1}{h_{n}^{2}}\left( \bar{r}_{n}\left( x\right)-r\left( x\right) \right) \stackrel{\mathbb{P}}{\rightarrow }\frac{1-q}{1-q-2a}m^{\left(2\right)}\left(x\right).
\end{eqnarray*}
\end{enumerate}
\end{theor}

The combination of Parts 1 and 2 of Theorem \ref{aAv:ra} ensures that the 
optimal weak pointwise convergence rate of $\bar{r}_n$ is obtained by choosing $a=1/5$, and equals $n^{2/5}$.

\begin{theor}[Strong pointwise convergence rate of $\bar{r}_n$]\label{aT:rba} $ $\\
Let Assumptions $\left( A1\right)-\left( A3\right)$, $\left(A5\right)$ and $\left(A6\right)$ hold for $x\in \mathbb{R}$ such that $f\left(x\right)\not=0$.
\begin{enumerate}

\item If there exists $c_{1} \geq 0 $ such that $nh_n^5/\ln \ln n\rightarrow c_{1}$, then, with probability one, the sequence
\[
\left(\sqrt{\frac{nh_n}{2\ln \ln n}}\left( \bar{r}_{n}\left( x\right)-r\left( x\right) \right) \right) 
\]
is relatively compact and its limit set is the interval
\begin{eqnarray*}
\lefteqn{\left[c_1^{\frac{1}{2}}\frac{1-q}{1-q-2a}m^{\left(2\right)}\left( x\right)-\sqrt{\frac{\left(1-q\right)^2}{1+a-2q}\frac{Var\left[Y/X=x\right]}{f\left(x\right)}\int_{\mathbb{R}} K^2\left(z\right)dz},\right.}\\
&&\left.c_1^{\frac{1}{2}}\frac{1-q}{1-q-2a}m^{\left(2\right)}\left( x\right)+\sqrt{\frac{\left(1-q\right)^2}{1+a-2q}\frac{Var\left[Y/X=x\right]}{f\left(x\right)}\int_{\mathbb{R}} K^2\left(z\right)dz}\right].
\end{eqnarray*}
\item If $nh_n^5/\ln \ln n \rightarrow \infty $, then  
\begin{eqnarray*}
\lim_{n\rightarrow \infty}\frac{1}{h_{n}^{2}}\left( \bar{r}_{n}\left( x\right)-r\left( x\right) \right) = \frac{1-q}{1-q-2a}m^{\left(2\right)}\left( x\right)\quad a.s. 
\end{eqnarray*}
\end{enumerate}
\end{theor}

\begin{theor}[Strong uniform convergence rate of $\bar{r}_n$]\label{aT:Ac} $ $\\
Let $I$ be a bounded open interval of $\mathbb{R}$ on which $\inf_{x\in I}f\left(x\right)>0$, and let Assumptions $\left( A1\right)-\left( A6\right) $ hold for all $x\in I$. 
\begin{enumerate}
\item If the sequence $\left(nh_n^5/[\ln n]^2\right)$ is bounded, and if 
$\alpha>\left(3a+3\right)/4$, then 
\begin{eqnarray*}
\sup_{x \in I}\left|\bar{r}_n\left(x\right)-r\left(x\right)\right|=O\left(\sqrt{n^{-1}h_n^{-1}} \ln n\right) \quad a.s.
\end{eqnarray*}
\item If $\lim_{n\to\infty}\left(nh_n^5/[\ln n]^2\right)=\infty$, and if, in the case 
$a\in[\al/5,1/5]$, $\alpha>\left(4a+1\right)/2$, then 
\begin{eqnarray*}
\sup_{x \in I}\left|\bar{r}_n\left(x\right)-r\left(x\right)\right|=O\left(h_n^2\right) \quad a.s.
\end{eqnarray*}
\end{enumerate}
\end{theor}

Whatever the choices of the stepsize $\left(\gamma_n\right)$ and of the 
weight $\left(q_n\right)$ are, the convergence rate of the averaged 
R\'ev\'esz's estimator has the same order as that of the generalized 
R\'ev\'esz's estimator defined with a stepsize $\left(\gamma_n\right)$ 
satisfying $\lim_{n\tend \infty}(n\g_n)^{-1}\neq 0$ (and, in particular, 
of R\'ev\'esz's estimator). The main advantage of 
the averaged R\'ev\'esz's estimator on its nonaveraged version is that 
its convergence rate is obtained without tedious conditions on the 
marginal density $f$.\\

Now, to compare the averaged R\'ev\'esz's estimator with the 
nonrecursive Nadaraya-Watson's estimator defined as
\begin{eqnarray*}
\tilde{r}_{n}\left(x\right)=\frac{\sum_{k=1}^nY_kK\left(h_n^{-1}[x-X_k]\right)}
{\sum_{k=1}^nK\left(h_n^{-1}[x-X_k]\right)},
\end{eqnarray*}
let us consider the estimation by confidence intervals point of view. 
In the context of density estimation, 
Hall (1992) shows that, to construct confidence intervals, slightly
undersmoothing is more efficient than bias estimation; in the framework of regression estimation, 
the method of 
undersmoothing to construct 
confidence regions is used in particular by 
Neumann and Polzehl (1998) and Claeskens and Van Keilegom (2003).
To undersmooth,
we choose $\left(h_n\right)$ such that $\lim_{n\tend\infty}nh_n^5=0$
(and thus $a\geq 1/5$). Moreover, to construct a confidence interval
for $r(x)$, it is advised to choose the weight $\left(q_n\right)$,
which minimizes the asymptotic variance of $\bar{r}_{n}$. For a given
$a$, the function $q\mapsto (1-q)^2/(1+a-2q)$ reaching its minimum at
the point $q=a$, we can state the following corollary.

\begin{coro}\label{aC:ra} $ $\\
Let Assumptions $\left( A1\right)-\left( A3\right)$, $\left(
A5\right)$ and $\left(A6\right)$ hold for $x \in \mathbb{R}$ such that
$f\left(x\right)\not=0$, and with $a\geq 1/5$. To minimize the
asymptotic variance of $\bar{r}_n$, $q$ must be chosen equal to $a$.
Moreover, if $\lim_{n\tend\infty}nh_n^5=0$, we then have
$$
\sqrt{nh_n}\left( \bar{r}_{n}\left( x\right)-r\left( x\right) \right)
\stackrel{\mathcal{D}}{\rightarrow}\mathcal{N}\left(0
,\left(1-a\right)\frac{Var\left[Y\vert X=x\right]}{f\left(x\right)}\int_{\mathbb{R}} K^2\left(z\right)dz\right).
$$
\end{coro}

Let us recall that, when the bandwidth $(h_n)$ is chosen such that
$\lim_{n\tend\infty}nh_n^5=0$, Nadaraya-Watson's estimator satisfies the 
central limit theorem 
\begin{equation}
\label{aN-W}
\sqrt{nh_n}\left( \tilde{r}_{n}\left( x\right)-r\left( x\right) \right)
\stackrel{\mathcal{D}}{\rightarrow}\mathcal{N}\left(0
,\frac{Var\left[Y\vert X=x\right]}{f\left(x\right)}\int_{\mathbb{R}} K^2\left(z\right)dz\right).
\end{equation}
It turns out that the averaged R\'ev\'esz's estimator defined with a 
weight $(q_n)$ belonging to ${\cal GS}(-a)$ has a smaller asymptotic 
variance than Nadaraya-Watson's estimator. According to 
the estimation by confidence intervals point of view, it is thus better 
to use the averaged R\'ev\'esz's estimator rather than  Nadaraya-Watson's one. 
This superiority of the recursive averaged R\'ev\'esz's estimator on the 
classical nonrecursive Nadaraya-Watson's estimator must be related to that 
of recursive density estimators on the classical nonrecursive Rosenblatt's 
estimator, and can be explained more easily in the framework of density 
estimation: roughly speaking, Rosenblatt's estimator can be seen as the 
average of $n$ independent random variables, which all share the same variance $v_n$, whereas 
recursive estimators appear as the average of $n$ independent random variables whose variances $v_k^*$, $1\leq k\leq n$, satisfy $v_k^*<v_n$ for all $k<n$ and $v_n^*=v_n$.

\section{Simulations}\label{asection simul}
The object of this section is to provide a simulations study comparing 
Nadaraya-Watson's estimator and the averaged R\'ev\'esz's estimator. 
We consider the regression model
\begin{eqnarray*}
Y=r\left(X\right)+d\varepsilon
\end{eqnarray*}
where $d>0$ and $\varepsilon$ is $\mathcal{N}\left(0,1\right)$-distributed. 
Whatever the estimator is, we choose the kernel 
$K\left(x\right)=\left(2\pi\right)^{-1/2}\exp\left(-x^2/2\right)$, 
and the bandwidth equal to $\left(h_n\right)=n^{-1/5}\left(\ln n\right)^{-1}$ 
(which corresponds to a slight undersmoothing). The confidence intervals 
of $r\left(x\right)$ we consider are the following.
\begin{itemize}
\item When Nadaraya-Watson's estimator $\tilde{r}_n$ is used, we set
\begin{eqnarray*}
\tilde{I}_n&=&\left[\tilde{r}_n\left(x\right)-1.96\sqrt{\frac{\sum_{i=1}^n\left(Y_i-\tilde{r}_n\left(X_i\right)\right)^2}{n^2h_n\tilde{f}_n\left(x\right)}\int_{\mathbb{R}}K^2\left(z\right)dz},\right.\\
&&\left.\,\tilde{r}_n\left(x\right)+1.96\sqrt{\frac{\sum_{i=1}^n\left(Y_i-\tilde{r}_n\left(X_i\right)\right)^2}{n^2h_n\tilde{f}_n\left(x\right)}\int_{\mathbb{R}}K^2\left(z\right)dz}\right],
\end{eqnarray*}
where $\tilde{f}_n\left(x\right)=(nh_n)^{-1}\sum_{k=1}^n
K\left(h_n^{-1}[x-X_k]\right)$  
is Rosenblatt's density estimator. In view of~(\ref{aN-W}), the 
asymptotic confidence level of $\tilde{I}_n$ is $95\%$.
\item To define the averaged R\'ev\'esz's estimator $\bar{r}_n$, 
we choose the weights $\left(q_n\right)$ equal to $\left(h_n\right)$. 
This choice guarantees that $\left(h_n\right)$ and $\left(q_n\right)$ 
both belong to $\mathcal{GS}\left(-a\right)$ (with $a=1/5$ here), so that, 
in view of Corollary~\ref{aC:ra}, the asymptotic variance of $\bar{r}_n$ 
is minimal. We also let $\left(\gamma_n\right)=\left(n^{-0.9}\right)$ 
(this choice being allowed by our assumptions). Moreover, we estimate the 
density $f$ by the recursive estimator  
$\hat{f}_n\left(x\right)$ defined as 
\begin{equation*}
\hat f_n(x)=(1-\be_n)\hat f_{n-1}(x)+\be_nh_n^{-1}
K\left(\frac{x-X_n}{h_n}\right),
\end{equation*}
where $\left(\be_n\right)=\left(\frac{4}{5}n^{-1}\right)$; this choice 
of the stepsize $\left(\be_n\right)$ is known to minimize the variance of $\hat f_n$ 
(see Mokkadem et al. (2008)). 
Finally, replacing $\tilde{r}_n$ and $\tilde{f}_n$ by the recursive 
estimators $\bar{r}_n$ and 
$\hat{f}_n$ in the definition of $\tilde{I}_n$, we get the recursive 
confidence interval 
\begin{eqnarray*}
\bar{I}_n&=&\left[\bar{r}_n\left(x\right)-1.96\sqrt{\frac{\sum_{i=1}^n\left(Y_i-\bar{r}_n\left(X_i\right)\right)^2}{n^2h_n\hat{f}_n\left(x\right)}\int_{\mathbb{R}}K^2\left(z\right)dz},\right.\\
&&\left.\,\bar{r}_n\left(x\right)+1.96\sqrt{\frac{\sum_{i=1}^n\left(Y_i-\bar{r}_n\left(X_i\right)\right)^2}{n^2h_n\hat{f}_n\left(x\right)}\int_{\mathbb{R}}K^2\left(z\right)dz}\right].
\end{eqnarray*}
The widths of the intervals $\tilde{I}_n$ and $\bar{I}_n$ are of the same order, but 
the asymptotic level of $\bar{I}_n$ is larger than that of $\tilde{I}_n$. More precisely, 
let $\Phi$ denote the distribution function of the standard normal; 
the application of Corollary~\ref{aC:ra} ensures that the asymptotic 
level of  $\bar{I}_n$ is $2\Phi\left(1.96/\sqrt{4/5}\right)-1=97.14\%$.
\end{itemize}

 We consider three sample sizes $n=50$, $n=100$ and $n=200$, 
three regression functions $r(x)=\cos(x)$,  
$r\left(x\right)=0.3\exp(-4\left(x+1\right)^2)
+0.7\exp(-16\left(x-1\right)^2)$, and $r\left(x\right)=1+0.4x$, 
three points $x=-0.5$, $x=0$ and $x=0.5$, two values of $d$, $d=1$ and $d=2$, 
and three densities of $X$, standard normal, normal mixture and student with $6$ degrees of freedom. In each case the number of simulations is $N=5000$. In each table, the first line corresponds to the use of Nadaraya-Watson's estimator $\tilde{r}_n$ and gives the empirical levels $\#\{r\left(x\right)\in \tilde{I}_n
\}/N$; the second line corresponds to the use of the averaged R\'ev\'esz's estimator $\bar{r}_n$ and gives the empirical levels $\#\{r\left(x\right)\in \bar{I}_n\}/N$.  For convenience, we recall the theoretical levels in the last column CL.\\

The simulations results confirm the theoretical ones: the coverage error of the intervals built up with the averaged R\'ev\'esz's estimator is smaller than the coverage error of the intervals built up with Nadaraya-Watson's estimator. 

\paragraph{Model $\mathbf{r\left(x\right)=\cos\left(x\right)}$}
\begin{eqnarray*}
\begin{tabular}[t]{llllllllllllll}
& \multicolumn{9}{c}{Distribution of $X$: $\mathcal{N}\left(0,1\right)$}\\
 & \multicolumn{3}{c}{$x=-0.5$}  & \multicolumn{3}{c}{$x=0$}  & \multicolumn{3}{c}{$x=0.5$} & CL \\ 
&$n=50$ &$n=100$ &$n=200$ &$n=50$ &$n=100$ &$n=200$& $n=50$ &$n=100$  &$n=200$ &\\
$$ & $$ &$$ & $$ &$$ &$d=1$ &$$ &$$  &$$ & $$& \\
& $96.5\%$ & $96.76\%$ & $96.5\%$ &$96.44\%$ & $96.62\%$ & $96.84\%$ &$96.7\%$ & $97.04\%$ & $96.92\%$ & $95\%$ \\ 
& $99.82\%$ & $99.9\%$ & $99.92\%$ & $99.8\%$ &$99.68\%$ & $99.76\%$ & $99.94\%$ & $99.86\%$ & $99.88\%$  &$97.14\%$\\ 
 
$$ & $$ &$$ & $$ &$$ &$d=2$ &$$ &$$  &$$ & $$& \\
& $95.42\%$ & $95.32\%$ & $95.7\%$ &$94.94\%$ & $95.44\%$ & $95.08\%$ &$95.4\%$ & $95.44\%$ & $96.2\%$ & $95\%$ \\ 
& $99.82\%$ & $99.86\%$ & $99.76\%$ &$99.66\%$ & $99.6\%$ & $99.44\%$ &$99.82\%$ & $99.9\%$ & $99.98\%$ & $97.14\%$ \\ 

\end{tabular}
\end{eqnarray*}

\paragraph{Model $\mathbf{r\left(x\right)=0.3\exp\left(-4\left(x+1\right)^2\right)+0.7\exp\left(-16\left(x-1\right)^2\right)}$}
\begin{eqnarray*}
\begin{tabular}[t]{llllllllllllll}
& \multicolumn{9}{c}{Distribution of $X$: $\mathcal{N}\left(0,1\right)$}\\
 & \multicolumn{3}{c}{$x=-0.5$}  & \multicolumn{3}{c}{$x=0$}  & \multicolumn{3}{c}{$x=0.5$} & CL \\ 
&$n=50$ &$n=100$ &$n=200$ &$n=50$ &$n=100$ &$n=200$& $n=50$ &$n=100$  &$n=200$ &\\
$$ & $$ &$$ & $$ &$$ &$d=1$ &$$ &$$  &$$ & $$& \\
& $95.04\%$ & $94.74\%$ & $95.08\%$ &$95.06\%$ & $95.28\%$ & $95.4\%$ &$95.44\%$ & $95.44\%$ & $95.84\%$ & $95\%$ \\ 
& $99.8\%$ & $99.62\%$ & $99.46\%$ & $99.24\%$ &$99.34\%$ & $99.06\%$ & $99.34\%$ & $99.34\%$ & $99.12\%$  &$97.14\%$\\ 

 $$ & $$ &$$ & $$ &$$ &$d=2$ &$$ &$$  &$$ & $$& \\
& $95.26\%$ & $95.14\%$ & $95.34\%$ &$94.74\%$ & $94.88\%$ & $95.06\%$ &$94.48\%$ & $95.56\%$ & $95.62\%$ & $95\%$ \\ 
 & $99.86\%$ & $99.76\%$ & $99.72\%$ &$99.64\%$ & $99.52\%$ & $99.38\%$ & $99.62\%$ & $99.74\%$ & $99.6\%$ &$97.14\%$\\ 

\end{tabular}
\end{eqnarray*}

\paragraph{Model $\mathbf{r\left(x\right)=1+0.4x}$}
\begin{eqnarray*}
\begin{tabular}[t]{llllllllllllll}
& \multicolumn{9}{c}{Distribution of $X$: ${\mathcal{N}\left(0,1\right)}$}\\
 & \multicolumn{3}{c}{$x=-0.5$}  & \multicolumn{3}{c}{$x=0$}  & \multicolumn{3}{c}{$x=0.5$} & CL \\ 
&$n=50$ &$n=100$ &$n=200$ &$n=50$ &$n=100$ &$n=200$& $n=50$ &$n=100$  &$n=200$ &\\
$$ & $$ &$$ & $$ &$$ &$d=1$ &$$ &$$  &$$ & $$& \\
& $96.32\%$ & $95.94\%$ & $96.1\%$ &$96.24\%$ & $96.2\%$ & $96\%$ &$96.1\%$ & $96.24\%$ & $96.62\%$ & $95\%$ \\ 
& $99.84\%$ & $99.9\%$ & $99.6\%$ & $99.92\%$& $99.82\%$ &$99.72\%$ & $99.86\%$ & $99.8\%$ & $99.76\%$  &$97.14\%$\\  
  
$$ & $$ &$$ & $$ &$$ &$d=2$ &$$ &$$  &$$ & $$& \\
& $95.46\%$ & $94.76\%$ & $95.16\%$ &$95.56\%$ & $95.38\%$ & $95.54\%$ &$94.98\%$ & $94.96\%$ & $95.62\%$ & $95\%$ \\ 
& $99.82\%$ & $99.88\%$ & $99.62\%$ & $99.88\%$& $99.78\%$ &$99.68\%$ & $99.88\%$ & $99.82\%$ & $99.68\%$  &$97.14\%$\\  

\end{tabular}
\end{eqnarray*}
\paragraph{Model $\mathbf{r\left(x\right)=\cos\left(x\right)}$}
\begin{eqnarray*}
\begin{tabular}[t]{llllllllllllll}
& \multicolumn{9}{c}{Distribution of $X$: ${1/2\mathcal{N}\left(-1/2,1\right)+1/2\mathcal{N}\left(1/2,1\right)}$}\\
 & \multicolumn{3}{c}{$x=-0.5$}  & \multicolumn{3}{c}{$x=0$}  & \multicolumn{3}{c}{$x=0.5$} & CL \\ 
&$n=50$ &$n=100$ &$n=200$ &$n=50$ &$n=100$ &$n=200$& $n=50$ &$n=100$  &$n=200$ &\\ 
$$ & $$ &$$ & $$ &$$ &$d=1$ &$$ &$$  &$$ & $$& \\
& $96.96\%$ & $97.06\%$ & $97.12\%$ &$97.26\%$ & $96.8\%$ & $97.1\%$ &$97.46\%$ & $96.94\%$ & $96.94\%$ & $95\%$ \\ 
& $99.96\%$ & $99.92\%$ & $99.88\%$ & $99.86\%$ &$99.8\%$ & $99.66\%$ & $99.96\%$ & $99.96\%$ & $99.8\%$  &$97.14\%$\\ 

$$ & $$ &$$ & $$ &$$ &$d=2$ &$$ &$$  &$$ & $$& \\
& $95.6\%$ & $95.32\%$ & $95.56\%$ &$95.08\%$ & $95.36\%$ & $95.64\%$ &$96.38\%$ & $95.7\%$ & $95.34\%$ & $95\%$ \\ 
& $99.82\%$ & $99.92\%$ & $99.74\%$ & $99.94\%$ &$99.78\%$ & $99.64\%$ & $99.96\%$ & $99.9\%$ & $99.64\%$  &$97.14\%$\\ 

\end{tabular}
\end{eqnarray*}

\paragraph{Model $\mathbf{r\left(x\right)=0.3\exp\left(-4\left(x+1\right)^2\right)+0.7\exp\left(-16\left(x-1\right)^2\right)}$}
\begin{eqnarray*}
\begin{tabular}[t]{llllllllllllll}
& \multicolumn{9}{c}{Distribution of $X$: ${1/2\mathcal{N}\left(-1/2,1\right)+1/2\mathcal{N}\left(1/2,1\right)}$}\\
 & \multicolumn{3}{c}{$x=-0.5$}  & \multicolumn{3}{c}{$x=0$}  & \multicolumn{3}{c}{$x=0.5$} & CL \\ 
&$n=50$ &$n=100$ &$n=200$ &$n=50$ &$n=100$ &$n=200$& $n=50$ &$n=100$  &$n=200$ &\\
$$ & $$ &$$ & $$ &$$ &$d=1$ &$$ &$$  &$$ & $$& \\
& $94.9\%$ & $95.38\%$ & $95.3\%$ &$95.56\%$ & $94.56\%$ & $94.86\%$ &$95.24\%$ & $95.24\%$ & $95.48\%$ & $95\%$ \\ 
& $99.74\%$ & $99.62\%$ & $99.58\%$ & $99.44\%$ &$99.22\%$ & $99.1\%$ & $99.34\%$ & $99.28\%$ & $99.06\%$  &$97.14\%$\\  

 $$ & $$ &$$ & $$ &$$ &$d=2$ &$$ &$$  &$$ & $$& \\
& $94.54\%$ & $95.34\%$ & $94.92\%$ &$95.2\%$ & $94.4\%$ & $94.82\%$ &$95.24\%$ & $95.06\%$ & $95.14\%$ & $95\%$ \\ 
& $99.82\%$ & $99.78\%$ & $99.74\%$ & $99.84\%$ &$99.74\%$ & $99.6\%$ & $99.8\%$ & $99.78\%$ & $99.58\%$  &$97.14\%$\\  

\end{tabular}
\end{eqnarray*}

\paragraph{Model $\mathbf{r\left(x\right)=1+0.4x}$}
\begin{eqnarray*}
\begin{tabular}[t]{llllllllllllll}
& \multicolumn{9}{c}{Distribution of $X$: ${1/2\mathcal{N}\left(-1/2,1\right)+1/2\mathcal{N}\left(1/2,1\right)}$}\\
 & \multicolumn{3}{c}{$x=-0.5$}  & \multicolumn{3}{c}{$x=0$}  & \multicolumn{3}{c}{$x=0.5$} & CL \\ 
&$n=50$ &$n=100$ &$n=200$ &$n=50$ &$n=100$ &$n=200$& $n=50$ &$n=100$  &$n=200$ &\\ 
$$ & $$ &$$ & $$ &$$ &$d=1$ &$$ &$$  &$$ & $$& \\
& $96.32\%$ & $96.66\%$ & $96.84\%$ &$96.46\%$ & $96.74\%$ & $96.64\%$ &$96.6\%$ & $96.72\%$ & $97.2\%$ & $95\%$ \\ 
& $99.92\%$ & $99.88\%$ & $99.8\%$ & $99.94\%$& $99.98\%$ &$99.84\%$ & $99.88\%$ & $99.9\%$ & $99.86\%$  &$97.14\%$\\

$$ & $$ &$$ & $$ &$$ &$d=2$ &$$ &$$  &$$ & $$& \\
& $95.18\%$ & $95.46\%$ & $96.1\%$ &$95.08\%$ & $95.52\%$ & $95.6\%$ &$95.58\%$ & $95.44\%$ & $95.74\%$ & $95\%$ \\ 
& $99.94\%$ & $99.86\%$ & $99.78\%$ & $99.88\%$& $99.96\%$ &$99.7\%$ & $99.9\%$ & $99.86\%$ & $99.8\%$  &$97.14\%$\\
 
\end{tabular}
\end{eqnarray*}
\paragraph{Model $\mathbf{r\left(x\right)=\cos\left(x\right)}$}
\begin{eqnarray*}
\begin{tabular}[t]{llllllllllllll}
& \multicolumn{9}{c}{Distribution of $X$: ${\mathcal{T}\left(6\right)}$}\\
 & \multicolumn{3}{c}{$x=-0.5$}  & \multicolumn{3}{c}{$x=0$}  & \multicolumn{3}{c}{$x=0.5$} & CL \\ 
&$n=50$ &$n=100$ &$n=200$ &$n=50$ &$n=100$ &$n=200$& $n=50$ &$n=100$  &$n=200$ &\\
$$ & $$ &$$ & $$ &$$ &$d=1$ &$$ &$$  &$$ & $$& \\
& $96.98\%$ & $97.54\%$ & $97.64\%$ &$97.02\%$ & $97.28\%$ & $97.52\%$ &$97.6\%$ & $97.1\%$ & $96.98\%$ & $95\%$ \\ 
& $99.9\%$ & $99.84\%$ & $99.62\%$ & $99.74\%$ &$99.86\%$ & $99.88\%$ & $99.98\%$ & $99.9\%$ & $99.86\%$  &$97.14\%$\\ 
 
$$ & $$ &$$ & $$ &$$ &$d=2$ &$$ &$$  &$$ & $$& \\
& $95.6\%$ & $95.96\%$ & $95.94\%$ &$95.4\%$ & $95.84\%$ & $96.06\%$ &$96.26\%$ & $95.62\%$ & $95.24\%$ & $95\%$ \\ 
& $99.88\%$ & $99.78\%$ & $99.82\%$ & $99.74\%$ &$99.72\%$ & $99.8\%$ & $99.98\%$ & $99.82\%$ & $99.68\%$  &$97.14\%$\\ 

\end{tabular}
\end{eqnarray*}

\paragraph{Model $\mathbf{r\left(x\right)=0.3\exp\left(-4\left(x+1\right)^2\right)+0.7\exp\left(-16\left(x-1\right)^2\right)}$}
\begin{eqnarray*}
\begin{tabular}[t]{llllllllllllll}
& \multicolumn{9}{c}{Distribution of $X$: ${\mathcal{T}\left(6\right)}$}\\
 & \multicolumn{3}{c}{$x=-0.5$}  & \multicolumn{3}{c}{$x=0$}  & \multicolumn{3}{c}{$x=0.5$} & CL \\ 
&$n=50$ &$n=100$ &$n=200$ &$n=50$ &$n=100$ &$n=200$& $n=50$ &$n=100$  &$n=200$ &\\
$$ & $$ &$$ & $$ &$$ &$d=1$ &$$ &$$  &$$ & $$& \\
& $95.3\%$ & $94.88\%$ & $95.08\%$ &$95.5\%$ & $95.06\%$ & $95.02\%$ &$95.28\%$ & $95.48\%$ & $95.56\%$ & $95\%$ \\ 
& $99.8\%$ & $99.68\%$ & $99.46\%$ & $99.16\%$ &$99.26\%$ & $99.18\%$ & $99.4\%$ & $99.24\%$ & $99.18\%$  &$97.14\%$\\  

 $$ & $$ &$$ & $$ &$$ &$d=2$ &$$ &$$  &$$ & $$& \\
& $94.88\%$ & $94.5\%$ & $94.8\%$ &$95.28\%$ & $94.8\%$ & $94.64\%$ &$95.06\%$ & $95.3\%$ & $95.3\%$ & $95\%$ \\ 
& $99.84\%$ & $99.82\%$ & $99.58\%$ & $99.64\%$ &$99.66\%$ & $99.58\%$ & $99.8\%$ & $99.7\%$ & $99.7\%$  &$97.14\%$\\  

\end{tabular}
\end{eqnarray*}

\paragraph{Model $\mathbf{r\left(x\right)=1+0.4x}$}
\begin{eqnarray*}
\begin{tabular}[t]{llllllllllllll}
& \multicolumn{9}{c}{Distribution of $X$: ${\mathcal{T}\left(6\right)}$}\\
 & \multicolumn{3}{c}{$x=-0.5$}  & \multicolumn{3}{c}{$x=0$}  & \multicolumn{3}{c}{$x=0.5$} & CL \\ 
&$n=50$ &$n=100$ &$n=200$ &$n=50$ &$n=100$ &$n=200$& $n=50$ &$n=100$  &$n=200$ &\\ 
$$ & $$ &$$ & $$ &$$ &$d=1$ &$$ &$$  &$$ & $$& \\
& $96.62\%$ & $97.04\%$ & $97\%$ &$97.2\%$ & $97.08\%$ & $97.02\%$ &$96.36\%$ & $97.14\%$ & $97.22\%$ & $95\%$ \\ 
& $99.84\%$ & $99.9\%$ & $99.92\%$ & $99.94\%$& $99.88\%$ &$99.82\%$ & $99.86\%$ & $99.84\%$ & $99.86\%$  &$97.14\%$\\
 
$$ & $$ &$$ & $$ &$$ &$d=2$ &$$ &$$  &$$ & $$& \\
& $95.04\%$ & $95.62\%$ & $95.54\%$ &$95.96\%$ & $95.58\%$ & $95.88\%$ &$94.94\%$ & $96.14\%$ & $95.86\%$ & $95\%$ \\ 
& $99.82\%$ & $99.9\%$ & $99.82\%$ & $99.84\%$& $99.78\%$ &$99.66\%$ & $99.86\%$ & $99.84\%$ & $99.76\%$  &$97.14\%$\\

\end{tabular}
\end{eqnarray*}

\section{Outlines of the proofs} \label{asection 3}
From now on, we set $n_0\geq 3$ such that $\forall k\geq n_0$, $\gamma_k \leq \left(2\|f\|_{\infty}\right)^{-1}$ and $\gamma_kh_k^{-1}\|K\|_{\infty}\leq 1$. 
Moreover, we introduce the following notations:
\begin{eqnarray*}
s_n&=&\sum_{k=n_0}^n \gamma_k,\\
Z_n\left(x\right)&=&h_n^{-1}K\left(\frac{x-X_n}{h_n}\right),\\
W_n\left(x\right)&=&h_n^{-1}Y_nK\left(\frac{x-X_n}{h_n}\right),\\
\end{eqnarray*}
and, for $s>0$,
\begin{eqnarray*}
\Pi_n\left(s\right)&=&\prod_{j=n_0}^{n}\left(1-s\gamma_j\right),\\
U_{k,n}\left(s\right)&=&\Pi_n\left(s\right)\Pi_k^{-1}\left(s\right).\\
\end{eqnarray*}
Finally, we define the sequences $\left(m_n\right)$ and $\left(\tilde{m}_n\right)$ by 
setting 
\begin{eqnarray}\label{amn12}
\left(m_n\right)&= &\left\{\begin{array}{llll}\left(\sqrt{\gamma_nh_n^{-1}}\right) & \mbox{ if } \lim_{n\to\infty}\left(\gamma_nh_n^{-5}\right)=\infty,\\
\left(h_n^2\right) & \mbox{ otherwise. } 
\end{array}\right.\\
\label{amnt12}
\left(\tilde{m}_n\right)&=&\left\{\begin{array}{llll}\left(\sqrt{\gamma_nh_n^{-1}}\ln n\right) & \mbox{ if }  \lim_{n\to\infty}\left(\gamma_nh_n^{-5}\ln n\right)=\infty,\\
\left(h_n^2\right) & \mbox{ otherwise. } 
\end{array}\right.
\end{eqnarray}
Note that the sequences $\left(m_n\right)$ and $\left(\tilde{m}_n\right)$ belong to $\mathcal{GS}\left(-m^*\right)$ with 
\begin{eqnarray}\label{a10bis}
m^*=\min\left\{\frac{\alpha-a}{2},2a\right\}.
\end{eqnarray}

Before giving the outlines of the proofs, we state the following technical lemma, which is proved in Mokkadem et al. (2008), and which will be used throughout the demonstrations.
\begin{lemma}\label{alap} $ $\\
Let $\left(v_n\right)\in \mathcal{GS}\left(v^*\right)$, $\left(\gamma_n\right)\in\mathcal{GS}\left(-\alpha\right)$ with $\alpha>0$, and set $m>0$. If $ms-v^*\xi>0$ (where $\xi$ is defined in~\eqref{axi1p}), then
\begin{eqnarray*}
\lim_{n \to \infty}v_n\Pi_n^{m}\left(s\right)\sum_{k=n_0}^n\Pi_{k}^{-m}\left(s\right)\frac{\gamma_k}{v_k}=\frac{1}{ms-v^*\xi}. 
\end{eqnarray*}
Moreover, for all positive sequence $\left(\alpha_n\right)$ such that $\lim_{n \to \infty}\alpha_n=0$, and all $C$,
\begin{eqnarray*}
\lim_{n \to \infty}v_n\Pi_n^{m}\left(s\right)\left[\sum_{k=n_0}^n\Pi_{k}^{-m}\left(s\right) \frac{\gamma_k}{v_k}\alpha_k+C\right]=0.
\end{eqnarray*}
\end{lemma}
As explained in the introduction, we note that the stochastic approximation algorithm~\eqref{arun} can be rewritten as:
\begin{eqnarray}
r_n(x) & =&\left(1-\gamma_nZ_n\left(x\right)\right)r_{n-1}(x)+\gamma_nW_n\left(x\right)\label{arnMar}\\
&=&\left(1-\gamma_nf\left(x\right)\right)r_{n-1}(x)+\gamma_n\left(f\left(x\right)-Z_n\left(x\right)\right)r_{n-1}\left(x\right)+\gamma_nW_n\left(x\right)\label{arn1}.
\end{eqnarray}
To establish the asymptotic behaviour of $\left(r_n\right)$ and $\left(\bar{r}_n\right)$, we introduce the auxiliary stochastic approximation algorithm defined by setting $\rho_n\left(x\right)=r\left(x\right)$ for all $n\leq n_0-2$, $\rho_{n_0-1}\left(x\right)=r_{n_0-1}\left(x\right)$, and, for $n\geq n_0$,
\begin{eqnarray}\label{arnhat}
\rho_n(x)&=&\left(1-\gamma_nf\left(x\right)\right)\rho_{n-1}(x)+\gamma_n\left(f\left(x\right)-Z_n\left(x\right)\right)r\left(x\right)+\gamma_nW_n\left(x\right).
\end{eqnarray}
We first give the asymptotic behaviour of $\left(\rho_n\right)$ and of $\left(\bar{\rho}_n\right)$ in Section \ref{arh} and \ref{arhbar} respectively (and refer to Section \ref{alem} for the proof of the different lemmas). Then, we show in Section \ref{asection 3.3} how the asymptotic behaviour of $\left(r_n\right)$ and $\left(\bar{r}_n\right)$ can be deduced from that of $\left(\rho_n\right)$ and $\left(\bar{\rho}_n\right)$ respectively.

\subsection{Asymptotic behaviour of $\rho_n$}\label{arh}
The aim of this section is to give the outlines of the proof of the three following lemmas.

\begin{lemma}[Weak pointwise convergence rate of $\rho_n$]\label{aL:rhatlc} $ $\\
Let Assumptions $\left( A1\right)-\left( A3\right) $ hold for $x \in \mathbb{R}$ such that $f\left(x\right)\not=0$.
\begin{enumerate}
\item If there exists $c \geq 0 $ such that $\gamma_n^{-1}h_n^5\rightarrow c$, and if $\lim_{n\to \infty}\left(n\gamma_n\right)>\left(1-a\right)/\left(2f\left(x\right)\right)$, then
\begin{eqnarray}\label{atlc:rb}
\lefteqn{\sqrt{\gamma_n^{-1}h_n}\left( \rho_{n}\left( x\right)-r\left( x\right) \right)}\nonumber\\
&&\stackrel{\mathcal{D}}{\rightarrow}\mathcal{N}\left(\frac{\sqrt{c}f\left(x\right)m^{\left(2\right)}\left(x\right)}{f\left(x\right)-2a\xi},\frac{Var\left[Y\vert X=x\right]f\left(x\right)}{\left(2f\left(x\right)-\left(\alpha-a\right)\xi \right)}\int_{\mathbb{R}} K^2\left(z\right)dz\right).
\end{eqnarray}
\item If  $\gamma_n^{-1}h_n^5\rightarrow \infty $, and if $\lim_{n\to \infty}\left(n\gamma_n\right)>2a/f\left(x\right)$, then  
\begin{eqnarray}\label{atlc:rc}
\frac{1}{h_{n}^{2}}\left( \rho_{n}\left( x\right)-r\left( x\right) \right) \stackrel{\mathbb{P}}{\rightarrow }\frac{f\left(x\right)m^{\left(2\right)}\left(x\right)}{\left(f\left(x\right)-2a \xi\right)} .
\end{eqnarray}
\end{enumerate}
\end{lemma}

\begin{lemma}[Strong pointwise convergence rate of $\rho_{n}$]\label{aL:rhatlog} $ $\\
Let Assumptions $\left( A1\right)-\left( A3\right) $ hold for $x \in \mathbb{R}$ such that $f\left(x\right)\not=0$.
\begin{enumerate}
\item If there exists $c\geq 0 $ such that $\gamma_n^{-1}h_n^5/\ln\left(s_n\right) \rightarrow c$, and if $\lim_{n\to \infty}\left(n\gamma_n\right)>\left(1-a\right)/\left(2f\left(x\right)\right)$, then, with probability one, the sequence
\[
\left(\sqrt{\frac{\gamma_n^{-1}h_n}{2 \ln\left(s_n\right)}}\left( \rho_{n}\left( x\right)-r\left( x\right) \right)\right)  
\]
is relatively compact and its limit set is the interval
\begin{eqnarray*}
\left[\sqrt{\frac{c}{2}}\frac{f\left(x\right)m^{\left(2\right)}\left( x\right)}{f\left(x\right)-2a\xi}-\sqrt{\frac{Var\left[Y\vert X=x\right]f\left(x\right)\int_{\mathbb{R}} K^2\left(z\right)dz}{\left(2f\left(x\right)-\left(\alpha-a\right)\xi\right)}},\right.\\
\left.\sqrt{\frac{c}{2}}\frac{f\left(x\right)m^{\left(2\right)}\left( x\right)}{f\left(x\right)-2a\xi}+\sqrt{\frac{Var\left[Y\vert X=x\right]f\left(x\right)\int_{\mathbb{R}} K^2\left(z\right)dz}{\left(2f\left(x\right)-\left(\alpha-a\right)\xi\right)}}\right].
\end{eqnarray*}
\item If $\gamma_n^{-1}h_n^5/\ln\left(s_n\right)\rightarrow \infty $, and if $\lim_{n\to \infty}\left(n\gamma_n\right)>2a/f\left(x\right)$ then, with probability one,
\begin{eqnarray*}
\lim_{n\rightarrow \infty}\frac{1}{h_{n}^{2}}\left( \rho_{n}\left( x\right)-r\left( x\right) \right) = \frac{f\left(x\right)m^{\left(2\right)}\left( x\right)}{f\left(x\right)-2a\xi}. 
\end{eqnarray*}
\end{enumerate}
\end{lemma}

\begin{lemma}[Strong uniform convergence rate of $\rho_n$]\label{aL:Tuhat} $ $\\
Let  $I$ be a bounded open interval on which $\varphi=\inf_{x \in I}f\left(x\right)>0$,  
and let Assumptions $\left( A1\right)-\left( A4\right) $ hold for all $x\in I$. 
If $\lim_{n\tend\infty}n\g_n>\min\left\{(1-a)/(2\varphi),2a/\varphi\right\}$, then 
\begin{eqnarray*}
\sup_{x \in I}\left|\rho_n\left( x\right)-r\left(x\right)\right|=O\left(\max\left\{\sqrt{\gamma_nh_n^{-1}} \ln n ,h_n^2\right\}\right) \quad a.s.
\end{eqnarray*}
\end{lemma}

To prove Lemmas \ref{aL:rhatlc} and \ref{aL:rhatlog}, we first remark that, in 
view of \eqref{arnhat}, we have, for $n\geq n_0$,  
\begin{eqnarray}\label{arnhat1}
\rho_n(x)-r\left(x\right)&=&\left(1-\gamma_nf\left(x\right)\right)\left(\rho_{n-1}(x)-r\left(x\right)\right)+\gamma_n\left(W_n\left(x\right)-r\left(x\right)Z_n\left(x\right)\right)\nonumber\\
&=&\Pi_n\left(f\left(x\right)\right)\sum_{k=n_0}^n\Pi_k^{-1}\left(f\left(x\right)\right)\gamma_k\left(W_k\left(x\right)-r\left(x\right)Z_k\left(x\right)\right)\nonumber\\
&&+\Pi_n\left(f\left(x\right)\right)\left(\rho_{n_0-1}\left(x\right)-r\left(x\right)\right)\nonumber\\
&=&\tilde{T}_n\left(x\right)+\tilde{R}_n\left(x\right),
\end{eqnarray}
with, since $\rho_{n_0-1}=r_{n_0-1}$, 
\begin{eqnarray*}
\tilde{T}_n\left(x\right)&=&\sum_{k=n_0}^nU_{k,n}\left(f\left(x\right)\right)\gamma_k\left(W_k\left(x\right)-r\left(x\right)Z_k\left(x\right)\right),\nonumber\\
\tilde{R}_n\left(x\right)&=&\Pi_n\left(f\left(x\right)\right)\left(r_{n_0-1}\left(x\right)-r\left(x\right)\right).
\end{eqnarray*}
Noting that $\left|r_{n_0-1}\left(x\right)-r\left(x\right)\right|=O\left(1\right)$ a.s. and applying Lemma~\ref{alap}, we get
\begin{eqnarray*}
\left|\tilde{R}_n\left(x\right)\right|&=&O\left(\Pi_n\left(f\left(x\right)\right)\right)\quad a.s.\\
&=&o\left(m_n\right)\quad a.s.
\end{eqnarray*}
Lemmas~\ref{aL:rhatlc} and~\ref{aL:rhatlog} are thus straightforward consequences of the following lemmas, which are proved in Sections~\ref{a4.0} and~\ref{aprlb}, respectively.
\begin{lemma}\label{a2bis}
The two parts of Lemma \ref{aL:rhatlc} hold when $\rho_n\left(x\right)-r\left(x\right)$ is replaced by $\tilde{T}_n\left(x\right)$.
\end{lemma}
\begin{lemma}\label{a3bis}
The two parts of Lemma \ref{aL:rhatlog} hold when $\rho_n\left(x\right)-r\left(x\right)$ is replaced by $\tilde{T}_n\left(x\right)$.
\end{lemma}
In the same way, we remark that
\begin{eqnarray*}
\sup_{x\in I}\left|\tilde{R}_n\left(x\right)\right|&=&O\left(\sup_{x\in I}\Pi_n\left(f\left(x\right)\right)\right)\quad a.s.\\
&=&O\left(\Pi_n\left(\varphi\right)\right)\quad a.s.\\
&=&o\left(m_n\right)\quad a.s.,
\end{eqnarray*}
so that Lemma \ref{aL:Tuhat} is a straightforward consequence of the following lemma, which is proved in Section \ref{a4.1}.
\begin{lemma}\label{a4bis}
Lemma \ref{aL:Tuhat} holds when $\rho_n-r$ is replaced by $\tilde{T}_n$.
\end{lemma}
\subsection{Asymptotic behaviour of $\bar{\rho}_n$}\label{arhbar}
The purpose of this section is to give the outlines of the proof of the three following lemmas.

\begin{lemma}[Weak pointwise convergence rate of $\bar{\rho}_n$]\label{aL:arb} $ $\\
Let Assumptions $\left( A1\right)-\left( A3\right)$, $\left(A5\right)$ and $\left(A6\right)$ hold for $x \in \mathbb{R}$ such that $f\left(x\right)\not=0$.
\begin{enumerate}
\item If there exists $c \geq 0 $ such that  $nh_n^5\rightarrow  c$, then
\begin{eqnarray*}
\lefteqn{\sqrt{nh_n}\left(\bar{\rho}_n\left(x\right)-r\left(x\right)\right)}\\
&&\stackrel{\mathcal{D}}{\rightarrow}\mathcal{N}\left(c^{\frac{1}{2}}\frac{1-q}{1-q-2a}m^{\left(2\right)}\left(x\right),\frac{\left(1-q\right)^2}{1+a-2q}\frac{Var\left[Y\vert X=x\right]}{f\left(x\right)}\int_{\mathbb{R}} K^2\left(z\right)dz\right).
\end{eqnarray*}
\item If  $nh_n^5\rightarrow \infty $, then  
\begin{eqnarray*}
h_{n}^{-2}\left(\bar{\rho}_n\left(x\right)-r\left(x\right)\right) \stackrel{\mathbb{P}}{\rightarrow }\frac{1-q}{1-q-2a}m^{\left(2\right)}\left(x\right).
\end{eqnarray*}
\end{enumerate}
\end{lemma}

\begin{lemma}[Strong pointwise convergence rate of $\bar{\rho}_n$]\label{aL:brb} $ $\\
Let Assumptions $\left( A1\right)-\left( A3\right)$, $\left(A5\right)$ and $\left(A6\right)$ hold for $x\in \mathbb{R}$ such that $f\left(x\right)\not=0$.
\begin{enumerate}
\item If there exists $c_{1} \geq 0 $ such that $nh_n^5/\ln \ln n\rightarrow c_{1}$, then, with probability one, the sequence
\[
\left(\sqrt{\frac{nh_n}{2\ln \ln n}}\left(\bar{\rho}_n\left(x\right)-r\left(x\right)\right)\right)
\]
is relatively compact and its limit set is the interval
\begin{eqnarray*}
\lefteqn{\left[\frac{1-q}{1-q-2a}\sqrt{\frac{c_1}{2}}m^{\left(2\right)}\left( x\right)-\sqrt{\frac{\left(1-q\right)^2}{1+a-2q}\frac{Var\left[Y/X=x\right]}{f\left(x\right)}\int_{\mathbb{R}} K^2\left(z\right)dz},\right.}\\
&&\left.\frac{1-q}{1-q-2a}\sqrt{\frac{c_1}{2}}m^{\left(2\right)}\left( x\right)+\sqrt{\frac{\left(1-q\right)^2}{1+a-2q}\frac{Var\left[Y/X=x\right]}{f\left(x\right)}\int_{\mathbb{R}} K^2\left(z\right)dz}\right].
\end{eqnarray*}
\item If $nh_n^5/\ln \ln n \rightarrow \infty $, then  
\begin{eqnarray*}
\lim_{n\rightarrow \infty}h_{n}^{-2}\left(\bar{\rho}_n\left(x\right)-r\left(x\right)\right) = \frac{1-q}{1-q-2a}m^{\left(2\right)}\left( x\right)\quad a.s. 
\end{eqnarray*}
\end{enumerate}
\end{lemma}

\begin{lemma}[Strong uniform convergence rate of $\bar{\rho}_n$]\label{aL:Tubar} $ $\\
Let $I$ be a bounded open interval on which $\varphi=\inf_{x \in I}f\left(x\right)>0$ and let Assumptions $\left( A1\right)-\left( A6\right) $ hold for all $x\in I$. We have 
\begin{eqnarray*}
\sup_{x \in I}\left|\bar{\rho}_n\left(x\right)-r\left(x\right)\right|=O\left(\max\left\{\sqrt{n^{-1}h_n^{-1}} \ln n ,h_n^2\right\}\right) \quad a.s.
\end{eqnarray*}
\end{lemma} 

To prove Lemmas \ref{aL:arb}-\ref{aL:Tubar}, we note that \eqref{arnhat} gives, 
for $n\geq n_0$, 
\begin{eqnarray*}
\rho_{n}\left(x\right)-\rho_{n-1}\left(x\right)=-\gamma_nf\left(x\right)\left[\rho_{n-1}(x)-r\left(x\right)\right]+\gamma_n\left[W_n\left(x\right)-r\left(x\right)Z_n\left(x\right)\right],
\end{eqnarray*}
and thus
\begin{eqnarray*}
\rho_{n-1}(x)-r\left(x\right)&=&\frac{1}{f\left(x\right)}\left[W_n\left(x\right)-r\left(x\right)Z_n\left(x\right)\right]-\frac{1}{\gamma_nf\left(x\right)}\left[\rho_{n}(x)-\rho_{n-1}(x)\right].
\end{eqnarray*}
It follows that
\begin{eqnarray}\label{a24jul}
\bar{\rho}_n\left(x\right)-r\left(x\right)&=&\frac{1}{\sum_{k=1}^nq_k}\sum_{k=1}^nq_k\left[\rho_k\left(x\right)-r\left(x\right)\right]\nonumber\\
&=&\frac{1}{f\left(x\right)}T_n\left(x\right)-\frac{1}{f\left(x\right)}R_n^{\left(0\right)}\left(x\right)
\end{eqnarray}
with 
\begin{eqnarray*}
T_n\left(x\right)&=&\frac{1}{\sum_{k=1}^nq_k}\sum_{k=n_0-1}^nq_k\left[W_{k+1}\left(x\right)-r\left(x\right)Z_{k+1}\left(x\right)\right],\\
R_n^{\left(0\right)}\left(x\right)&=&\frac{1}{\sum_{k=n_0-1}^nq_k}\sum_{k=1}^n\frac{q_k}{\gamma_{k+1}}\left[\rho_{k+1}(x)-\rho_{k}(x)\right].
\end{eqnarray*}
Let us note that $R_n^{\left(0\right)}$ can be rewritten as
\begin{eqnarray*}
R_n^{\left(0\right)}\left(x\right)&=&\frac{1}{\sum_{k=1}^nq_k}\sum_{k=n_0-1}^{n}\frac{q_k}{\gamma_{k+1}}\left[\left(\rho_{k+1}\left(x\right)-r\left(x\right)\right)-\left(\rho_{k}\left(x\right)-r\left(x\right)\right)\right]\nonumber\\
&=&\frac{1}{\sum_{k=1}^nq_k}\sum_{k=n_0}^{n}\left(\frac{q_{k-1}}{\gamma_k}-\frac{q_k}{\gamma_{k+1}}\right)\left(\rho_k\left(x\right)-r\left(x\right)\right)\nonumber\\
&&+\frac{1}{\sum_{k=1}^nq_k}\frac{q_n}{\gamma_{n+1}}\left(\rho_{n+1}\left(x\right)-r\left(x\right)\right)-\frac{1}{\sum_{k=1}^nq_k}\frac{q_{n_0-1}}{\gamma_{n_0}}
\left(\rho_{n_0-1}\left(x\right)-r\left(x\right)\right)\nonumber\\
&=&\frac{1}{\sum_{k=1}^nq_k}\sum_{k=n_0}^{n}\frac{q_{k-1}}{\gamma_k}\left[1-\frac{q_{k-1}^{-1}\gamma_k}{q_{k}^{-1}\gamma_{k+1}}\right]\left(\rho_k\left(x\right)-r\left(x\right)\right)\nonumber\\
&&+\frac{1}{\sum_{k=1}^nq_k}\frac{q_n}{\gamma_{n+1}}
\left(\rho_{n+1}\left(x\right)-r\left(x\right)\right)
-\frac{1}{\sum_{k=1}^nq_k}\frac{q_{n_0-1}}{\gamma_{n_0}}
\left(\rho_{n_0-1}\left(x\right)-r\left(x\right)\right).
\end{eqnarray*}
Since $\left(q_{k-1}^{-1}\gamma_k\right) \in \mathcal{GS}\left(q-\alpha\right)$, we have
\begin{eqnarray*}
\left[1-\frac{q_{k-1}^{-1}\gamma_k}{q_{k}^{-1}\gamma_{k+1}}\right]&=&1-\left(1-\frac{\left(q-\alpha \right)}{k}+o\left(\frac{1}{k}\right)\right)\nonumber\\
&=&O\left(k^{-1}\right),
\end{eqnarray*}
and thus
\begin{eqnarray}\label{aR0rbeq}
\left|R_n^{\left(0\right)}\left(x\right)\right|
&=&O\left(\frac{1}{\sum_{k=1}^nq_k}\left[
\sum_{k=n_0}^{n}k^{-1}q_{k-1}\gamma_k^{-1}\left|\rho_k\left(x\right)-r\left(x\right)\right|
+\frac{q_n}{\gamma_{n+1}}\left|\rho_{n+1}\left(x\right)-r\left(x\right)\right|
\right.\right. \nonumber \\ 
& & \mbox{ } \left.\left.
+\frac{q_{n_0-1}}{\gamma_{n_0}}
\left|\rho_{n_0-1}\left(x\right)-r\left(x\right)\right|\right]\right).
\end{eqnarray}
The application of Lemma~\ref{aL:rhatlog} ensures that
\begin{eqnarray*}
\left|R_n^{\left(0\right)}\left(x\right)\right| 
&=&O\left(\frac{1}{\sum_{k=1}^nq_k}\left[\sum_{k=2}^{n}\left(k^{-1}q_k\gamma_{k}^{-1}\right)\left(\left(\gamma_kh_k^{-1}\ln \left(s_k\right)\right)^{\frac{1}{2}}+h_k^2\right)
\right.\right. \nonumber \\ 
& & \mbox{ } \left.\left.
+\frac{q_n}{\gamma_{n+1}}\left(\left(\gamma_nh_n^{-1}\ln \left(s_n\right)\right)^{\frac{1}{2}}+h_n^2\right)+1\right]\right)\quad a.s.
\end{eqnarray*}
Now, let us recall that, if $\left(u_n\right)\in \mathcal{GS}\left(-u^*\right)$ 
with $u^*<1$, then we have, for any fixed $k_0\geq 1$,
\begin{eqnarray}
\lim_{n \to \infty}\frac{nu_n}{\sum_{k=k_0}^nu_k}&=&1-u^*,\label{abt}
\end{eqnarray}
and, if $u^* \geq 1$, then for all $\epsilon>0$, $u_n=O\left(n^{-1+\epsilon}\right)$ and thus
\begin{eqnarray}
\sum_{k=1}^nu_k&=&O\left(n^{\epsilon}\right).\label{abteps}
\end{eqnarray}
Now, set $\epsilon \in \left]0,\min\left\{1-q-2a,\left(1+a\right)/2-q\right\}\right[$ 
(the existence of such an $\epsilon$ being ensured by $\left(A6\right)$); 
in view of $\left(A5\right)$, we get
\begin{eqnarray*}
\left|R_n^{\left(0\right)}\left(x\right)\right|&=&O\left(\frac{1}{nq_n}\left(n^{\epsilon}+q_n\gamma_{n}^{-\frac{1}{2}}h_n^{-\frac{1}{2}}\left(\ln \left(s_n\right)\right)^{\frac{1}{2}}+q_n\gamma_{n}^{-1}h_n^2\right)+\frac{1}{n\gamma_n}\left(\left(\gamma_nh_n^{-1}\ln \left(s_n\right)\right)^{\frac{1}{2}}+h_n^2\right)\right)a.s.\nonumber\\
&=&O\left(\frac{n^{\epsilon}}{nq_n}+\frac{\sqrt{n^{-1}h_n^{-1}}}{\sqrt{n\gamma_n\left(\ln \left(s_n\right)\right)^{-1}}}+\frac{h_n^2}{n\gamma_n}\right)a.s.\nonumber\\
&=&o\left(\sqrt{n^{-1}h_n^{-1}}+h_n^2\right)a.s.
\end{eqnarray*}
In view of~\eqref{a24jul}, Lemmas~\ref{aL:arb} and~\ref{aL:brb} are thus straightforward consequences of the two following lemmas, which are proved in Sections~\ref{a4.2} and~\ref{a4.3} respectively.
\begin{lemma}[Weak pointwise convergence rate of $T_n$]\label{aL:rrb} $ $\\
The two parts of Lemma \ref{aL:arb} hold when $\bar{\rho}_n\left(x\right)-r\left(x\right)$ is replaced by $\left[f\left(x\right)\right]^{-1}T_n\left(x\right)$ .
\end{lemma}

\begin{lemma}[Strong pointwise convergence rate of $T_n$]\label{aL:rrc} $ $\\
The two parts of Lemma \ref{aL:brb} hold when $\bar{\rho}_n\left(x\right)-r\left(x\right)$ is replaced by $\left[f\left(x\right)\right]^{-1}T_n\left(x\right)$ .
\end{lemma}

Now, in view of~\eqref{aR0rbeq}, the application of Lemma~\ref{aL:Tuhat} ensures that 
\begin{eqnarray*}
\lefteqn{\sup_{x\in I}\left|R_n^{\left(0\right)}\left(x\right)\right|}\\
&= &O\left(\frac{1}{\sum_{k=1}^nq_k}\left[\sum_{k=n_0}^{n}k^{-1}q_k\gamma_k^{-1}\sup_{x\in I}\left|\rho_k\left(x\right)-r\left(x\right)\right|+\frac{q_n}{\gamma_{n+1}}\sup_{x\in I}\left|\rho_{n+1}\left(x\right)-r\left(x\right)\right|\right]\right)\quad a.s.\\
&=&O\left(\frac{1}{\sum_{k=1}^nq_k}\left[\sum_{k=n_0}^{n}\left(k^{-1}q_k\gamma_{k}^{-1}\right)\left(\left(\gamma_kh_k^{-1}\right)^{\frac{1}{2}}\ln k+h_k^2\right)+\frac{q_n}{\gamma_{n+1}}\left(\left(\gamma_nh_n^{-1}\right)^{\frac{1}{2}}\ln n+h_n^2\right)\right]\right)\quad a.s.
\end{eqnarray*}
Setting $\epsilon \in \left]0,\min\left\{1-q-2a,\left(1+a\right)/2-q\right\}\right[$ again, 
we get, in view of $\left(A5\right)$, 
\begin{eqnarray*}
\sup_{x\in I}\left|R_n^{\left(0\right)}\left(x\right)\right|&=&O\left(\frac{1}{nq_n}\left(n^{\epsilon}+q_n\gamma_{n}^{-\frac{1}{2}}h_n^{-\frac{1}{2}}\ln n+q_n\gamma_{n}^{-1}h_n^2\right)+\frac{1}{n\gamma_n}\left(\left(\gamma_nh_n^{-1}\right)^{\frac{1}{2}}\ln n+h_n^2\right)\right)\quad a.s.\nonumber\\
&=&O\left(\frac{n^{\epsilon}}{nq_n}+\frac{\sqrt{n^{-1}h_n^{-1}}\,\ln n}{\sqrt{n\gamma_n}}+\frac{h_n^2}{n\gamma_n}\right)a.s.\nonumber\\
&=&o\left(\sqrt{n^{-1}h_n^{-1}}\, \ln n+h_n^2\right)\quad a.s.
\end{eqnarray*}

In view of~\eqref{a24jul}, Lemma~\ref{aL:Tubar} is thus a straightforward consequence of the following lemma, which is proved in Section~\ref{a4.4}.
\begin{lemma}[Strong uniform convergence rate of $T_n$]\label{aL:Turbar} $ $\\
Lemma~\ref{aL:Tubar} holds when $\bar{\rho}_n-r$ is replaced by $T_n$ .
\end{lemma} 
\subsection{How to deduce the asymptotic behaviour of $r_n$ and $\bar{r}_n$ from that of $\rho_n$ and $\bar{\rho}_n$} \label{asection 3.3}
Set
\begin{eqnarray*}
\Delta_n\left(x\right)&=&r_n\left(x\right)-\rho_n\left(x\right)
\end{eqnarray*}
and
\begin{eqnarray*}
\bar{\Delta}_n\left(x\right)&=&\frac{1}{\sum_{k=1}^nq_k}\sum_{k=1}^nq_k\Delta_k\left(x\right)\\
&=&\bar{r}_n\left(x\right)-\bar{\rho}_n\left(x\right).
\end{eqnarray*}
To deduce the asymptotic behaviour of $r_n$ (respectively $\bar{r}_n$) from that of $\rho_n$ (respectively $\bar{\rho}_n$), we prove that $\Delta_n$ (respectively $\bar{\Delta}_n$) is negligible in front of $\rho_n$ (respectively $\bar{\rho}_n$). Note that, in view of~\eqref{arn1} and~\eqref{arnhat}, and since $\rho_{n_0-1}\left(x\right)=r_{n_0-1}\left(x\right)$, we have, for $n\geq n_0$, 
\begin{eqnarray}\label{a26mar}
\Delta_n\left(x\right)&=&\left(1-\gamma_nf\left(x\right)\right)\Delta_{n-1}\left(x\right)+\gamma_n\left(f\left(x\right)-Z_n\left(x\right)\right)\left(r_{n-1}\left(x\right)-r\left(x\right)\right)\nonumber\\
&=&\sum_{k=n_0}^nU_{k,n}\left(f\left(x\right)\right)\left(f\left(x\right)-Z_k\left(x\right)\right)\left(r_{k-1}\left(x\right)-r\left(x\right)\right).
\end{eqnarray}
The difficulty which appears here is that $\Delta_n$ is expressed in function of the terms $r_k-r$, so that an upper bound of $r_n-r$ is necessary for the obtention of an upper bound of $\Delta_n$. Now, the key to overcome this difficulty is the following property $\left(\mathcal{P}\right)$ : if $\left(r_n-r\right)$ is known to be bounded almost surely by a sequence $\left(w_n\right)$, then it can be shown that $\left(\Delta_n\right)$ is bounded almost surely by a sequence $\left(w'_n\right)$ such that $\lim_{n\to\infty}w'_nw_n^{-1}=0$, which may allow to upper bound $r_n-r$ by a sequence smaller than $\left(w_n\right)$.
To deduce the asymptotic behaviour of $r_n$ (respectively $\bar{r}_n$) from that of $\rho_n$ (respectively $\bar{\rho}_n$), we thus proceed as follows. We first establish a rudimentary upper bound of $\left(r_n-r\right)$. Then, applying Property $\left(\mathcal{P}\right)$ several times, we successively improve our upper bound of $\left(r_n-r\right)$, and this until we obtain an upper bound, which allows to prove that $\Delta_n$ (respectively $\bar{\Delta}_n$) is negligible in front of $\rho_n$ (respectively $\bar{\rho}_n$).

We first establish the pointwise results on $r_n$ and $\bar{r}_n$ (that is, 
Theorems \ref{aP:ra},~\ref{aT:rb},~\ref{aAv:ra}, and~\ref{aT:rba}) in Section \ref{asection 3.3.1}, 
and then the uniform ones (that is, Theorems \ref{aT:c} and~\ref{aT:Ac}) in 
Section \ref{asection 3.3.2}.

\subsubsection{Proof of Theorems~\ref{aP:ra},~\ref{aT:rb},~\ref{aAv:ra} and~\ref{aT:rba}}
\label{asection 3.3.1}
The proof of Theorems~\ref{aP:ra},~\ref{aT:rb},~\ref{aAv:ra} and~\ref{aT:rba} relies on the repeted application of the following lemma, which is proved in Section~\ref{a4.5}.
\begin{lemma}\label{aL:rwn} 
Let Assumptions $\left(A1\right)-\left(A3\right)$ hold, and assume that there exists $\left(w_n\right)\in\mathcal{GS}\left(w^*\right)$ such that $\left|r_n\left(x\right)-r\left(x\right)\right|=O\left(w_n\right)$ a.s. 
\begin{enumerate}
\item If the sequence $\left(n\gamma_n\right)$ is bounded, if 
$\lim_{n\tend\infty}n\g_n>\min\left\{(1-a)/(2f(x)),2a/f(x)\right\}$, 
and if $w^*\geq 0$, then, for all $\delta>0$, 
\begin{eqnarray*}
\left|\Delta_n\left(x\right)\right|=O\left(m_nw_n\left(\ln n\right)^{\frac{\left(1+\delta \right)}{2}}\right)+o\left(m_n\right)\quad a.s.
\end{eqnarray*}
\item If $\lim_{n\to \infty}\left(n\gamma_n\right)=\infty$, then, for all $\delta>0$, 
\begin{eqnarray*}
\left|\Delta_n\left(x\right)\right|=O\left(m_nw_n\left(n^{1+\delta}\gamma_n\right)^{\frac{\left(1+\delta \right)}{2}}\right)\quad a.s.
\end{eqnarray*}
\end{enumerate}
\end{lemma}
We first establish a preliminary upper bound for $r_n\left(x\right)-r\left(x\right)$. Then, we successively prove Theorems~\ref{aP:ra} and~\ref{aT:rb} in the case $\left(n\gamma_n\right)$ is bounded, Theorems~\ref{aP:ra} and~\ref{aT:rb} in the case $\lim_{n\to\infty}\left(n\gamma_n\right)=\infty$, and finally Theorems~\ref{aAv:ra} and~\ref{aT:rba}.
\paragraph{Preliminary upper bound of $r_n\left(x\right)-r\left(x\right)$} $ $\\
Since $0\leq 1-\gamma_nZ_n\left(x\right)\leq 1$ for all $n\geq n_0$, it follows from~\eqref{arnMar} that, for $n\geq n_0$,
\begin{eqnarray}
\left|r_n\left(x\right)\right|&\leq &\left|r_{n-1}\left(x\right)\right|+\gamma_n\left|Y_n\right|h_n^{-1}\|K\|_{\infty}\nonumber \\
& \leq &\left|r_{n_0-1}\left(x\right)\right|+\left(\sup_{k\leq n}\left|Y_k\right|\right)\|K\|_{\infty}\sum_{k=1}^n\gamma_kh_k^{-1}\label{abis22}.
\end{eqnarray}
Since
\begin{eqnarray*}
\mathbb{P}\left(\sup_{k\leq n}\left|Y_k\right|>n^2\right)\leq n\mathbb{P}\left(\left|Y\right|>n^2\right)\leq n^{-3}\mathbb{E}\left(\left|Y\right|^2\right),
\end{eqnarray*}
we have $\sup_{k\leq n}\left|Y_k\right| \leq n^2$ a.s. Moreover, since $\left(\gamma_n h_n^{-1}\right) \in \mathcal{GS}\left(-\alpha+a\right)$ with $1-\alpha+a>0$, we note that $\sum_{k=1}^n\gamma_kh_k^{-1}=O\left(n\gamma_nh_n^{-1}\right)$. We thus deduce that
\begin{eqnarray}\label{arn3gh}
\left|r_n\left(x\right)-r\left(x\right)\right|=O\left(n^3\gamma_nh_n^{-1}\right)\quad a.s.
\end{eqnarray}
\paragraph{Proof of Theorems~\ref{aP:ra} and~\ref{aT:rb} in the case the sequence $\left(n\gamma_n\right)$ is bounded.} $ $\\
In this case, $\alpha=1$, and Lemmas~\ref{aL:rhatlog} and~\ref{aL:rwn} imply that:
\begin{eqnarray}\label{alm2r}
&\bullet& \left|\rho_n\left(x\right)-r\left(x\right)\right|=O\left(m_n\ln n\right)\quad a.s.\\
&\bullet& \mbox{If there exists }  \left(w_n\right)\in \mathcal{GS}\left(w^*\right), w^*\geq 0, \mbox{such that } \left|r_n\left(x\right)-r\left(x\right)\right|=O\left(w_n\right) a.s.,\nonumber\\
& &\mbox{then }\left|\Delta_n\left(x\right)\right|
=O\left(m_nw_n\ln n\right)+o\left(m_n\right)\quad a.s.\label{alm5d}
\end{eqnarray}

Set $p_0=\max\left\{p \mbox{ such that}-m^*p+2+a\geq 0\right\}$, set $j \in \left\{0,1,\ldots,p_0-1\right\}$, and assume that
\begin{eqnarray}\label{a12ju}
\left|r_n\left(x\right)-r\left(x\right)\right|&=&O\left(m_n^j\left(n^3\gamma_nh_n^{-1}\right)\left(\ln n\right)^j\right)\quad a.s.
\end{eqnarray}
Since the sequence $\left(w_n\right)=\left(m_n^j\left(n^3\gamma_nh_n^{-1}\right)\left(\ln n\right)^j\right)$ belongs to $\mathcal{GS}\left(-m^*j+2+a\right)$ with $-m^*j+2+a>0$, the application of~\eqref{alm5d} implies that
\begin{eqnarray*}
\left|\Delta_n\left(x\right)\right|=O\left(m_n^{j+1}\left(n^3\gamma_nh_n^{-1}\right)\left(\ln n\right)^{j+1}\right)+o\left(m_n\right)\quad a.s.
\end{eqnarray*}
Since $\left(m_n^{j+1}\left(n^3\gamma_nh_n^{-1}\right)\left(\ln n\right)^{j+1}\right) \in \mathcal{GS}\left(-m^*\left(j+1\right)+2+a\right)$ with $-m^*\left(j+1\right)+2+a\geq 0$, whereas $\left(m_n\right) \in \mathcal{GS}\left(-m^*\right)$ with $-m^*<0$, it follows that
\begin{eqnarray*}
\left|\Delta_n\left(x\right)\right|=O\left(m_n^{j+1}\left(n^3\gamma_nh_n^{-1}\right)\left(\ln n\right)^{j+1}\right)\quad a.s.,
\end{eqnarray*}
and the application of~\eqref{alm2r} leads to
\begin{eqnarray*}
\left|r_n\left(x\right)-r\left(x\right)\right|&\leq&\left|\rho_n\left(x\right)-r\left(x\right)\right|+\left|\Delta_n\left(x\right)\right|\\
&=&O\left(m_n^{j+1}\left(n^3\gamma_nh_n^{-1}\right)\left(\ln n\right)^{j+1}\right)\quad a.s.
\end{eqnarray*}
Since~\eqref{arn3gh} ensures that~\eqref{a12ju} is satisfied for $j=0$, we have proved by induction that
\begin{eqnarray*}
\left|r_n\left(x\right)-r\left(x\right)\right|=O\left(m_n^{p_0}\left(n^3\gamma_nh_n^{-1}\right)\left(\ln n\right)^{p_0}\right)\quad a.s.
\end{eqnarray*}
Applying~\eqref{alm5d} with $\left(w_n\right)=\left(m_n^{p_0}\left(n^3\gamma_nh_n^{-1}\right)\left(\ln n\right)^{p_0}\right)$ and then~\eqref{alm2r}, we obtain
\begin{eqnarray*}
\left|r_n\left(x\right)-r\left(x\right)\right|=O\left(m_n^{p_0+1}\left(n^3\gamma_nh_n^{-1}\right)\left(\ln n\right)^{p_0+1}\right)+O\left(m_n\ln n\right)\quad a.s.
\end{eqnarray*}
Since the sequences $\left(m_n^{p_0+1}\left(n^3\gamma_nh_n^{-1}\right)\left(\ln n\right)^{p_0+1}\right)$ and $\left(m_n\ln n\right)$ are in $\mathcal{GS}\left(-m^*\left(p_0+1\right)+2+a\right)$ with $-m^*\left(p_0+1\right)+2+a<0$, and $\mathcal{GS}\left(-m^*\right)$ with $-m^*<0$ respectively, it follows that
\begin{eqnarray*}
\left|r_n\left(x\right)-r\left(x\right)\right|=O\left(\left(\ln n\right)^{-2}\right)\quad a.s.
\end{eqnarray*}
Applying once more~\eqref{alm5d} with $\left(w_n\right)=\left(\left(\ln n\right)^{-2}\right)\in\mathcal{GS}\left(0\right)$, we get
\begin{eqnarray*}
\left|\Delta_n\left(x\right)\right|=O\left(m_n\left(\ln n\right)^{-1}\right)+o\left(m_n\right)=o\left(m_n\right)\quad a.s.
\end{eqnarray*}
Theorem~\ref{aP:ra} (respectively Theorem~\ref{aT:rb}) in the case $\left(n\gamma_n\right)$ is bounded then follows from the application of Lemma~\ref{aL:rhatlc} (respectively Lemma~\ref{aL:rhatlog}).
\paragraph{Proof of Theorems~\ref{aP:ra} and~\ref{aT:rb} in the case $\lim_{n\to\infty}\left(n\gamma_n\right)=\infty$} $ $\\
In this case, Lemmas~\ref{aL:rhatlog} and~\ref{aL:rwn} imply that, for all $\delta>0$,
\begin{eqnarray}\label{alm2ra}
&\bullet& \left|\rho_n\left(x\right)-r\left(x\right)\right|=O\left(\tilde{m}_n\right)\quad a.s.\\
&\bullet& \mbox{If there exists} \left(w_n\right)\in \mathcal{GS}\left(w^*\right) \mbox{such that } \left|r_n\left(x\right)-r\left(x\right)\right|=O\left(w_n\right) a.s.,\nonumber\\
& &\mbox{then }\left|\Delta_n\left(x\right)\right|
=O\left(m_n\left(n^{1+\delta}\gamma_n\right)^{\frac{1+\delta}{2}}w_n\right)\quad a.s.\label{alm5da}
\end{eqnarray}
Now, set $\delta>0$ such that $c\left(\delta\right)=-m^*+\left(1+\delta\right)\left(1+\delta-\alpha\right)/2<0$ (the existence of such a $\delta$ being ensured by $\left(A2\right)$). In view of~\eqref{arn3gh}, the application of~\eqref{alm5da} with $\left(w_n\right)=\left(n^3\gamma_nh_n^{-1}\right)$ ensures that
\begin{eqnarray*}
\left|\Delta_n\left(x\right)\right|
&=&O\left(m_n\left(n^{1+\delta}\gamma_n\right)^{\frac{1+\delta}{2}}n^3\gamma_nh_n^{-1}\right)\quad a.s.,
\end{eqnarray*}
and, in view of~\eqref{alm2ra}, it follows that
\begin{eqnarray*}
\left|r_n\left(x\right)-r\left(x\right)\right|&=&O\left(\tilde{m}_n\right)+O \left(m_n\left(n^{1+\delta}\gamma_n\right)^{\frac{1+\delta}{2}}n^3\gamma_nh_n^{-1}\right)\quad a.s.
\end{eqnarray*}
Set $p\geq 1$, and assume that
\begin{eqnarray*}
\left|r_n\left(x\right)-r\left(x\right)\right|&=&O\left(\tilde{m}_n\right)+O\left(m_n^p\left(n^{1+\delta}\gamma_n\right)^{p\left(\frac{1+\delta}{2}\right)}n^3\gamma_nh_n^{-1}\right)\quad a.s.
\end{eqnarray*}
The application of~\eqref{alm5da} with $\left(w_n\right)=\left(\tilde{m}_n\right)$ and with $\left(w_n\right)=\left(m_n^p\left(n^{1+\delta}\gamma_n\right)^{p\left(\frac{1+\delta}{2}\right)}n^3\gamma_nh_n^{-1}\right)$ ensures that
\begin{eqnarray*}
\left|\Delta_n\left(x\right)\right|&=&O\left(m_n\left(n^{1+\delta}\gamma_n\right)^{\frac{1+\delta}{2}}\tilde{m}_n\right)+O\left(m_n^{p+1}\left(n^{1+\delta}\gamma_n\right)^{\left(p+1\right)\frac{1+\delta}{2}}n^3\gamma_nh_n^{-1}\right)\quad a.s.
\end{eqnarray*}
The sequence $\left(m_n\left(n^{1+\delta}\gamma_n\right)^{\frac{1+\delta}{2}}\right)$ being in $\mathcal{GS}\left(c\left(\delta\right)\right)$ with $c\left(\delta\right)<0$, it follows that
\begin{eqnarray*}
\left|\Delta_n\left(x\right)\right|&=&o\left(\tilde{m}_n\right)+O\left(m_n^{p+1}\left(n^{1+\delta}\gamma_n\right)^{\left(p+1\right)\frac{1+\delta}{2}}n^3\gamma_nh_n^{-1}\right)\quad a.s.
\end{eqnarray*}
and, in view of~\eqref{alm2ra}, we obtain
\begin{eqnarray*}
\left|r_n\left(x\right)-r\left(x\right)\right|&=&O\left(\tilde{m}_n\right)+O\left(m_n^{p+1}\left(n^{1+\delta}\gamma_n\right)^{\left(p+1\right)\frac{1+\delta}{2}}n^3\gamma_nh_n^{-1}\right)\quad a.s.
\end{eqnarray*}
We have thus proved by induction that, for all $p\geq 1$,
\begin{eqnarray*}
\left|r_n\left(x\right)-r\left(x\right)\right|&=&O\left(\tilde{m}_n\right)+O\left(m_n^{p}\left(n^{1+\delta}\gamma_n\right)^{p\left(\frac{1+\delta}{2}\right)}n^3\gamma_nh_n^{-1}\right)\quad a.s.
\end{eqnarray*}
By setting $p$ large enough, we deduce that
\begin{eqnarray*}
\left|r_n\left(x\right)-r\left(x\right)\right|&=&O\left(\tilde{m}_n\right)\quad a.s.
\end{eqnarray*}
Applying once more~\eqref{alm5da} with $\left(w_n\right)=\left(\tilde{m}_n\right)$, we get
\begin{eqnarray}\label{a25julys}
\left|\Delta_n\left(x\right)\right| & = & 
O\left(m_n\left(n^{1+\delta}\gamma_n\right)^{\frac{1+\delta}{2}}\tilde{m}_n\right)
\quad a.s.\\
 & = & 
o\left(m_n\right)\quad a.s. \nonumber
\end{eqnarray}
Theorem~\ref{aP:ra} (respectively Theorem~\ref{aT:rb}) in the case $\lim_{n\to\infty}\left(n\gamma_n\right)=\infty$ then follows from the application of Lemma~\ref{aL:rhatlc} (respectively Lemma~\ref{aL:rhatlog}).
\paragraph{Proof of Theorems~\ref{aAv:ra} and~\ref{aT:rba}} $ $\\
In view of~\eqref{a25julys}, and applying~\eqref{abt} and~\eqref{abteps}, we get, for all $\delta>0$,
\begin{eqnarray}\label{a26jul}
\left|\bar{\Delta}_n\left(x\right)\right|&=&O\left(\frac{1}{\sum_{k=1}^nq_k}
\sum_{k=1}^nq_k\tilde{m}_k^2 \left(k^{1+\delta}\gamma_k\right)^{\frac{1+\delta}{2}}\right)\quad a.s.\nonumber\\
&=&O\left(\frac{1}{nq_n}\left[n^{\delta}+nq_n\tilde{m}_n^2 \left(n^{1+\delta}\gamma_n\right)^{\frac{1+\delta}{2}}\right]\right)\quad a.s.\nonumber\\
&=&O\left(n^{\delta-1}q_n^{-1}+\tilde{m}_n^2 \left(n^{1+\delta}\gamma_n\right)^{\frac{1+\delta}{2}}\right)\quad a.s.
\end{eqnarray}
 $\bullet$ Let us first consider the case when the sequence $\left(nh_n^5\right)$ is bounded. In this case, we have $a\geq 1/5$, so that $a\geq \alpha/5$ and $m^*=\left(\alpha-a\right)/2$. Noting that $\left(A2\right)$ implies $a<3\alpha-2$, and applying $\left(A6\right)$, we can set $\delta>0$ such that 
\begin{eqnarray*}
\delta-\frac{\left(1+a\right)}{2}+q<0 \mbox{ and } \frac{\left(1-a\right)}{2}-2m^*+\frac{\left(1+\delta\right)}{2}\left(1+\delta-\alpha\right)<0.
\end{eqnarray*}
In view of~\eqref{a26jul}, we then obtain
\begin{eqnarray*}
\sqrt{nh_n}\left|\bar{\Delta}_n\left(x\right)\right|=o\left(1\right)\quad a.s.
\end{eqnarray*}
The first part of Theorem~\ref{aAv:ra} (respectively of Theorem~\ref{aT:rba}) then follows from the application of the first part of Lemma~\ref{aL:arb} (respectively of Lemma~\ref{aL:brb}).\\
$\bullet$ Let us now consider the case when $\lim_{n\to\infty}\left(nh_n^5\right)=\infty$. 
Noting that $\left(A2\right)$ then ensures that $6a<3\alpha-1$, and applying $\left(A6\right)$, we can set $\delta>0$ such that 
\begin{eqnarray*}
2a+\delta-1+q<0 \mbox{ and } 2a-2m^*+\frac{\left(1+\delta\right)}{2}\left(1+\delta-\alpha\right)<0.
\end{eqnarray*}
It then follows from~\eqref{a26jul} that
\begin{eqnarray*}
h_n^{-2}\left|\bar{\Delta}_n\left(x\right)\right|=o\left(1\right)\quad a.s.
\end{eqnarray*}
The second part of Theorem~\ref{aAv:ra} (respectively of Theorem~\ref{aT:rba}) then follows from the application of the second part of Lemma~\ref{aL:arb} (respectively of Lemma~\ref{aL:brb}).
\subsubsection{Proof of Theorems~\ref{aT:c} and~\ref{aT:Ac}} \label{asection 3.3.2}
Set
\begin{eqnarray}\label{a30bis}
B_n=n\gamma_nh_n^{-1}\ln n.
\end{eqnarray}

The proof of Theorems~\ref{aT:c} and~\ref{aT:Ac} relies on the repeted application of the following lemma, which is proved in Section~\ref{a4.6}.
\begin{lemma}\label{aL:rwn2}
Let $I$ be a bounded open interval on which $\varphi=\inf_{x\in I}f\left(x\right)>0$, let Assumptions $\left(A1\right)-\left(A4\right)$ hold for all $x \in I$, and assume that there exists $\left(w_n\right)\in \mathcal{GS}\left(w^*\right)$ such that $\sup_{x\in I}\left|r_n\left(x\right)-r\left(x\right)\right|=O\left(w_n\right)$ a.s. 
Moreover,\\
$\bullet$ in the case when $\left(n\gamma_n\right)$ is bounded, assume that $\lim_{n\to \infty} \left(n\gamma_n\right)>m^*/\varphi$ and that $w^*\geq0$;\\
$\bullet$ in the case when $\lim_{n\to\infty}\left(n\gamma_n\right)=\infty$, assume that the sequence $\left(w_n^{-1}B_n\sqrt{\gamma_nh_n^{-1}\ln n}\right)$ is bounded. Then, we have
\begin{eqnarray*}
\sup_{x\in I}\left|\Delta_n\left(x\right)\right|=O\left(m_nw_n\sqrt{\ln n}\right)\quad a.s.
\end{eqnarray*}
\end{lemma}
We first establish a preliminary upper bound of $\sup_{x\in I}\left(\left|r_n\left(x\right)-r\left(x\right)\right|\right)$ (which is better than the pointwise upper bound~\eqref{arn3gh} since the random variables $Y_k$ are assumed to have a finite exponential moment in Theorems~\ref{aT:c} and~\ref{aT:Ac}). Then, we successively prove Theorem~\ref{aT:c} in the case when $\left(n\gamma_n\right)$ is bounded, Theorem~\ref{aT:c} in the case when $\lim_{n\to\infty}\left(n\gamma_n\right)=\infty$, and finally Theorem~\ref{aT:Ac}.
\paragraph{Preliminary upper bound.} $ $\\
Proceeding as for the proof of~\eqref{arn3gh}, we note that, for all $n\geq n_0$,\begin{eqnarray*}
\sup_{x\in I}\left|r_n\left(x\right)\right|\leq \sup_{x\in I}\left|r_{n_0-1}\left(x\right)\right| +\left(\sup_{k \leq n}\left|Y_k\right|\right)\|K\|_{\infty}\sum_{k=1}^n\gamma_kh_k^{-1},
\end{eqnarray*}
with, this time, in view of $\left(A4\right)$,
\begin{eqnarray*}
\mathbb{P}\left[\sup_{k\leq n}\left|Y_k\right|>\frac{3}{t^*}\ln n\right]
 \leq n\mathbb{P}\left[\exp\left(t^*\left|Y\right|\right)>n^3\right]\leq n^{-2}\mathbb{E}\left(\exp\left(t^*|Y|\right)\right).
\end{eqnarray*}
We deduce that
\begin{eqnarray*}
\sup_{x\in I}\left|r_n\left(x\right)-r\left(x\right)\right|=O\left(B_n\right)\quad a.s.
\end{eqnarray*}
\paragraph{Proof of Theorem~\ref{aT:c} in the case $\left(n\gamma_n\right)$ is bounded} $ $\\
In this case, we have $\alpha=1$, $\left(B_n\right)\in\mathcal{GS}\left(a\right)$ (with $a>0$); the application of Lemma~\ref{aL:rwn2} with $\left(w_n\right)=\left(B_n\right)$ ensures that
\begin{eqnarray*}
\sup_{x\in I}\left|\Delta_n\left(x\right)\right|=O\left(m_nB_n\sqrt{\ln n}\right)\quad a.s.
\end{eqnarray*}
Applying Lemma~\ref{aL:Tuhat}, we get
\begin{eqnarray*}
\sup_{x \in I}\left|r_n\left( x\right)-r\left(x\right)\right|& \leq &\sup_{x \in I}\left|\rho_n\left( x\right)-r\left(x\right)\right|+\sup_{x \in I}\left|\Delta_n\left( x\right)\right|\\
&=&O\left(m_nB_n\sqrt{\ln n}\right) \quad a.s.
\end{eqnarray*}
Since $\left(m_nB_n\sqrt{\ln n}\right)\in \mathcal{GS}\left(-m^*+a\right)$ (with $-m^*+a<0$), it follows that
\begin{eqnarray*}
\sup_{x \in I}\left|r_n\left( x\right)-r\left(x\right)\right|&=&O\left(\left[\ln n\right]^{-1}\right) \quad a.s.
\end{eqnarray*}
Applying once more Lemma~\ref{aL:rwn2} with $\left(w_n\right)=\left(\left[\ln n\right]^{-1}\right)$, we get
\begin{eqnarray*}
\sup_{x \in I}\left|\Delta_n\left( x\right)\right|&=&O\left(m_n\left(\ln n\right)^{-\frac{1}{2}}\right)\quad a.s.\nonumber\\
&=&o\left(\tilde{m}_n\right)\quad a.s.
\end{eqnarray*}
Theorem~\ref{aT:c} in the case when $\left(n\gamma_n\right)$ is bounded then follows from the application of Lemma~\ref{aL:Tuhat}.
\paragraph{Proof of Theorem~\ref{aT:c} in the case $\lim_{n\to\infty}\left(n\gamma_n\right)=\infty$} $ $\\
The sequence $\left(\sqrt{\gamma_nh_n^{-1}\ln n}\right)$ being clearly bounded, we can apply Lemma~\ref{aL:rwn2} with $\left(w_n\right)=\left(B_n\right)$; we then obtain 
\begin{eqnarray*}
\sup_{x\in I}\left|\Delta_n\left(x\right)\right|=O\left(m_nB_n\sqrt{\ln n}\right)\quad a.s.
\end{eqnarray*}
The application of Lemma~\ref{aL:Tuhat} then ensures that
\begin{eqnarray*}
\sup_{x \in I}\left|r_n\left( x\right)-r\left(x\right)\right|
&=&O\left(\tilde{m}_n\right)+O\left(m_nB_n\sqrt{\ln n}\right) \quad a.s.\\
&=&O\left(m_nB_n\sqrt{\ln n}\right) \quad a.s.
\end{eqnarray*}
Since $\left(m_nB_n\sqrt{\ln n}\right)^{-1}B_n\sqrt{\gamma_nh_n^{-1}\ln n}=m_n^{-1}\sqrt{\gamma_nh_n^{-1}}=O\left(1\right)$, we can apply once more Lemma~\ref{aL:rwn2} with $\left(w_n\right)=\left(m_nB_n\sqrt{\ln n}\right)$; we get
\begin{eqnarray}\label{a25julymp}
\sup_{x \in I}\left|\Delta_n\left( x\right)\right|&=&O\left(m_n^2B_n\ln n\right) \quad a.s.
\end{eqnarray}
Noting that $\left(m_nB_n\ln n\right)\in \mathcal{GS}\left(-m^*+1-\alpha+a\right)$ with, in view of $\left(A4\right)iii)$, $-m^*+1-\alpha+a<0$, it follows that
\begin{eqnarray*}
\sup_{x \in I}\left|\Delta_n\left( x\right)\right|&=&o\left(m_n\right) \quad a.s.
\end{eqnarray*}
Theorem~\ref{aT:c} in the case when $\lim_{n\to\infty}\left(n\gamma_n\right)=\infty$ then follows from the application of Lemma~\ref{aL:Tuhat}.
\paragraph{Proof of Theorem~\ref{aT:Ac}} $ $\\
\begin{itemize}
\item In the case when the sequence $\left(nh_n^5/\ln n\right)$ is bounded, we have, in view of~\eqref{a25julymp},
\begin{eqnarray*}
\sqrt{nh_n}\left(\ln n\right)^{-1}\sup_{x\in I}\left|\Delta_n\left(x\right)\right|=O\left(\sqrt{nh_n}m_n^2B_n\right)\quad a.s.
\end{eqnarray*}
Now, in this case, we have $a\geq1/5\geq \alpha/5$ and thus $m^*=\left(\alpha-a\right)/2$. It follows that $\left(\sqrt{nh_n}m_n^2B_n\right)\in\mathcal{GS}\left(3\left(1+a\right)/2-2\alpha\right)$ with $3\left(1+a\right)/2-2\alpha<0$, and thus 
\begin{eqnarray*}
\sup_{x\in I}\left|\Delta_n\left(x\right)\right|=O\left(\sqrt{n^{-1}h_n^{-1}}\,\ln n\right)\quad a.s.
\end{eqnarray*}
The first part of Theorem~\ref{aT:Ac} then follows from the application of Lemma~\ref{aL:Tubar}.\\
\item In the case when $\lim_{n\to\infty}\left(nh_n^5/\ln n\right)=\infty$,~\eqref{a25julymp} 
ensures that
\begin{eqnarray*}
h_n^{-2}\sup_{x\in I}\left|\Delta_n\left(x\right)\right|=O\left(h_n^{-2}m_n^2B_n\ln n\right)\quad a.s.,
\end{eqnarray*}
with $\left(h_n^{-2}m_n^2B_n\ln n\right)\in\mathcal{GS}\left(3a-2m^*+1-\alpha\right)$. Noting that the assumptions of Theorem~\ref{aT:Ac} ensure that $3a-2m^*+1-\alpha<0$, we deduce that
\begin{eqnarray*}
\sup_{x\in I}\left|\Delta_n\left(x\right)\right|=O\left(h_n^{2}\right)\quad a.s.\end{eqnarray*}
The second part of Theorem~\ref{aT:Ac} then follows from the application of Lemma~\ref{aL:Tubar}.
\end{itemize}

\section{Proof of Lemmas}\label{alem}
\subsection{Proof of Lemma~\ref{a2bis}}\label{a4.0}
We establish that, under the condition $\lim_{n\to\infty}\left(n\gamma_n\right)>\left(\alpha-a\right)/\left(2f\left(x\right)\right)$,
\begin{eqnarray}\label{aTrevtlc}
&\bullet& \mbox{ if $a\geq\alpha/5$, then }\nonumber\\
&&\sqrt{\gamma_n^{-1}h_n} \left(\tilde{T}_{n}\left( x\right)-\mathbb{E}\left(\tilde{T}_{n}\left(x\right)\right)\right)\stackrel{\mathcal{D}}{\rightarrow}\mathcal{N}\left(0,\frac{Var\left[Y\vert X=x\right]f\left(x\right)}{\left(2f\left(x\right)-\left(\alpha-a\right)\xi\right)}\int_{\mathbb{R}} K^2\left(z\right)dz\right),\\
&\bullet&  \mbox{ if $a>\alpha/5$, then } \sqrt{\gamma_n^{-1}h_n}\mathbb{E}\left(\tilde{T}_n\left(x\right)\right)\rightarrow 0,\label{aTnE1}
\end{eqnarray}
and prove that, under the condition $\lim_{n\to\infty}\left(n\gamma_n\right)>2a/f\left(x\right)$,
\begin{eqnarray}
&\bullet&  \mbox{ if $a\leq\alpha/5$, then } h_n^{-2}\mathbb{E}\left(\tilde{T}_n\left(x\right)\right)\rightarrow \frac{f\left(x\right)m^{\left(2\right)}\left(x\right)}{f\left(x\right)-2a\xi},\label{aTnE}\\
&\bullet & \mbox{ if $a<\alpha/5$, then } h_n^{-2}\left(\tilde{T}_n\left(x\right)-\mathbb{E}\left(\tilde{T}_n\left(x\right)\right)\right)\stackrel{\mathbb{P}}{\rightarrow}0\label{aTnP}.
\end{eqnarray}
As a matter of fact the combination of~\eqref{aTrevtlc} and~\eqref{aTnE1} (respectively of~\eqref{aTrevtlc} and~\eqref{aTnE}) gives Part 1 of Lemma~\ref{a2bis} in the case $a>\alpha/5$ (respectively $a=\alpha/5$), that of~\eqref{aTnE} and~\eqref{aTnP} (respectively of~\eqref{aTrevtlc} and~\eqref{aTnE}) gives Part 2 of Lemma~\ref{a2bis} in the case $a<\alpha/5$ (respectively $a=\alpha/5$). We prove~\eqref{aTrevtlc},~\eqref{aTnP},~\eqref{aTnE}, and~\eqref{aTnE1} successively.
\paragraph{Proof of~\eqref{aTrevtlc}}
Set
\begin{eqnarray}\label{aynk}
\tilde{\eta}_{k}\left(x\right)=\Pi_k^{-1}\left(f\left(x\right)\right)\gamma_k\left[\left(W_k\left(x\right)-r\left(x\right)Z_{k}\left(x\right)\right)\right],
\end{eqnarray}
so that $\tilde{T}_n\left(x\right)-\mathbb{E}\left(\tilde{T}_n\left(x\right)\right)=\Pi_n\left(f\left(x\right)\right) \sum_{k=n_{0}}^n\left[\tilde{\eta}_k\left(x\right)-\mathbb{E}\left(\tilde{\eta_k}\left(x\right)\right)\right]$. We have
\begin{eqnarray*}
Var\left(\tilde{\eta}_{k}\left(x\right)\right)&=&\Pi_k^{-2}\left(f\left(x\right)\right)\gamma_k^2\left[Var\left(W_k\left(x\right)\right)+r^2\left(x\right)Var\left(Z_k\left(x\right)\right)-2r\left(x\right)Cov\left(W_k\left(x\right),Z_k\left(x\right)\right)\right].
\end{eqnarray*}
In view of $\left(A3\right)$, classical computations give
\begin{eqnarray}
Var\left(W_k\left(x\right)\right)&=&\frac{1}{h_k}\left[\mathbb{E}\left[Y^2\vert X=x\right]f\left(x\right)\int_{\mathbb{R}}K^2\left(z\right)dz+o\left(1\right)\right],\label{avarw}\\
Var\left(Z_k\left(x\right)\right)&=&\frac{1}{h_k}\left[f\left(x\right)\int_{\mathbb{R}}K^2\left(z\right)dz+o\left(1\right)\right],\label{avarz}\\
Cov\left(W_k\left(x\right),Z_k\left(x\right)\right)&=&\frac{1}{h_k}\left[r\left(x\right)f\left(x\right)\int_{\mathbb{R}}K^2\left(z\right)dz+o\left(1\right)\right]\label{acovwz}.
\end{eqnarray}
It follows that
\begin{eqnarray}\label{a28jun}
Var\left(\tilde{\eta}_{k}\left(x\right)\right)&=&\frac{\Pi_k^{-2}\left(f\left(x\right)\right)\gamma_k^2}{h_k}\left[Var\left[Y\vert X=x\right]f\left(x\right)\int_{\mathbb{R}} K^2\left(z\right)dz+o\left(1\right)\right],
\end{eqnarray}
and, since $\lim_{n\to\infty}\left(n\gamma_n\right)>\left(\alpha-a\right)/\left(2f\left(x\right)\right)$, Lemma~\ref{alap} ensures that
\begin{eqnarray}\label{avnti} 
v_n^2&=&\sum_{k=n_{0}}^nVar\left(\tilde{\eta}_{k}\left(x\right)\right)\nonumber\\
&=&\sum_{k=n_{0}}^n\frac{\Pi_k^{-2}\left(f\left(x\right)\right)\gamma_k^2}{h_k}\left[Var\left[Y\vert X=x\right]f\left(x\right)\int_{\mathbb{R}}K^2\left(z\right)dz+o\left(1\right)\right]\nonumber\\
&=&\frac{1}{\Pi_n^2\left(f\left(x\right)\right)}\frac{\gamma_n}{h_n}\frac{1}{2f\left(x\right)-\left(\alpha-a\right)\xi}\left[Var\left[Y\vert X=x\right]f\left(x\right)\int_{\mathbb{R}}K^2\left(z\right)dz+o\left(1\right)\right].
\end{eqnarray}
For all $p\in\left]0,1\right]$ and in view of $\left(A3\right)$, we have
\begin{eqnarray}\label{aMok}
\lefteqn{\mathbb{E}\left(\left|Y_k-r\left(x\right)\right|^{2+p}K^{2+p}\left(\frac{x-X_k}{h_k}\right)\right)}\nonumber\\
&=&h_k\int_{\mathbb{R}^2}\left|y-r\left(x\right)\right|^{2+p}K^{2+p}\left(s\right)g\left(x-h_ks,y\right)dyds\nonumber\\
&\leq &2^{1+p}h_k\int_{\mathbb{R}}\left\{\int_{\mathbb{R}}\left|y\right|^{2+p}g\left(x-h_ks,y\right)dy+\left|r\left(x\right)\right|^{2+p}\int_{\mathbb{R}}g\left(x-h_ks,y\right)dy\right\}K^{2+p}\left(s\right)ds\nonumber\\
&=&O\left(h_k\right).
\end{eqnarray}
Now, set $p\in\left]0,1\right]$ such that $\lim_{n\to\infty}\left(n\gamma_n\right)>\left(1+p\right)\left(\alpha-a\right)/\left(\left(2+p\right)f\left(x\right)\right)$. Applying Lemma~\ref{alap}, we get 
\begin{eqnarray}\label{aet2p}
\sum_{k=n_0}^n\mathbb{E}\left[\left|\tilde{\eta}_{k}\left(x\right)\right|^{2+p}\right]
& = & O\left(\sum_{k=n_0}^n \frac{\Pi_k^{-2-p}\left(f\left(x\right)\right)\gamma_k^{2+p}}{h_k^{2+p}}\mathbb{E}\left(\left|Y_k-r\left(x\right)\right|^{2+p}K^{2+p}\left(\frac{x-X_k}{h_k}\right)\right)\right)\nonumber\\
& = & O\left(\sum_{k=n_0}^n \frac{\Pi_k^{-2-p}\left(f\left(x\right)\right)\gamma_k^{2+p}}{h_k^{1+p}}\right)\nonumber\\
&=&O\left(\frac{1}{\Pi_n^{2+p}\left(f\left(x\right)\right)}\frac{\gamma_n^{1+p}}{h_n^{1+p}}\right).
\end{eqnarray}
Using~\eqref{avnti}, we deduce that
\begin{eqnarray*}
\frac{1}{v_n^{2+p}}\sum_{k=n_0}^n\mathbb{E}\left[\left|\tilde{\eta_{k}}\left(x\right)\right|^{2+p}\right]& = &O\left(\left(\frac{\gamma_n}{h_n}\right)^{\frac{p}{2}}\right)\\
&=&o\left(1\right),
\end{eqnarray*}
and~\eqref{aTrevtlc} follows by application of Lyapounov Theorem.
\paragraph{Proof of~\eqref{aTnP}}
In view of~\eqref{a28jun}, and since $a<\al/5$ and $\lim_{n\to\infty}\left(n\gamma_n\right)>2a/f\left(x\right)$, the application of Lemma~\ref{alap} ensures that
\begin{eqnarray*}
Var\left(\tilde{T}_n\left(x\right)\right)&=&\Pi_n^2\left(f\left(x\right)\right)\sum_{k=n_0}^n\frac{\Pi_k^{-2}\left(f\left(x\right)\right)\gamma_k^2}{h_k}\left[Var\left[Y\vert X=x\right]f\left(x\right)\int_{\mathbb{R}}K^2\left(z\right)dz+o\left(1\right)\right]\nonumber\\
&=&\Pi_n^2\left(f\left(x\right)\right)\sum_{k=n_0}^n\Pi_k^{-2}\left(f\left(x\right)\right)\gamma_ko\left(h_k^4\right)\nonumber\\
&=&o\left(h_n^4\right),
\end{eqnarray*}
which gives~\eqref{aTnP}.
\paragraph{Proof of~\eqref{aTnE}}
We have
\begin{eqnarray*}
\mathbb{E}\left(\tilde{T}_n\left(x\right)\right)&=&\Pi_n\left(f\left(x\right)\right)\sum_{k=n_0}^{n}\Pi_k^{-1}\left(f\left(x\right)\right)\gamma_k\left[\left(\mathbb{E}\left(W_k\left(x\right)\right)-a\left(x\right)\right)-r\left(x\right)\left(\mathbb{E}\left(Z_{k}\left(x\right)\right)-f\left(x\right)\right)\right].
\end{eqnarray*}
In view of $\left(A3\right)$ we obtain
\begin{eqnarray}
\mathbb{E}\left(W_k\left(x\right)\right)-a\left(x\right)&=&\frac{1}{2}h_k^2\int_{\mathbb{R}}y\frac{\partial^2 g}{\partial x^2}\left(x,y\right)dy\left[1+o\left(1\right)\right]\int_{\mathbb{R}}z^2K\left(z\right)dz,\label{aespw}\\
\mathbb{E}\left(Z_k\left(x\right)\right)-f\left(x\right)&=&\frac{1}{2}h_k^2\int_{\mathbb{R}}\frac{\partial^2 g}{\partial x^2}\left(x,y\right)dy\left[1+o\left(1\right)\right]\int_{\mathbb{R}}z^2K\left(z\right)dz\label{aespz}.
\end{eqnarray}
Since $\lim_{n\to\infty}\left(n\gamma_n\right)>2a/f\left(x\right)$, it follows from the application of Lemma~\ref{alap} that
\begin{eqnarray*}
\mathbb{E}\left(\tilde{T}_n\left(x\right)\right)&=&\Pi_n\left(f\left(x\right)\right)\sum_{k=n_0}^{n} \Pi_k^{-1}\left(f\left(x\right)\right)\gamma_kh_k^2\nonumber\\
&&\left[\frac{1}{2}\left(\int_{\mathbb{R}}y\frac{\partial^2 g}{\partial x^2}\left(x,y\right)dy-r\left(x\right)\int_{\mathbb{R}}\frac{\partial^2 g}{\partial x^2}\left(x,y\right)dy\right)+o\left(1\right)\right]\int_{\mathbb{R}}z^2K\left(z\right)dz\nonumber\\
&=&\Pi_n\left(f\left(x\right)\right)\sum_{k=n_0}^{n}\Pi_k^{-1}\left(f\left(x\right)\right)\gamma_kh_k^2\left[m^{\left(2\right)}\left(x\right)f\left(x\right)+o\left(1\right)\right]\nonumber\\
&=&\frac{1}{f\left(x\right)-2a\xi}h_n^2\left[m^{\left(2\right)}\left(x\right)f\left(x\right)+o\left(1\right)\right],
\end{eqnarray*}
which gives~\eqref{aTnE}.
\paragraph{Proof of~\eqref{aTnE1}}
Since $a>\alpha/5$ and $\lim_{n\to\infty}\left(n\gamma_n\right)>\left(1-a\right)/\left(2f\left(x\right)\right)$, we have
\begin{eqnarray*}
\mathbb{E}\left(\tilde{T}_n\left(x\right)\right)&=&\Pi_n\left(f\left(x\right)\right)\sum_{k=n_0}^{n}\Pi_k^{-1}\left(f\left(x\right)\right)\gamma_ko\left(\sqrt{\gamma_kh_k^{-1}}\right)\nonumber\\
&=&o\left(\sqrt{\gamma_nh_n^{-1}}\right),
\end{eqnarray*}
which gives~\eqref{aTnE1}.
\subsection{Proof of Lemma~\ref{a3bis}}\label{aprlb}
Set
\begin{eqnarray*}
S_n\left(x\right)&=&\sum_{k=n_0}^n\left[\tilde{\eta_{k}}\left(x\right)-\mathbb{E}\left(\tilde{\eta_{k}}\left(x\right)\right)\right].
\end{eqnarray*}
where $\tilde{\eta_k}$ is defined in~\eqref{aynk}.\\
$\bullet$ Let us first consider the case $a\geq \alpha/5$ (in which case $\lim_{n\to\infty}\left(n\gamma_n\right)>\left(\alpha-a\right)/\left(2f\left(x\right)\right)$). We set $H_n^2\left(f\left(x\right)\right)=\Pi_n^2\left(f\left(x\right)\right)\gamma_n^{-1}h_n$, and note that, since $\left(\gamma_n^{-1}h_n\right)\in \mathcal{GS}\left(\alpha-a\right)$, we have
\begin{eqnarray}\label{a23ju}
\ln \left(H_n^{-2}\left(f\left(x\right)\right)\right)&=&-2\ln \left(\Pi_n\left(f\left(x\right)\right)\right)+\ln\left(\prod_{k=n_{0}}^n\frac{\gamma_{k-1}^{-1}h_{k-1}}{\gamma_k^{-1}h_k}\right)+\ln \left(\gamma_{n_0-1}h_{n_0-1}^{-1}\right)\nonumber\\
&=&-2\sum_{k=n_0}^n\ln\left(1-f\left(x\right)\gamma_k\right)+\sum_{k=n_0}^n\ln\left(1-\frac{\alpha-a}{k}+o\left(\frac{1}{k}\right)\right)+\ln \left(\gamma_{n_0-1}h_{n_0-1}^{-1}\right)\nonumber\\
&=&\sum_{k=n_0}^n\left(2f\left(x\right)\gamma_k+o\left(\gamma_k\right)\right)-\sum_{k=n_0}^n\left(\left(\alpha-a\right)\xi\gamma_k+o\left(\gamma_k\right)\right)+\ln \left(\gamma_{n_0-1}h_{n_0-1}^{-1}\right)\nonumber\\
&=&\left(2f\left(x\right)-\xi\left(\alpha-a\right)\right)s_n+o\left(s_n\right).
\end{eqnarray}
Since $2f\left(x\right)-\xi\left(\alpha-a\right)>0$, it follows in particular that $\lim_{n\to\infty}H_n^{-2}\left(f\left(x\right)\right)=\infty$. Moreover, we clearly have $\lim_{n\to\infty}H_n^2\left(f\left(x\right)\right)/H_{n-1}^2\left(f\left(x\right)\right)=1$, and by~\eqref{avnti},
\begin{eqnarray*}
\lim_{n\to\infty}H_n^2\left(f\left(x\right)\right)\sum_{k=n_0}^nVar\left[\tilde{\eta}_k\left(x\right)\right]=\left[2f\left(x\right)-\left(\alpha-a\right)\xi\right]^{-1}Var\left[Y\vert X=x\right]f\left(x\right)\int_{\mathbb{R}}K^2\left(z\right)dz,
\end{eqnarray*}
and, in view of~\eqref{aMok}, $\mathbb{E}\left[\left|\tilde{\eta}_{n}\left(x\right)\right|^{3}\right]=O\left(\Pi_n^{-3}\left(f\left(x\right)\right)\gamma_n^3h_n^{-2}\right)$. Now, since $\left(\gamma_n^{-1}h_n\right)\in \mathcal{GS}\left(\alpha-a\right)$, applying Lemma~\ref{alap} and~\eqref{a23ju}, we get
\begin{eqnarray*}
\frac{1}{n\sqrt{n}}\sum_{k=n_0}^n\mathbb{E}\left(\left|H_n\left(f\left(x\right)\right)\tilde{\eta}_k\left(x\right)\right|^3\right)&=&O\left(\frac{H_n^3\left(f\left(x\right)\right)}{n\sqrt{n}}\left(\sum_{k=n_0}^n \frac{\Pi_k^{-3}\left(f\left(x\right)\right)\gamma_k^{3}}{h_k^{2}}\right)\right)\\
&=&O\left(\frac{\Pi_n^3\left(f\left(x\right)\right)\gamma_n^{-\frac{3}{2}}h_n^{\frac{3}{2}}}{n\sqrt{n}}\left(\sum_{k=n_0}^n \Pi_k^{-3}\left(f\left(x\right)\right)\gamma_ko\left(\left(\gamma_kh_k^{-1}\right)^{\frac{3}{2}}\right)\right)\right)\\
&=&o\left(\frac{1}{n\sqrt{n}}\right)\\
&=&o\left(\left[\ln \left(H_n^{-2}\left(f\left(x\right)\right)\right)\right]^{-1}\right).
\end{eqnarray*}
The application of Theorem 1 in Mokkadem and Pelletier (2008) then ensures that, with probability one, the sequence

\begin{eqnarray*}
\left(\frac{H_n\left(f\left(x\right)\right)S_n\left(x\right)}{\sqrt{2\ln \ln \left(H_n^{-2}\left(f\left(x\right)\right)\right)}}\right)=\left(\frac{\sqrt{\gamma_n^{-1}h_n}\left(\tilde{T}_{n}\left( x\right)-\mathbb{E}\left(\tilde{T}_n\left(x\right)\right)\right)}{\sqrt{2\ln \ln \left(H_n^{-2}\left(f\left(x\right)\right)\right)}}\right)
\end{eqnarray*}
is relatively compact and its limit set is the interval 
\begin{eqnarray}\label{a10jul}
\left[-\sqrt{\frac{Var\left[Y\vert X=x\right]f\left(x\right)}{\left(2f\left(x\right)-\left(\alpha-a\right)\xi\right)}\int_{\mathbb{R}} K^2\left(z\right)dz},\sqrt{\frac{Var\left[Y\vert X=x\right]f\left(x\right)}{\left(2f\left(x\right)-\left(\alpha-a\right)\xi\right)}\int_{\mathbb{R}} K^2\left(z\right)dz}\right].
\end{eqnarray} 
In view of~\eqref{a23ju}, we have $\lim_{n\to\infty}\ln \ln \left(H_n^{-2}\left(f\left(x\right)\right)\right)/ \ln s_n=1$. It follows that, with probability one, the sequence $\left(\sqrt{\gamma_n^{-1}h_n}\left(\tilde{T}_n\left(x\right)-\mathbb{E}\left(\tilde{T}_n\left(x\right)\right)\right)/\sqrt{2\ln s_n}\right)$ is relatively compact, and its limit set is the interval given in~\eqref{a10jul}. The application of~\eqref{aTnE1} (respectively of~\eqref{aTnE}) concludes the proof of Lemma~\ref{a3bis} in the case $a>\alpha/5$ (respectively $a=\alpha/5$).\\
$\bullet$ Let us now consider the case $a<\alpha/5$ (in which case $\lim_{n\to\infty}\left(n\gamma_n\right)>2a/f\left(x\right)$). Set $H_n^{-2}\left(f\left(x\right)\right)=\Pi_n^{-2}\left(f\left(x\right)\right)h_n^{4}\left(\ln \ln \left(\Pi_n^{-2}\left(f\left(x\right)\right)h_n^4\right)\right)^{-1}$, and note that, since $\left(h_n^{-4}\right)\in \mathcal{GS}\left(4a\right)$, we have
\begin{eqnarray}\label{a29ju}
\ln \left(\Pi_n^{-2}\left(f\left(x\right)\right)h_n^{4}\right)&=&-2\ln \left(\Pi_n\left(f\left(x\right)\right)\right)+\ln\left(\prod_{k=n_0}^n\frac{h_{k-1}^{-4}}{h_k^{-4}}\right)+\ln \left(h_{n_0-1}^4\right)\nonumber\\
&=&-2\sum_{k=n_0}^n\ln\left(1-\gamma_kf\left(x\right)\right)+\sum_{k=n_0}^n\ln\left(1-\frac{4a}{k}+o\left(\frac{1}{k}\right)\right)+\ln \left(h_{n_0-1}^4\right)\nonumber\\
&=&\sum_{k=n_0}^n\left(2\gamma_kf\left(x\right)+o\left(\gamma_k\right)\right)-\sum_{k=n_0}^n\left(4a\xi\gamma_k+o\left(\gamma_k\right)\right)+\ln \left(h_{n_0-1}^4\right)\nonumber\\
&=&\left(2f\left(x\right)-4a\xi\right)s_n+o\left(s_n\right).
\end{eqnarray}
Since $2f\left(x\right)-4a\xi>0$, it follows in particular that $\lim_{n\to\infty}\Pi_n^{-2}\left(f\left(x\right)\right)h_n^4=\infty$, and thus \\$\lim_{n\to\infty}H_n^{-2}\left(f\left(x\right)\right)=\infty$. Moreover, we clearly have $\lim_{n\to\infty}H_n^2\left(f\left(x\right)\right)/H_{n-1}^2\left(f\left(x\right)\right)=1$. Now, set $\epsilon \in \left]0,\alpha-5a\right[$ such that $\lim_{n\to\infty}\left(n\gamma_n\right)>2a/f\left(x\right)+\epsilon/2$; in view of~\eqref{avnti}, and applying Lemma~\ref{alap}, we get
\begin{eqnarray*}
H_n^2\left(f\left(x\right)\right)\sum_{k=n_0}^nVar\left[\tilde{\eta}_k\left(x\right)\right]&=&O\left(\Pi_n^2\left(f\left(x\right)\right)h_n^{-4}\ln \ln \left(\Pi_n^{-2}\left(f\left(x\right)\right)h_n^{4}\right)\sum_{k=n_0}^n\frac{\Pi_k^{-2}\left(f\left(x\right)\right)\gamma_k^2}{h_k}\right)\\
&=&O\left(\Pi_n^2\left(f\left(x\right)\right)h_n^{-4}\ln \ln \left(\Pi_n^{-2}\left(f\left(x\right)\right)h_n^{4}\right)\sum_{k=n_0}^n\Pi_k^{-2}\left(f\left(x\right)\right)\gamma_ko\left(h_k^4k^{-\epsilon}\right)\right)\\
&=&o\left(1\right).
\end{eqnarray*}
Moreover, in view of~\eqref{aMok} we have $\mathbb{E}\left[\left|\tilde{\eta}_{n}\left(x\right)\right|^{3}\right]=O\left(\Pi_n^{-3}\left(f\left(x\right)\right)\gamma_n^3h_n^{-2}\right)$, and thus in view of~\eqref{a29ju}, we get
\begin{eqnarray*}
\lefteqn{\frac{1}{n\sqrt{n}}\sum_{k=n_0}^n\mathbb{E}\left(\left|H_n\left(f\left(x\right)\right)\tilde{\eta}_k\left(x\right)\right|^3\right)}\\
&=&O\left(\frac{\Pi_n^3\left(f\left(x\right)\right)h_n^{-6}}{n\sqrt{n}}\left(\ln \ln \left(\Pi_n^{-2}\left(f\left(x\right)\right)h_n^4\right)\right)^{\frac{3}{2}}\left(\sum_{k=n_0}^n \frac{\Pi_k^{-3}\left(f\left(x\right)\right)\gamma_k^{3}}{h_k^{2}}\right)\right)\\
&=&O\left(\frac{\Pi_n^3\left(f\left(x\right)\right)h_n^{-6}}{n\sqrt{n}}\left(\ln \ln \left(\Pi_n^{-2}\left(f\left(x\right)\right)h_n^4\right)\right)^{\frac{3}{2}}\left(\sum_{k=n_0}^n \Pi_k^{-3}\left(f\left(x\right)\right)\gamma_ko\left(h_k^{6}\right)\right)\right)\\
&=&o\left(\frac{\left(\ln \ln \left(\Pi_n^{-2}\left(f\left(x\right)\right)h_n^4\right)\right)^{\frac{3}{2}}}{n\sqrt{n}}\right)\\
&=&o\left(\left[\ln \left(H_n^{-2}\left(f\left(x\right)\right)\right)\right]^{-1}\right).
\end{eqnarray*}
The application of Theorem 1 in Mokkadem and Pelletier (2008) then ensures that, with probability one, 
\begin{eqnarray*}
\lim_{n\to\infty}\frac{H_n\left(f\left(x\right)\right)S_n\left(x\right)}{\sqrt{2\ln \ln \left(H_n^{-2}\left(f\left(x\right)\right)\right)}}=\lim_{n\to\infty}h_n^{-2}\frac{\sqrt{\ln \ln \left(\Pi_n^{-2}\left(f\left(x\right)\right)h_n^4\right)}}{\sqrt{2\ln \ln \left(H_n^{-2}\left(f\left(x\right)\right)\right)}}\left(\tilde{T}_n\left(x\right)-\mathbb{E}\left(\tilde{T}_n\left(x\right)\right)\right)=0.
\end{eqnarray*}
Noting that~\eqref{a29ju} ensures that $\lim_{n\to\infty}\ln \ln \left(H_n^{-2}\left(f\left(x\right)\right)\right)/\ln \ln \left(\Pi_n^{-2}\left(f\left(x\right)\right)h_n^4\right)=1$, we get
\begin{eqnarray*}
\lim_{n\to\infty}h_n^{-2}\left[\tilde{T}_{n}\left( x\right)-\mathbb{E}\left(\tilde{T}_n\left(x\right)\right)\right]=0\quad a.s.,
\end{eqnarray*}
and Lemma~\ref{a3bis} in the case $a<\alpha/5$ follows from~\eqref{aTnE}. 
\subsection{Proof of Lemma~\ref{a4bis}}\label{a4.1}
 Let us write 
\begin{eqnarray*}
\tilde{T}_n\left(x\right)=\tilde{T}_n^{\left(1\right)}\left(x\right)-r\left(x\right)\tilde{T}_n^{\left(2\right)}\left(x\right)
\end{eqnarray*}
with
\begin{eqnarray*}
\tilde{T}_n^{\left(1\right)}\left(x\right)&=&\sum_{k=n_0}^nU_{k,n}\left(f\left(x\right)\right)\gamma_kh_k^{-1}Y_kK\left(\frac{x-X_k}{h_k}\right)\\
\tilde{T}_n^{\left(2\right)}\left(x\right)&=&\sum_{k=n_0}^nU_{k,n}\left(f\left(x\right)\right)\gamma_kh_k^{-1}K\left(\frac{x-X_k}{h_k}\right).
\end{eqnarray*}
Lemma~\ref{a4bis} is proved by showing that, for $i\in \left\{1,2\right\}$, under the condition $\lim_{n\to\infty}\left(n\gamma_n\right)>\left(\alpha-a\right)/\left(2\varphi \right)$,
\begin{eqnarray}\label{a1/08/a}
&\bullet& \mbox{ if $a\geq\alpha/5$, then } \sup_{x\in I} \left|\tilde{T}_{n}^{\left(i\right)}\left( x\right)-\mathbb{E}\left(\tilde{T}_{n}^{\left(i\right)}\left(x\right)\right)\right|=O\left(\sqrt{\gamma_nh_n^{-1}}\,\ln n\right)\quad a.s.,\\
&\bullet&  \mbox{ if $a>\alpha/5$, then } \sup_{x \in I}\left|\mathbb{E}\left(\tilde{T}_{n}\left(x\right)\right)\right|=o\left(\sqrt{\gamma_nh_n^{-1}}\,\ln n\right),\label{a1/08/b}
\end{eqnarray}
and by proving that, under the condition $\lim_{n\to\infty}\left(n\gamma_n\right)>2a/\varphi$,
\begin{eqnarray}
&\bullet&  \mbox{ if $a<\alpha/5$, then }\sup_{x\in I} \left|\tilde{T}_{n}^{\left(i\right)}\left( x\right)-\mathbb{E}\left(\tilde{T}_{n}^{\left(i\right)}\left(x\right)\right)\right|=o\left(h_n^2\right)\quad a.s.,\label{a1/08/c}\\
&\bullet & \mbox{ if $a\leq \alpha/5$, then } \sup_{x\in I} \mathbb{E}\left(\tilde{T}_{n}\left(x\right)\right)=O\left(h_n^2\right)\label{a1/08/d}.
\end{eqnarray}
As a matter of fact, Lemma~\ref{a4bis} follows from the combination of~\eqref{a1/08/a} and~\eqref{a1/08/b} in the case $a>\alpha/5$, from that of~\eqref{a1/08/a} and~\eqref{a1/08/d} in the case $a=\alpha/5$, and from that of~\eqref{a1/08/c} and~\eqref{a1/08/d} in the case $a<\alpha/5$. \\

The proof of~\eqref{a1/08/b} and~\eqref{a1/08/d} is similar to that of~\eqref{aTnE1} and~\eqref{aTnE}, and is omitted. Moreover the proof of~\eqref{a1/08/a} and~\eqref{a1/08/c} for $i=2$ is similar to that for $i=1$, and is omitted too. To prove simultaneously~\eqref{a1/08/a} and~\eqref{a1/08/c} for $i=1$, we introduce the sequence $\left(v_n\right)$ defined as
\begin{eqnarray}\label{avnf12}
\left(v_n\right)&=&\left\{\begin{array}{llll}\left(\sqrt{\gamma_n^{-1}h_n}\right) & \mbox{ if }  a \geq \alpha/5,\\
\left(h_n^{-2}\left(\ln n\right)^2\right) & \mbox{ if } a < \alpha/5 .
\end{array}\right.
\end{eqnarray}
As a matter of fact,~\eqref{a1/08/a} and~\eqref{a1/08/c} are proved for $i=1$ by establishing that
\begin{eqnarray}\label{a1/08/e}
\sup_{x\in I} \left|\tilde{T}_{n}^{\left(1\right)}\left( x\right)-\mathbb{E}\left(\tilde{T}_{n}^{\left(1\right)}\left(x\right)\right)\right|=O\left(v_n^{-1}\ln n\right)\quad a.s.
\end{eqnarray}
To this end we first state the following lemma.
\begin{lemma}\label{aL:ET} 
There exists $s>0$ such that, for all $C>0$,
\begin{eqnarray*}
\sup_{x\in I}\mathbb{P}\left[\frac{v_n}{\ln n}\left|\tilde{T}_n^{\left(1\right)}\left(x\right)-\mathbb{E}\left(\tilde{T}_n^{\left(1\right)}\left(x\right)\right)\right|\geq C\right]=O\left(n^{-\frac{C}{s}}\right).
\end{eqnarray*}
\end{lemma}
We first show how~\eqref{a1/08/e} can be deduced from Lemma~\ref{aL:ET}. Set $\left(M_n\right)\in \mathcal{GS}\left(\tilde{m}\right)$ with $\tilde{m}>0$, and note that, for all $C>0$, we have
\begin{eqnarray}\label{amot}
\lefteqn{\mathbb{P}\left[\frac{v_n}{\ln n}\sup_{x \in I}\left|\tilde{T}_n^{\left(1\right)}\left(x\right)-\mathbb{E}\left(\tilde{T}_n^{\left(1\right)}\left(x\right)\right)\right| \geq C\right]}\nonumber\\
& \leq & \mathbb{P}\left[\frac{v_n}{\ln n}\sup_{x \in I}\left|\tilde{T}_n^{\left(1\right)}\left(x\right)-\mathbb{E}\left(\tilde{T}_n^{\left(1\right)}\left(x\right)\right)\right| \geq C \mbox{ and } \sup_{k \leq n}\left|Y_k\right|\leq M_n\right]\nonumber\\
&&+\mathbb{P}\left[\sup_{k \leq n}\left|Y_k\right|\geq M_n\right].
\end{eqnarray}
Let $(d_n)$ be a positive sequence such that $d_n<1$ for all $n$, and such that 
$\lim_{n\to\infty}\gamma_n^{-1}d_n=0$. 
Let $I_i^{\left(n\right)}$ be $N\left(n\right)$ intervals of length $d_n$ such that $\cup_{i=1}^{N\left(n\right)}I_i^{\left(n\right)}=I$, and for all $i\in\left\{1,\ldots,N\left(n\right)\right\}$, set $x_i^{\left(n\right)} \in I_i^{\left(n\right)}$. We have
\begin{eqnarray*}
\lefteqn{\mathbb{P}\left[\frac{v_n}{\ln n}\sup_{x \in I}\left|\tilde{T}_n^{\left(1\right)}\left(x\right)-\mathbb{E}\left(\tilde{T}_n^{\left(1\right)}\left(x\right)\right)\right| \geq C \mbox{ and } \sup_{k \leq n}\left|Y_k\right|\leq M_n\right]}\\
&& \leq \sum_{i=1}^{N\left(n\right)}\mathbb{P}\left[\frac{v_n}{\ln n}\sup_{x \in I_i^{\left(n\right)}}\left|\tilde{T}_n^{\left(1\right)}\left(x\right)-\mathbb{E}\left(\tilde{T}_n^{\left(1\right)}\left(x\right)\right)\right|\geq C \mbox{ and } \sup_{k \leq n}\left|Y_k\right|\leq M_n\right].
\end{eqnarray*}
Let us prove that there exists $c^*$ such that, for all $x,y\in I$ such that $\left|x-y\right|\leq d_n$, and on $\left\{\sup_{k \leq n}\left|Y_k\right|\leq M_n\right\}$, 
\begin{eqnarray}\label{aTMS}
\left|\tilde{T}_n^{\left(1\right)}\left(x\right)-\tilde{T}_n^{\left(1\right)}\left(y\right)\right|\leq c^*M_nh_n^{-1}\gamma_n^{-1}d_n.
\end{eqnarray}
To this end, we write
\begin{eqnarray*}
\left|\tilde{T}_n^{\left(1\right)}\left(x\right)-\tilde{T}_n^{\left(1\right)}\left(y\right)\right|&\leq &A_{n,1}\left(x,y\right)+A_{n,2}\left(x,y\right)
\end{eqnarray*}
with
\begin{eqnarray*}
A_{n,1}\left(x,y\right)&=&\left|\sum_{k=n_0}^nU_{k,n}\left(f\left(x\right)\right)\gamma_kh_k^{-1}Y_k\left[K\left(\frac{x-X_k}{h_k}\right)-K\left(\frac{y-X_k}{h_k}\right)\right]\right|,\\
A_{n,2}\left(x,y\right)&=&\left|\sum_{k=n_0}^nU_{k,n}\left(f\left(y\right)\right)\gamma_kh_k^{-1}Y_kK\left(\frac{y-X_k}{h_k}\right)\left[\frac{U_{k,n}\left(f\left(x\right)\right)}{U_{k,n}\left(f\left(y\right)\right)}-1\right]\right|.
\end{eqnarray*}
$\bullet$ Since $K$ is Lipschitz-continuous, and by application of Lemma~\ref{alap}, there exist $k^*,c_1^*>0$ such that, for all $x,y \in \mathbb{R}$ satisfying $\left|x-y\right| \leq d_n$ and on $\left\{\sup_{k \leq n}\left|Y_k\right|\leq M_n\right\}$, we have:
\begin{eqnarray}\label{a57bis}
A_{n,1}\left(x,y\right)& \leq & k^*M_n\sum_{k=n_0}^nU_{k,n}\left(\varphi\right)\gamma_kh_k^{-2}d_n\nonumber\\
& \leq & c_1^*M_nh_n^{-2}d_n.
\end{eqnarray}
$\bullet$ Now, let $c_i^*$ be positive constants; for $j \geq n_0$ and for $p\in \left\{1,2\right\}$ we have  
\begin{eqnarray}\label{aju8}
\left(\frac{1-\gamma_jf\left(x\right)}{1-\gamma_jf\left(y\right)}\right)^p&=&\left(1+\frac{\gamma_j\left(f\left(y\right)-f\left(x\right)\right)}{1-\gamma_jf\left(y\right)}\right)^p\nonumber\\
& \leq & \left(1+\frac{c_2^*\gamma_jd_n}{1-\gamma_j\|f\|_{\infty}}\right)^p\nonumber\\
& \leq & \left(1+2c_2^*\gamma_jd_n\right)^p\nonumber\\
& \leq &1+c_3^*\gamma_jd_n.
\end{eqnarray}
We deduce that, for $k$ and $n$ such that $n_0\leq k\leq n$,
\begin{eqnarray}\label{a9ju}
\left[\frac{U_{k,n}^p\left(f\left(x\right)\right)}{U_{k,n}^p\left(f\left(y\right)\right)}-1\right]& = &\left( \prod_{j=k+1}^{n}\frac{1-\gamma_jf\left(x\right)}{1-\gamma_jf\left(y\right)}\right)^p-1\nonumber\\
& \leq & \left(\prod_{j=k+1}^{n}\exp\left(c_3^*\gamma_jd_n\right)\right)-1\nonumber\\
& \leq & \exp\left[c_3^*d_n\gamma_n^{-1}\sum_{j=k+1}^{n}\gamma_j\gamma_n\right]-1\nonumber\\
& \leq & \exp\left[c_4^*d_n\gamma_n^{-1}\sum_{n\geq 1}\gamma_n^2\right]-1\nonumber\\
& \leq & \exp\left[c_5^*d_n\gamma_n^{-1}\right]-1\nonumber\\
& \leq & c_5^*d_n\gamma_n^{-1}\exp\left[c_5^*d_n\gamma_n^{-1}\right]\nonumber\\
& \leq & c_6^*d_n\gamma_n^{-1}.
\end{eqnarray}
The application of~\eqref{a9ju} with $p=1$ and of Lemma~\ref{alap} ensures that, for all $x,y \in I$ satisfying $\left|x-y\right| \leq d_n$, and on $\left\{\sup_{k \leq n}\left|Y_k\right|\leq M_n\right\}$, we have:
\begin{eqnarray}\label{a59bis}
A_{n,2}\left(x,y\right)& \leq & \|K\|_{\infty}M_n\sum_{k=n_0}^nU_{k,n}\left(\varphi\right)\gamma_kh_k^{-1}\left(c_6^*d_n\gamma_n^{-1}\right)\nonumber\\
& \leq & c_7^*M_nh_n^{-1}\gamma_n^{-1}d_n.
\end{eqnarray}
The upper bound~\eqref{aTMS} follows from the combination of~\eqref{a57bis} and~\eqref{a59bis}.\\
Now, it follows from~\eqref{aTMS} that, for all $x \in I_i^{\left(n\right)}$ and on $\left\{\sup_{k \leq n}\left|Y_k\right|\leq M_n\right\}$, we have
\begin{eqnarray*}
\lefteqn{\left|\tilde{T}_n^{\left(1\right)}\left(x\right)-\mathbb{E}\left(\tilde{T}_n^{\left(1\right)}\left(x\right)\right)\right|}\\
 &\leq & \left|\tilde{T}_n^{\left(1\right)}\left(x\right)-\tilde{T}_n^{\left(1\right)}\left(x_i^{\left(n\right)}\right)\right|+\left|\tilde{T}_n^{\left(1\right)}\left(x_i^{\left(n\right)}\right)-\mathbb{E}\left(\tilde{T}_n^{\left(1\right)}\left(x_i^{\left(n\right)}\right)\right)\right|+\left|\mathbb{E}\left(\tilde{T}_n^{\left(1\right)}\left(x_i^{\left(n\right)}\right)\right)-\mathbb{E}\left(\tilde{T}_n^{\left(1\right)}\left(x\right)\right)\right| \\
&\leq & 2c^*M_nh_n^{-1}\gamma_n^{-1}d_n+\left|\tilde{T}_n^{\left(1\right)}\left(x_i^{\left(n\right)}\right)-\mathbb{E}\left(\tilde{T}_n^{\left(1\right)}\left(x_i^{\left(n\right)}\right)\right)\right|.
\end{eqnarray*}
In view of~\eqref{amot}, we obtain, for all $C>0$, 
\begin{eqnarray*}
\lefteqn{\mathbb{P}\left[\frac{v_n}{\ln n}\sup_{x \in I}\left|\tilde{T}_n^{\left(1\right)}\left(x\right)-\mathbb{E}\left(\tilde{T}_n^{\left(1\right)}\left(x\right)\right)\right| \geq C \right]}\nonumber\\
& \leq &\sum_{i=1}^{N\left(n\right)}\mathbb{P}\left[\frac{v_n}{\ln n}\left|T_n^{\left(1\right)}\left(x_i^{\left(n\right)}\right)-\mathbb{E}\left(T_n^{\left(1\right)}\left(x_i^{\left(n\right)}\right)\right)\right|+2c^*M_n\frac{v_n}{\ln n}h_n^{-1}\gamma_n^{-1}d_n \geq C \right]\nonumber\\
&&+n\mathbb{P}\left[\left|Y\right|\geq M_n\right].
\end{eqnarray*}
Now, note that, in view of~\eqref{avnf12}, $\left(v_n\right)\in\mathcal{GS}\left(m^*\right)$ where $m^*$ is defined in~\eqref{a10bis}. Set $\left(d_n\right)\in \mathcal{GS}\left(-\left(\tilde{m}+m^*+a+\alpha\right)\right)$ such that, for all $n$, $2c^*M_nv_n\left(\ln n\right)^{-1}h_n^{-1}\gamma_n^{-1}d_n \leq C/2$; in view of Lemma~\ref{aL:ET} and Assumption $\left(A4\right)ii)$, there exists $s>0$ such that
\begin{eqnarray*}
\lefteqn{\mathbb{P}\left[\frac{v_n}{\ln n}\sup_{x \in I}\left|\tilde{T}_n^{\left(1\right)}\left(x\right)-\mathbb{E}\left(\tilde{T}_n^{\left(1\right)}\left(x\right)\right)\right| \geq C \right]}\nonumber\\
& \leq & N\left(n\right)\sup_{x \in I}\mathbb{P}\left[\frac{v_n}{\ln n}\left|\tilde{T}_n^{\left(1\right)}\left(x\right)-\mathbb{E}\left(\tilde{T}_n^{\left(1\right)}\left(x\right)\right)\right| \geq \frac{C}{2}\right]+n\exp \left(-t^*M_n\right)\mathbb{E}\left(\exp\left(t^*\left|Y\right|\right)\right)\nonumber\\
&=&O\left(d_n^{-1}n^{-\frac{C}{2s}}+n\exp\left(-t^*M_n\right)\right).
\end{eqnarray*}
Since $\left(M_n\right)\in \mathcal{GS}\left(\tilde{m}\right)$ with $\tilde{m}>0$, and since $\left(d_n^{-1}\right)\in \mathcal{GS}\left(\tilde{m}+m^*+a+\alpha\right)$, we can choose $C$ large enough so that 
\begin{eqnarray*}
\sum_{n \geq 0}\mathbb{P}\left[\frac{v_n}{\ln n}\sup_{x \in I}\left|\tilde{T}_n^{\left(1\right)}\left(x\right)-\mathbb{E}\left(\tilde{T}_n^{\left(1\right)}\left(x\right)\right)\right| \geq C \right]<\infty,
\end{eqnarray*}
which gives~\eqref{a1/08/e}.\\
It remains to prove Lemma~\ref{aL:ET}. For all $x \in I$ and all $s>0$, we have 
\begin{eqnarray}\label{a29mai}
\mathbb{P}\left[\frac{v_n}{\ln n}\left(\tilde{T}_n^{\left(1\right)}\left(x\right)-\mathbb{E}\left(\tilde{T}_n^{\left(1\right)}\left(x\right)\right)\right) \geq C\right]
&=&\mathbb{P}\left[\exp \left[s^{-1}v_n\left(\tilde{T}_n^{\left(1\right)}\left(x\right)-\mathbb{E}\left(\tilde{T}_n^{\left(1\right)}\left(x\right)\right)\right)\right] \geq n^{\frac{C}{s}}\right]\nonumber\\
&\leq& n^{-\frac{C}{s}}\mathbb{E}\left(\exp \left[s^{-1}v_n\left(\tilde{T}_n^{\left(1\right)}\left(x\right)-\mathbb{E}\left(\tilde{T}_n^{\left(1\right)}\left(x\right)\right)\right)\right]\right)\nonumber\\
&\leq& n^{-\frac{C}{s}}\prod_{k=n_0}^n\mathbb{E}\left(\exp \left(s^{-1}V_{k,n}\left(x\right)\right)\right)
\end{eqnarray}
with
\begin{eqnarray*}
V_{k,n}\left(x\right)&=&v_nU_{k,n}\left(f\left(x\right)\right)\gamma_kh_k^{-1}\left[Y_kK\left(\frac{x-X_k}{h_k}\right)-\mathbb{E}\left(Y_kK\left(\frac{x-X_k}{h_k}\right)\right)\right].
\end{eqnarray*}
For $k$ and $n$ such that $n_0 \leq k\leq n$, set 
\begin{eqnarray*}
\alpha_{k,n}&=&v_nU_{k,n}\left(\varphi\right)\gamma_kh_k^{-1}.
\end{eqnarray*}
We have, for all $x\in I$,
\begin{eqnarray*}
\lefteqn{\mathbb{E}\left(\exp \left[s^{-1}V_{k,n}\left(x\right)\right]\right)}\\
& \leq & 1+\frac{1}{2}\mathbb{E}\left[s^{-2}V_{k,n}^2\left(x\right)\right]+ \mathbb{E}\left(s^{-3}\left|V_{k,n}^3\left(x\right)\right|\right)\exp \left[\left|V_{k,n}\left(x\right)\right|\right]\\
& \leq & 1+\frac{1}{2}s^{-2}\alpha_{k,n}^2Var\left[Y_kK\left(\frac{x-X_k}{h_k}\right)\right]\\
&&+ s^{-3}\alpha_{k,n}^3\|K\|_{\infty}^3\mathbb{E}\left[\left(\left|Y_k\right|^3+\left(\mathbb{E}\left(\left|Y_k\right|\right)\right)^3\right)\exp\left(s^{-1}\alpha_{k,n}\|K\|_{\infty}\left(\left|Y_k\right|+\mathbb{E}\left(\left|Y_k\right|\right)\right)\right)\right].
\end{eqnarray*}
Now, note that $\alpha_{k,n}$ can be rewritten as:
\begin{eqnarray*}
\alpha_{k,n}=\frac{v_n\Pi_n\left(\varphi\right)}{v_k\Pi_k\left(\varphi\right)}v_k\gamma_kh_k^{-1}.
\end{eqnarray*}
Since $\left(v_n\right)\in\mathcal{GS}\left(m^*\right)$ with $\varphi-m^*\xi>0$ (where $\xi$ is defined in~\eqref{axi1p}), we have
\begin{eqnarray*}
\frac{\Pi_n\left(\varphi\right)}{\Pi_{n-1}\left(\varphi\right)}\frac{v_n}{v_{n-1}}
&=&\left(1-\gamma_n\varphi\right)\left(1+m^*\frac{1}{n}+o\left(\frac{1}{n}\right)\right)\\
&=&\left(1-\gamma_n\varphi\right)\left(1+m^*\xi \gamma_n+o\left(\gamma_n\right)\right)\\
&=&1-\left(\varphi-m^*\xi\right)\gamma_n+o\left(\gamma_n\right)\\
&\leq& 1\mbox{\quad\quad for $n$ large enough}.
\end{eqnarray*}
Writing
\begin{eqnarray*}
\frac{v_n\Pi_n\left(\varphi\right)}{v_k\Pi_k\left(\varphi\right)}=\prod_{i=k}^{n-1}\frac{v_{i+1}\Pi_{i+1}\left(\varphi\right)}{v_{i}\Pi_{i}\left(\varphi\right)},
\end{eqnarray*}
we obtain 
\begin{eqnarray*}
\sup_{n_0\leq k\leq n}\frac{v_n\Pi_n\left(\varphi\right)}{v_k\Pi_k\left(\varphi\right)}<\infty ,
\end{eqnarray*}
and, since $\lim_{n\to \infty}v_k\gamma_kh_k=0$, we deduce that $\sup_{n_0\leq k \leq n}\alpha_{k,n}<\infty$. Thus, in view of Assumption $\left(A4\right)ii)$, there exist $s>0$ and 
$c^*>0$ such that, for all $k$ and $n$ such that $n_0\leq k \leq n$,
\begin{eqnarray*}
\mathbb{E}\left[\left(\left|Y_k\right|^3+\left(\mathbb{E}\left(\left|Y_k\right|\right)\right)^3\right)\exp\left(s^{-1}\alpha_{k,n}\|K\|_{\infty}\left(\left|Y_k\right|+\mathbb{E}\left(\left|Y_k\right|\right)\right)\right)\right] \leq c^*.
\end{eqnarray*}
From classical computations, we have $\sup_{x \in I}Var\left[Y_kK\left(\left(x-X_k\right)h_k^{-1}\right)\right]=O\left(h_k\right)$. We then deduce that there exist $C_1^*,C_2^*>0$ such that, for all $x \in I$, for all $k$ and $n$ such that $n_0\leq k \leq n$,
\begin{eqnarray*}
\mathbb{E}\left(\exp \left[s^{-1}V_{k,n}\left(x\right)\right]\right)
& \leq & 1+C_1^*v_n^2U_{n,k}^{2}\left(\varphi\right)\gamma_k^2h_k^{-1}+C_2^*v_n^{3}U_{k,n}^{3}\left(\varphi\right)\gamma_k^3h_k^{-3}\\
& \leq & \exp \left[C_1^*v_n^2U_{k,n}^{2}\left(\varphi\right)\gamma_k^2h_k^{-1}+C_2^*v_n^3U_{k,n}^{3}\left(\varphi\right)\gamma_k^3h_k^{-3}\right].
\end{eqnarray*}
Applying Lemma~\ref{alap}, we deduce from~\eqref{a29mai} that, for all $C>0$,
\begin{eqnarray*}
\lefteqn{\sup_{x \in I}\mathbb{P}\left[\frac{v_n}{\ln n}\left(\tilde{T}_n^{\left(1\right)}\left(x\right)-\mathbb{E}\left(\tilde{T}_n^{\left(1\right)}\left(x\right)\right)\right) \geq C\right]}\\
&\leq & n^{-\frac{C}{s}}\exp \left[C_1^*v_n^2\sum_{k=n_0}^nU_{k,n}^{2}\left(\varphi\right)\gamma_kO\left(v_k^{-2}\right)+C_2^*v_n^3\sum_{k=n_0}^nU_{k,n}^{3}\left(\varphi\right)\gamma_kO\left(v_k^{-3}\right)\right]\\
& = & O\left(n^{-\frac{C}{s}}\right).
\end{eqnarray*}
We establish exactly in the same way that, for all $C>0$,
\begin{eqnarray*}
\sup_{x \in I}\mathbb{P}\left[\frac{v_n}{\ln n}\left(\mathbb{E}\left(\tilde{T}_n^{\left(1\right)}\left(x\right)\right)-\tilde{T}_n^{\left(1\right)}\left(x\right)\right) \geq C\right] &=&O\left(n^{-\frac{C}{s}}\right),
\end{eqnarray*}
which concludes the proof of Lemma~\ref{aL:ET}.
\subsection{Proof of Lemma~\ref{aL:rrb}}\label{a4.2}
Set 
\begin{eqnarray}\label{ayynk}
\eta_{k}\left(x\right)=\left(W_{k+1}\left(x\right)-r\left(x\right)Z_{k+1}\left(x\right)\right).
\end{eqnarray}
In order to prove Lemma \ref{aL:rrb}, we first establish a central limit theorem for
\begin{eqnarray*}
T_n\left(x\right)-\mathbb{E}\left(T_n\left(x\right)\right)=\frac{1}{\sum_{k=1}^nq_k}
\sum_{k=n_0-1}^nq_k\left[\eta_k\left(x\right)-\mathbb{E}\left(\eta_k\left(x\right)\right)\right].
\end{eqnarray*}
In view of~\eqref{avarw}-\eqref{acovwz} and since $h_k/h_{k+1}=1+o\left(1\right)$, we have
\begin{eqnarray*}
Var\left(\eta_{k}\left(x\right)\right)&=&Var\left(W_{k+1}\left(x\right)\right)+r^2\left(x\right)Var\left(Z_{k+1}\left(x\right)\right)-2r\left(x\right)Cov\left(W_{k+1}\left(x\right),Z_{k+1}\left(x\right)\right),\\
&=&\frac{1}{h_{k}}\left[Var\left[Y\vert X=x\right]f\left(x\right)\int_{\mathbb{R}} K^2\left(z\right)dz+o\left(1\right)\right].
\end{eqnarray*}
Noting that $\left(q_n^2h_n^{-1}\right)\in\mathcal{GS}\left(-2q+a\right)$ with $q<\left(1+a\right)/2$, and using~\eqref{abt}, we get
\begin{eqnarray}\label{avnmoy} 
v_n^2&=&\sum_{k=n_0-1}^nq_k^2Var\left(\eta_{k}\left(x\right)\right)\nonumber\\
&=&\sum_{k=n_0-1}^n\frac{q_k^2}{h_{k}}\left[Var\left[Y\vert X=x\right]f\left(x\right)\int_{\mathbb{R}}K^2\left(z\right)dz+o\left(1\right)\right]\nonumber\\
&=&\frac{nq_n^2h_n^{-1}}{1-2q+a}\left[Var\left[Y\vert X=x\right]f\left(x\right)\int_{\mathbb{R}}K^2\left(z\right)dz+o\left(1\right)\right].
\end{eqnarray}
Now, set $p\in\left]0,1\right]$ such that $q<\left(1+a\left(1+p\right)\right)/\left(2+p\right)$; it follows from~\eqref{aMok} that
\begin{eqnarray*}
\sum_{k=n_0-1}^nq_k^{2+p}\mathbb{E}\left[\left|\eta_{k}\left(x\right)\right|^{2+p}\right]
& = & O\left(\sum_{k=n_0-1}^n \frac{q_k^{2+p}}{h_k^{2+p}}\mathbb{E}\left(\left|Y_k-r\left(x\right)\right|^{2+p}K^{2+p}\left(\frac{x-X_k}{h_k}\right)\right)\right)\nonumber\\
& = & O\left(\sum_{k=n_0-1}^n \frac{q_k^{2+p}}{h_k^{1+p}}\right).
\end{eqnarray*}
In view of~\eqref{avnmoy} and using~\eqref{abt}, we get
\begin{eqnarray*}
\frac{1}{v_n^{2+p}}\sum_{k=n_0-1}^nq_k^{2+p}\mathbb{E}\left[\left|\eta_{k}\left(x\right)\right|^{2+p}\right]&=&O\left(\frac{n q_n^{2+p}h_{n}^{-\left(1+p\right)}}{\left(n q_n^2h_{n}^{-1}\right)^{1+\frac{p}{2}}}\right)\\
&=&O\left(\frac{1}{n^{\frac{p}{2}}h_n^{\frac{p}{2}}}\right)\\
&=&o\left(1\right).
\end{eqnarray*}
The application of Lyapounov Theorem gives 
\begin{eqnarray*}
\frac{\sum_{k=1}^nq_k}{\sqrt{nq_n^2h_n^{-1}}}\left(T_{n}\left( x\right)-\mathbb{E}\left[T_n\left(x\right)\right]\right) \stackrel{\mathcal{D}}{\rightarrow}\mathcal{N}\left( 0,\frac{1}{1+a-2q}Var\left[Y\vert X=x\right]f\left(x\right)\int_{\mathbb{R}} K^2\left(z\right)dz\right),
\end{eqnarray*}
and applying~\eqref{abt}, we obtain
\begin{eqnarray}\label{arevtlc}
\sqrt{nh_n}\left(T_{n}\left( x\right)-\mathbb{E}\left[T_n\left(x\right)\right]\right) \stackrel{\mathcal{D}}{\rightarrow}\mathcal{N}\left( 0,\frac{\left(1-q\right)^2}{1+a-2q} Var\left[Y\vert X=x\right]f\left(x\right)\int_{\mathbb{R}} K^2\left(z\right)dz\right).
\end{eqnarray}
Now, note that
\begin{eqnarray*}
\mathbb{E}\left(T_n\left(x\right)\right)&=&\frac{1}{\sum_{k=1}^nq_k}\sum_{k=n_0-1}^{n}q_k\left[\left(\mathbb{E}\left(W_{k+1}\left(x\right)\right)-a\left(x\right)\right)-r\left(x\right)\left(\mathbb{E}\left(Z_{k+1}\left(x\right)\right)-f\left(x\right)\right)\right].
\end{eqnarray*}
Since $h_{n+1}/h_n=1+o\left(1\right)$, it follows from~\eqref{aespw} and~\eqref{aespz} that
\begin{eqnarray*}
\lefteqn{\mathbb{E}\left(T_n\left(x\right)\right)}\nonumber\\
&=&\frac{1}{\sum_{k=1}^nq_k}\sum_{k=n_0-1}^{n}q_kh_{k+1}^2\left[\frac{1}{2}\left(\int_{\mathbb{R}}y\frac{\partial^2 g}{\partial x^2}\left(x,y\right)dy-r\left(x\right)\int_{\mathbb{R}}\frac{\partial^2 g}{\partial x^2}\left(x,y\right)dy\right)+o\left(1\right)\right]\int_{\mathbb{R}}z^2K\left(z\right)dz\nonumber\\
&=&\frac{1}{\sum_{k=1}^nq_k}\sum_{k=n_0-1}^{n}q_kh_{k}^2\left[m^{\left(2\right)}\left(x\right)f\left(x\right)+o\left(1\right)\right].
\end{eqnarray*}
Applying~\eqref{abt}, we obtain
\begin{eqnarray}\label{aTna*31}
\lim_{n \to \infty}\frac{1}{h_n^2}\mathbb{E}\left(T_n\left(x\right)\right)
&=&\frac{1-q}{1-2a-q}m^{\left(2\right)}\left(x\right)f\left(x\right), 
\end{eqnarray}
and Lemma~\ref{aL:rrb} follows from the combination of~\eqref{arevtlc} and~\eqref{aTna*31}.
\subsection{Proof of Lemma~\ref{aL:rrc}}\label{a4.3}
Set
\begin{eqnarray*}
S_n\left(x\right)&=&\sum_{k=n_0-1}^nq_k\left[\eta_{k}\left(x\right)-\mathbb{E}\left(\eta_{k}\left(x\right)\right)\right]
\end{eqnarray*}
where $\eta_k$ is defined in~\eqref{ayynk}, and $H_n^{-2}=nh_n^{-1}q_n^2$. Let us first note that, since $\left(nh_n^{-1}q_n^2\right)\in \mathcal{GS}\left(1+a-2q\right)$ with $1+a-2q>0$, we have $\lim_{n\to\infty}H_n^{-2}=\infty$. Moreover, we have $\lim_{n\to\infty}H_n^2/H_{n-1}^2=1$, and, by~\eqref{avnmoy}, 
\[ 
\lim_{n\to\infty}H_n^2\sum_{k=n_0-1}^nq_k^2Var\left[\eta_k\left(x\right)\right]=\left[1+a-2q\right]^{-1}Var\left[Y\vert X=x\right]f\left(x\right)\int_{\mathbb{R}}K^2\left(z\right)dz
\]
and, by~\eqref{aMok}, $\mathbb{E}\left[\left|q_n\eta_{n}\left(x\right)\right|^{3}\right]=O\left(q_n^3h_n^{-2}\right)$. Since, for all $\epsilon >0$,
\begin{eqnarray*}
\frac{1}{n\sqrt{n}}\sum_{k=n_0-1}^n\mathbb{E}\left(\left|H_nq_k\eta_k\left(x\right)\right|^3\right)&=&O\left(\frac{H_n^3}{n\sqrt{n}}\sum_{k=n_0-1}^n\frac{q_k^3}{h_k^2}\right)\\
&=&O\left(n^{-3}h_n^{\frac{3}{2}}q_n^{-3}\left(n^{\epsilon}+nq_n^3h_n^{-2}\right)\right),
\end{eqnarray*}
we have
\begin{eqnarray*}
\frac{1}{n\sqrt{n}}
\sum_{k=n_0-1}^n\mathbb{E}\left(\left|H_nq_k\eta_k\left(x\right)\right|^3\right)
=o\left(\left[\ln \left(H_n^{-2}\right)\right]^{-1}\right).
\end{eqnarray*}
The application of Theorem 1 of Mokkadem and Pelletier (2008) ensures that, with probability one, the sequence
\begin{eqnarray*}
\left(\frac{H_nS_n\left(x\right)}{\sqrt{2\ln \ln \left(H_n^{-2}\right)}}\right)=\left(\frac{\sum_{k=1}^nq_k}{nq_n}\frac{\sqrt{nh_n}\left(T_{n}\left( x\right)-\mathbb{E}\left(T_n\left(x\right)\right)\right)}{\sqrt{2\ln \ln \left(H_n^{-2}\right)}}\right)
\end{eqnarray*}
is relatively compact and its limit set is the interval 
\begin{eqnarray*}
\left[-\sqrt{\frac{1}{1+a-2q}Var\left[Y\vert X=x\right]f\left(x\right)\int_{\mathbb{R}} K^2\left(z\right)dz},\sqrt{\frac{1}{1+a-2q}Var\left[Y\vert X=x\right]f\left(x\right)\int_{\mathbb{R}} K^2\left(z\right)dz}\right].
\end{eqnarray*}
Since $\lim_{n\to\infty}\ln \ln \left(H_n^{-2}\right)/\ln \ln n=1$, and using~\eqref{abt}, it follows that, with probability one, the sequence $\left(\sqrt{nh_n}\left(T_{n}\left( x\right)-\mathbb{E}\left(T_n\left(x\right)\right)\right)/\sqrt{2\ln \ln n}\right)$ is relatively compact, and its limit set is the interval
\begin{eqnarray*}
\left[-\sqrt{\frac{\left(1-q\right)^2}{1+a-2q}Var\left[Y\vert X=x\right]f\left(x\right)\int_{\mathbb{R}} K^2\left(z\right)dz},\sqrt{\frac{\left(1-q\right)^2}{1+a-2q}Var\left[Y\vert X=x\right]f\left(x\right)\int_{\mathbb{R}} K^2\left(z\right)dz}\right].
\end{eqnarray*}
The application of~\eqref{aTna*31} concludes the proof of Lemma~\ref{aL:rrc}.

\subsection{Proof of Lemma~\ref{aL:Turbar}}\label{a4.4} 
 Let us write $T_n\left(x\right)$ as
\begin{eqnarray*}
T_n\left(x\right)=T_{n,1}\left(x\right)-r\left(x\right)T_{n,2}\left(x\right)
\end{eqnarray*}
with
\begin{eqnarray*}
T_{n,1}\left(x\right)&=&\frac{1}{\sum_{k=1}^nq_k}\sum_{k=n_0-1}^n\frac{q_k}{h_{k+1}}Y_{k+1}K\left(\frac{x-X_{k+1}}{h_{k+1}}\right)\\
T_{n,2}\left(x\right)&=&\frac{1}{\sum_{k=1}^nq_k}\sum_{k=n_0-1}^n\frac{q_k}{h_{k+1}}K\left(\frac{x-X_{k+1}}{h_{k+1}}\right).
\end{eqnarray*}
Lemma~\ref{aL:Turbar} is proved by showing that, for $i\in \left\{1,2\right\}$,
\begin{eqnarray}\label{aT12bar}
\sup_{x\in I}\left|T_{n,i}\left(x\right)-\mathbb{E}\left(T_{n,i}\left(x\right)\right)\right|=O\left(\sqrt{n^{-1}h_n^{-1}}\ln n\right)\quad a.s.,
\end{eqnarray}
and that
\begin{eqnarray}\label{aET12bar}
\sup_{x\in I}\left|\mathbb{E}\left(T_{n,1}\left(x\right)\right)-r\left(x\right)\mathbb{E}\left(T_{n,2}\left(x\right)\right)\right|=O\left(h_n^{2}\right).
\end{eqnarray}
The proof of~\eqref{aET12bar} relies on classical computations and is omitted. Moreover the proof of~\eqref{aT12bar} for $i=2$ is similar to that for $i=1$, and is omitted too. We now prove~\eqref{aT12bar} for $i=1$. To this end, we first state the following lemma.
\begin{lemma}\label{aL:ETbar} 
There exists $s>0$ such that, for all $C>0$,
\begin{eqnarray*}
\sup_{x\in I}\mathbb{P}\left[\frac{\sqrt{nh_n}}{\ln n}\left|T_{n,1}\left(x\right)-\mathbb{E}\left(T_{n,1}\left(x\right)\right)\right|\geq C\right]=O\left(n^{-\frac{C}{s}}\right).
\end{eqnarray*}
\end{lemma}
We first show how~\eqref{aT12bar} for $i=1$ can be deduced from Lemma~\ref{aL:ETbar}, and then prove Lemma~\ref{aL:ETbar}. Set $\left(M_n\right)\in \mathcal{GS}\left(\tilde m\right)$ with 
$\tilde m>0$, and note that, for all $C>0$, we have
\begin{eqnarray}
\lefteqn{\mathbb{P}\left[\frac{\sqrt{nh_n}}{\ln n}\sup_{x \in I}\left|T_{n,1}\left(x\right)-\mathbb{E}\left(T_{n,1}\left(x\right)\right)\right| \geq C\right]}\nonumber\\
& \leq & \mathbb{P}\left[\frac{\sqrt{nh_n}}{\ln n}\sup_{x \in I}\left|T_{n,1}\left(x\right)-\mathbb{E}\left(T_{n,1}\left(x\right)\right)\right| \geq C \mbox{ and } \sup_{k \leq n}\left|Y_{k+1}\right|\leq M_n\right]\nonumber\\
&&+\mathbb{P}\left[\sup_{k \leq n}\left|Y_{k+1}\right|\geq M_n\right]. \nonumber
\end{eqnarray}
Let $I_i^{\left(n\right)}$ be $N\left(n\right)$ intervals of length $d_n$ such that $\cup_{i=1}^{N\left(n\right)}I_i^{\left(n\right)}=I$, and for all $i\in\left\{1,\ldots,N\left(n\right)\right\}$, set $x_i^{\left(n\right)} \in I_i^{\left(n\right)}$. We have
\begin{eqnarray}
\lefteqn{\mathbb{P}\left[\frac{\sqrt{nh_n}}{\ln n}\sup_{x \in I}\left|T_{n,1}\left(x\right)-\mathbb{E}\left(T_{n,1}\left(x\right)\right)\right| \geq C\right]}\nonumber\\
&\leq&  \sum_{i=1}^{N\left(n\right)}\mathbb{P}\left[\frac{\sqrt{nh_n}}{\ln n}\sup_{x \in I_i^{\left(n\right)}}\left|T_{n,1}\left(x\right)-\mathbb{E}\left(T_{n,1}\left(x\right)\right)\right|\geq C \mbox{ and } \sup_{k \leq n}\left|Y_{k+1}\right|\leq M_n\right]\nonumber\\
&& \mbox{ }+\mathbb{P}\left[\sup_{k \leq n}\left|Y_{k+1}\right|\geq M_n\right]. \label{amotbar}
\end{eqnarray}
Since $K$ is Lipschitz-continuous, there exist $k^*,c^*>0$, such that, for all $x,y \in \mathbb{R}$ satisfying $\left|x-y\right| \leq d_n$ and on $\left\{\sup_{k \leq n}\left|Y_{k+1}\right|\leq M_n\right\}$, we have:
\begin{eqnarray*}
\left|T_{n,1}\left(x\right)-T_{n,1}\left(y\right)\right|&=&\left|\frac{1}{\sum_{k=1}^{n}q_k}\sum_{k=n_0-1}^nq_kh_{k+1}^{-1}Y_{k+1}\left[K\left(\frac{x-X_{k+1}}{h_{k+1}}\right)-K\left(\frac{y-X_{k+1}}{h_{k+1}}\right)\right]\right|\nonumber\\
& \leq & k^*M_nd_n\frac{1}{\sum_{k=1}^{n}q_k}\sum_{k=1}^nq_kh_{k+1}^{-2}\nonumber\\
& \leq & c^*M_nh_n^{-2}d_n.
\end{eqnarray*}
It follows that, for all $x \in I_i^{\left(n\right)}$, on $\left\{\sup_{k \leq n}\left|Y_{k+1}\right|\leq M_n\right\}$, we have
\begin{eqnarray*}
\lefteqn{\left|T_{n,1}\left(x\right)-\mathbb{E}\left(T_{n,1}\left(x\right)\right)\right|}\\
 &\leq & \left|T_{n,1}\left(x\right)-T_{n,1}\left(x_i^{\left(n\right)}\right)\right|+\left|T_{n,1}\left(x_i^{\left(n\right)}\right)-\mathbb{E}\left(T_{n,1}\left(x_i^{\left(n\right)}\right)\right)\right|+\left|\mathbb{E}\left(T_{n,1}\left(x_i^{\left(n\right)}\right)\right)-\mathbb{E}\left(T_{n,1}\left(x\right)\right)\right| \\
&\leq & 2c^*M_nh_n^{-2}d_n+\left|T_{n,1}\left(x_i^{\left(n\right)}\right)-\mathbb{E}\left(T_{n,1}\left(x_i^{\left(n\right)}\right)\right)\right|.
\end{eqnarray*}
In view of~\eqref{amotbar}, we obtain, for all $C>0$, 
\begin{eqnarray*}
\lefteqn{\mathbb{P}\left[\frac{\sqrt{nh_n}}{\ln n}\sup_{x \in I}\left|T_{n,1}\left(x\right)-\mathbb{E}\left(T_{n,1}\left(x\right)\right)\right| \geq C \right]}\nonumber\\
& \leq &\sum_{i=1}^{N\left(n\right)}\mathbb{P}\left[\frac{\sqrt{nh_n}}{\ln n}\left|T_{n,1}\left(x_i^{\left(n\right)}\right)-\mathbb{E}\left(T_{n,1}\left(x_i^{\left(n\right)}\right)\right)\right|+2c^*M_n\frac{\sqrt{nh_n}}{\ln n}h_n^{-2}d_n \geq C \right]\nonumber\\
&&+n\mathbb{P}\left[\left|Y\right|\geq M_n\right].
\end{eqnarray*}
Now, set $\left(d_n\right)\in \mathcal{GS}\left(-\frac{1}{2}-\frac{3}{2}a-\tilde m\right)$ such that, for all $n$, $2c^*M_n\sqrt{nh_n}\left(\ln n\right)^{-1}h_n^{-2}d_n \leq C/2$; in view of Lemma~\ref{aL:ETbar} and Assumption $\left(A4\right)ii)$, there exists $s>0$ such that
\begin{eqnarray*}
\lefteqn{\mathbb{P}\left[\frac{\sqrt{nh_n}}{\ln n}\sup_{x \in I}\left|T_{n,1}\left(x\right)-\mathbb{E}\left(T_{n,1}\left(x\right)\right)\right| \geq C \right]}\nonumber\\
& \leq & N\left(n\right)\sup_{x \in I}\mathbb{P}\left[\frac{\sqrt{nh_n}}{\ln n}\left|T_{n,1}\left(x\right)-\mathbb{E}\left(T_{n,1}\left(x\right)\right)\right| \geq \frac{C}{2}\right]+n\exp \left(-t^*M_n\right)\mathbb{E}\left(\exp\left(t^*\left|Y\right|\right)\right)\nonumber\\
&=&O\left(d_n^{-1}n^{-\frac{C}{2s}}+n\exp\left(-t^*M_n\right)\right).
\end{eqnarray*}
Since $\left(d_n^{-1}\right)\in \mathcal{GS}\left(\frac{1}{2}+\frac{3}{2}a+\tilde m\right)$ and since $\left(M_n\right)\in \mathcal{GS}\left(\tilde m\right)$ with $\tilde m>0$, we can choose $C$ large enough so that 
\begin{eqnarray*}
\sum_{n \geq 0}\mathbb{P}\left[\frac{\sqrt{nh_n}}{\ln n}\sup_{x \in I}\left|T_{n,1}\left(x\right)-\mathbb{E}\left(T_{n,1}\left(x\right)\right)\right| \geq C \right]<\infty,
\end{eqnarray*}
which gives~\eqref{aT12bar} for $i=1$.\\
It remains to prove Lemma~\ref{aL:ETbar}. For all $x \in I$ and all $s>0$, we have 
\begin{eqnarray}\label{a29maibar}
\mathbb{P}\left[\frac{\sqrt{nh_n}}{\ln n}\left(T_{n,1}\left(x\right)-\mathbb{E}\left(T_{n,1}\left(x\right)\right)\right) \geq C\right]
&=&\mathbb{P}\left[\exp \left[s^{-1}\sqrt{nh_n}\left(T_{n,1}\left(x\right)-\mathbb{E}\left(T_{n,1}\left(x\right)\right)\right)\right] \geq n^{\frac{C}{s}}\right]\nonumber\\
&\leq& n^{-\frac{C}{s}}\mathbb{E}\left(\exp \left[s^{-1}\sqrt{nh_n}\left(T_{n,1}\left(x\right)-\mathbb{E}\left(T_{n,1}\left(x\right)\right)\right)\right]\right)\nonumber\\
&\leq& n^{-\frac{C}{s}}\prod_{k=n_0-1}^n\mathbb{E}\left(\exp \left(s^{-1}U_{k,n}\left(x\right)\right)\right)
\end{eqnarray}
with
\begin{eqnarray*}
U_{k,n}\left(x\right)&=&\frac{\sqrt{nh_n}}{\sum_{k=1}^nq_k}\frac{q_k}{h_{k+1}}\left[Y_{k+1}K\left(\frac{x-X_{k+1}}{h_{k+1}}\right)-\mathbb{E}\left(\left(Y_{k+1}K\left(\frac{x-X_{k+1}}{h_{k+1}}\right)\right)\right)\right].
\end{eqnarray*}
For $k$ and $n$ such that $k\leq n$, set 
\begin{eqnarray*}
\alpha_{k,n}&=&\frac{\sqrt{nh_n}}{\sum_{k=1}^nq_k}\frac{q_k}{h_{k+1}}.
\end{eqnarray*}
We have, for all $x\in I$,
\begin{eqnarray*}
\lefteqn{\mathbb{E}\left(\exp \left[s^{-1}U_{k,n}\left(x\right)\right]\right)}\\
& \leq & 1+\frac{1}{2}\mathbb{E}\left[s^{-2}U_{k,n}^2\left(x\right)\right]+ \mathbb{E}\left[s^{-3}\left|U_{k,n}^3\left(x\right)\right|\right]\exp \left[\left|U_{k,n}\left(x\right)\right|\right]\\
& \leq & 1+\frac{1}{2}s^{-2}\alpha_{k,n}^2Var\left[Y_{k+1}K\left(\frac{x-X_{k+1}}{h_{k+1}}\right)\right]\\
&&+ s^{-3}\alpha_{k,n}^3\|K\|_{\infty}^3\mathbb{E}\left[\left(\left|Y_{k+1}\right|^3+\left(\mathbb{E}\left(\left|Y_{k+1}\right|\right)\right)^3\right)\exp\left(s^{-1}\alpha_{k,n}\|K\|_{\infty}\left(\left|Y_{k+1}\right|+\mathbb{E}\left(\left|Y_{k+1}\right|\right)\right)\right)\right].
\end{eqnarray*}
Now, note that
\begin{eqnarray*}
\alpha_{k,n}&=&\left(\frac{nq_n}{\sum_{k=1}^nq_k}\right)\frac{q_k\sqrt{kh_k^{-1}}}{q_n\sqrt{nh_n^{-1}}}\sqrt{k^{-1}h_kh_{k+1}^{-2}}\nonumber\\
&=&\left(\frac{nq_n}{\sum_{k=1}^nq_k}\right)\left(\prod_{j=k}^{n-1}\frac{q_j\sqrt{jh_j^{-1}}}{q_{j+1}\sqrt{\left(j+1\right)h_{j+1}^{-1}}}\right)\sqrt{k^{-1}h_kh_{k+1}^{-2}}.
\end{eqnarray*}
Since $\left(q_j\sqrt{jh_j^{-1}}\right)\in\mathcal{GS}\left(-q+\left(1+a\right)/2\right)$ 
with $-q+\left(1+a\right)/2>0$, we have
\begin{eqnarray*}
\frac{q_j\sqrt{jh_j^{-1}}}{q_{j+1}\sqrt{\left(j+1\right)h_{j+1}^{-1}}}&=&1-\left(-q+\frac{1+a}{2}\right)\frac{1}{j}+o\left(\frac{1}{j}\right)\\
& \leq &1 \quad \mbox{for $j$ large enough.}
\end{eqnarray*}
It follows that $\sup_{k\leq n}\alpha_{k,n}<\infty$. Consequently, in view of Assumption $\left(A4\right)ii)$, there exist $s>0$ and $c^*>0$ such that, for all $k$ and $n$ such that 
$k\leq n$,
\begin{eqnarray*}
\mathbb{E}\left[\left(\left|Y_{k+1}\right|^3+\left(\mathbb{E}\left(\left|Y_{k+1}\right|\right)\right)^3\right)\exp\left(s^{-1}\alpha_{k,n}\|K\|_{\infty}\left(\left|Y_{k+1}\right|+\mathbb{E}\left(\left|Y_{k+1}\right|\right)\right)\right)\right] \leq c^*.
\end{eqnarray*}
Recall that $\sup_{x \in I}Var\left[Y_kK\left(\left(x-X_k\right)h_k^{-1}\right)\right]=O\left(h_k\right)$.
We then deduce that there exist positive constants $C_i^*$, such that, for all $x \in I$, and for all $k$ and $n$ such that $k \leq n$,
\begin{eqnarray*}
\mathbb{E}\left(\exp \left[s^{-1}U_{k,n}\left(x\right)\right]\right)
& \leq & 1+C_1^*\frac{nh_n}{\left(\sum_{k=1}^nq_k\right)^2}\frac{q_k^2}{h_k}+C_2^*\frac{\left(nh_n\right)^{\frac{3}{2}}}{\left(\sum_{k=1}^nq_k\right)^3}\frac{q_k^3}{h_k^3}\\
& \leq & \exp \left[C_3^*\frac{h_n}{nq_n^2}\frac{q_k^2}{h_k}+C_4^*\frac{h_n^{\frac{3}{2}}}{n^{\frac{3}{2}}q_n^3}\frac{q_k^3}{h_k^3}\right].
\end{eqnarray*}
Then, it follows from~\eqref{a29maibar} that, for all $C>0$,
\begin{eqnarray*}
\sup_{x \in I}\mathbb{P}\left[\frac{\sqrt{nh_n}}{\ln n}\left(T_{n,1}\left(x\right)-\mathbb{E}\left(T_{n,1}\left(x\right)\right)\right) \geq C\right]
&\leq & n^{-\frac{C}{s}}\exp \left[C_3^*\frac{h_n}{nq_n^2}\sum_{k=1}^n\frac{q_k^2}{h_k}+C_4^*\frac{h_n^{\frac{3}{2}}}{n^{\frac{3}{2}}q_n^3}\sum_{k=1}^n\frac{q_k^3}{h_k^3}\right]\\
& = & O\left(n^{-\frac{C}{s}}\right).
\end{eqnarray*}
We establish exactly in the same way that, for all $C>0$,
\begin{eqnarray*}
\sup_{x \in I}\mathbb{P}\left[\frac{\sqrt{nh_n}}{\ln n}\left(\mathbb{E}\left(T_{n,1}\left(x\right)\right)-T_{n,1}\left(x\right)\right) \geq C\right] &=&O\left(n^{-\frac{C}{s}}\right),
\end{eqnarray*}
which concludes the proof of Lemma~\ref{aL:ETbar}.
\subsection{Proof of Lemma~\ref{aL:rwn}}\label{a4.5}
In view of~\eqref{a26mar}, we have
\begin{eqnarray*}
\Delta_n\left(x\right)
&=&\Delta_n^{\left(1\right)}\left(x\right)+\Delta_n^{\left(2\right)}\left(x\right),
\end{eqnarray*}
with
\begin{eqnarray}
\Delta_n^{\left(1\right)}\left(x\right)&=&\sum_{k=n_0}^nU_{k,n}\left(f\left(x\right)\right)\gamma_k\left(\mathbb{E}\left[Z_k\left(x\right)\right]-Z_k\left(x\right)\right)\left(r_{k-1}(x)-r\left(x\right)\right),\label{adelt1}\\
\Delta_n^{\left(2\right)}\left(x\right)&=&\sum_{k=n_0}^nU_{k,n}\left(f\left(x\right)\right)\gamma_k\left(f\left(x\right)-\mathbb{E}\left[Z_k\left(x\right)\right]\right)\left(r_{k-1}(x)-r\left(x\right)\right).\label{adelt2}
\end{eqnarray}
Let us first note that, in view of~\eqref{aespz} and by application of Lemma~\ref{alap}, we have
\begin{eqnarray*}
\left|\Delta_n^{\left(2\right)}\left(x\right)\right|&=&O\left(\Pi_n\left(f\left(x\right)\right)\sum_{k=n_0}^n\Pi_k^{-1}\left(f\left(x\right)\right)\gamma_kh_k^2w_k\right)\quad a.s.\\
&=&O\left(\Pi_n\left(f\left(x\right)\right)\sum_{k=n_0}^n\Pi_k^{-1}\left(f\left(x\right)\right)\gamma_kO\left(m_k\right)w_k\right)\quad a.s.\\
&=&O\left(m_nw_n\right)\quad a.s.
\end{eqnarray*}
Let us now bound $\Delta_n^{\left(1\right)}\left(x\right)$. To this end, we set
\begin{eqnarray*}
\varepsilon_k\left(x\right)&=&\mathbb{E}\left(Z_k\left(x\right)\right)-Z_k\left(x\right),\\
G_{k}\left(x\right)&=&r_{k}\left(x\right)-r\left(x\right),\\
S_n\left(x\right)&=&\sum_{k=1}^n\Pi_k^{-1}\left(f\left(x\right)\right)\gamma_k\varepsilon_k\left(x\right)G_{k-1}\left(x\right),
\end{eqnarray*}
and $\mathcal{F}_k=\sigma\left(\left(X_1,Y_1\right),\ldots,\left(X_k,Y_k\right)\right)$. In view of~\eqref{avarz} and of Lemma~\ref{alap}, the increasing process of the martingale $\left(S_n\left(x\right)\right)$ satisfies
\begin{eqnarray*}
<S>_{n}\left(x\right)&=&\sum_{k=n_0}^{n}\mathbb{E}\left[\Pi_k^{-2}\left(f\left(x\right)\right)\gamma_k^2 \varepsilon_k^2\left(x\right)G_{k-1}^2\left(x\right)|\mathcal{F}_{k-1}\right]\nonumber\\
&=&\sum_{k=n_0}^{n}\Pi_k^{-2}\left(f\left(x\right)\right)\gamma_k^2G_{k-1}^2\left(x\right)Var\left[Z_k\left(x\right)\right]\nonumber\\
&=&O\left(\sum_{k=n_0}^n\Pi_k^{-2}\left(f\left(x\right)\right)\gamma_k^2w_k^2\frac{1}{h_k}\right)\quad a.s.\nonumber\\
&=&O\left(\sum_{k=n_0}^n\Pi_k^{-2}\left(f\left(x\right)\right)\gamma_km_k^2w_k^2\right)\quad a.s.\nonumber\\
&=&O\left(\Pi_n^{-2}\left(f\left(x\right)\right)m_n^2w_n^2\right)\quad a.s.
\end{eqnarray*} 
$\bullet$ Let us first consider the case the sequence $\left(n\gamma_n\right)$ is bounded. We then have $\left(\Pi_n^{-1}\left(f\left(x\right)\right)\in\mathcal{GS}\left(\xi^{-1}f\left(x\right)\right)\right)$, and thus $\ln \left(<S>_{n}\left(x\right)\right)=O\left(\ln n\right)$ a.s. Theorem 1.3.15 in Duflo (1997) then ensures that, for any $\delta >0$,
\begin{eqnarray*}
\left|S_n\left(x\right)\right|&=&o\left(<S>_n^{\frac{1}{2}} \left(x\right)\left( \ln<S>_n\left(x\right)\right)^{\frac{1+\delta}{2}}\right)+O\left(1\right)\quad a.s.\nonumber\\
&=&o\left(\Pi_n^{-1}\left(f\left(x\right)\right)m_nw_n\left(\ln n\right)^{\frac{1+\delta}{2}}\right)+O\left(1\right)\quad a.s.
\end{eqnarray*}
It follows that, for any $\delta >0$,
\begin{eqnarray*}
\left|\Delta_n^{\left(1\right)}\left(x\right)\right|
&=&o\left(m_nw_n\left(\ln n\right)^{\frac{1+\delta}{2}}\right)+O\left(\Pi_n\left(f\left(x\right)\right)\right)\quad a.s.\\
&=&o\left(m_nw_n\left(\ln n\right)^{\frac{1+\delta}{2}}\right)+o\left(m_n\right)\quad a.s.,
\end{eqnarray*}
which concludes the proof of Lemma~\ref{aL:rwn} in this case.\\
$\bullet$ Let us now consider the case $\lim_{n\to\infty}\left(n\gamma_n\right)=\infty$. In this case, for all $\delta>0$, we have
\begin{eqnarray*}
\ln \left(\Pi_n^{-2}\left(f\left(x\right)\right)\right)&=&\sum_{k=n_0}^n\ln \left(1-\gamma_kf\left(x\right)\right)^{-2}\\
&=&\sum_{k=n_0}^n\left(2\gamma_kf\left(x\right)+o\left(\gamma_k\right)\right)\\
&=&O\left(\sum_{k=1}^n\gamma_kk^{\delta}\right).
\end{eqnarray*}
Since $\left(\gamma_nn^{\delta}\right)\in\mathcal{GS}\left(-\left(\alpha-\delta\right)\right)$ with $\left(\alpha-\delta\right)<1$, we have
\begin{eqnarray*}
\lim_{n \to \infty}\frac{n\left(\gamma_nn^{\delta}\right)}{\sum_{k=1}^n\gamma_kk^{\delta}}=1-\left(\alpha-\delta \right).
\end{eqnarray*}
It follows that $\ln \left(\Pi_n^{-2}\left(f\left(x\right)\right)\right)=O\left(n^{1+\delta}\gamma_n\right)$. The sequence $\left(m_nw_n\right)$ being in $\mathcal{GS}\left(-m^*+w^*\right)$, we deduce that, for all $\delta>0$, we have
\begin{eqnarray*}
\ln \left(<S>_{n}\left(x\right)\right)&=&O\left(n^{1+\delta}\gamma_n\right)\quad a.s.
\end{eqnarray*}
Theorem 1.3.15 in Duflo (1997) then ensures that, for any $\delta >0$,
\begin{eqnarray*}
\left|S_n\left(x\right)\right|&=&o\left(<S>_n^{\frac{1}{2}}\left(x\right) \left( \ln<S>_n\left(x\right)\right)^{\frac{1+\delta}{2}}\right)+O\left(1\right)\quad a.s.\nonumber\\
&=&o\left(\Pi_n^{-1}\left(f\left(x\right)\right)m_nw_n\left(n^{1+\delta}\gamma_n\right)^{\frac{1+\delta}{2}}\right)+O\left(1\right)\quad a.s.
\end{eqnarray*}
It follows from the application of Lemma~\ref{alap} that, for any $\delta>0$,
\begin{eqnarray*}
\left|\Delta_n^{\left(1\right)}\left(x\right)\right|&=&o\left(m_nw_n\left(n^{1+\delta}\gamma_n\right)^{\frac{1+\delta}{2}}\right)+O\left(\Pi_n\left(f\left(x\right)\right)\right)\quad a.s.\\
&=&o\left(m_nw_n\left(n^{1+\delta}\gamma_n\right)^{\frac{1+\delta}{2}}\right)\quad a.s.,
\end{eqnarray*}
which concludes the proof of Lemma~\ref{aL:rwn}.
\subsection{Proof of Lemma~\ref{aL:rwn2}}\label{a4.6}
Let us first note that, in view of~\eqref{adelt2}, and by application of Lemma~\ref{alap}, we have
\begin{eqnarray*}
\sup_{x\in I}\left|\Delta_n^{\left(2\right)}\left(x\right)\right|&=&O\left(\sum_{k=n_0}^n\left(\sup_{x\in I}U_{k,n}\left(f\left(x\right)\right)\right)\gamma_kh_k^2w_k\right)\quad a.s.\\
&=&O\left(\sum_{k=n_0}^nU_{k,n}\left(\varphi\right)\gamma_km_kw_k\right)\quad a.s.\\
&=&O\left(m_nw_n\right)\quad a.s.
\end{eqnarray*}
Now, set 
\begin{eqnarray}\label{aAn*}
A_n=\frac{3}{t^*}\ln n
\end{eqnarray}
(where $t^*$ is defined in $\left(A4\right)ii)$), and write $\Delta_n^{\left(1\right)}$ 
(defined in (\ref{adelt1})) as
\begin{eqnarray*}
\Delta_n^{\left(1\right)}\left(x\right)=\Pi_n\left(f\left(x\right)\right)M_n^{\left(n\right)}\left(x\right)+\Pi_n\left(f\left(x\right)\right)S_n\left(x\right)
\end{eqnarray*}
with
\begin{eqnarray*}
S_n\left(x\right)&=&\sum_{k=n_0}^n\Pi_k^{-1}\left(f\left(x\right)\right)\gamma_k\left(\mathbb{E}\left[Z_k\left(x\right)\right]-Z_k\left(x\right)\right)\left(r_{k-1}(x)-r\left(x\right)\right)\mathds{1}_{\sup_{l\leq k-1}\left|Y_l\right|> A_n},\nonumber\\
M_k^{\left(n\right)}\left(x\right)&=&\sum_{j=n_0}^k\Pi_j^{-1}\left(f\left(x\right)\right)\gamma_j\left(\mathbb{E}\left[Z_j\left(x\right)\right]-Z_j\left(x\right)\right)\left(r_{j-1}(x)-r\left(x\right)\right)\mathds{1}_{\sup_{l\leq j-1}\left|Y_l\right|\leq A_n}.
\end{eqnarray*}
Let us first prove a uniform strong upper bound for $S_n$. For any $c>0$, we have
\begin{eqnarray*}
\sum_{n\geq 0}\mathbb{P}\left[\sup_{x \in I}m_n^{-1}w_n^{-1}\left|S_n\left(x\right)\right|\geq c\right]&=&O\left(\sum_{n\geq 0}\mathbb{P}\left(\sup_{l\leq n-1}\left|Y_l\right|>A_n\right)\right)\\
&=&O\left(\sum_{n\geq 0}n\mathbb{P}\left(\left|Y\right|>A_n\right)\right)\\
&=&O\left(\sum_{n\geq 0}n\exp\left(-t^*A_n\right)\right)\\
&<&\infty.
\end{eqnarray*}
It follows that
\begin{eqnarray*}
\sup_{x \in I}\left|S_n\left(x\right)\right|&=&O\left(m_nw_n\right)\quad a.s.
\end{eqnarray*}
To establish the strong uniform bound of $M_n^{\left(n\right)}$, we shall apply the following result given in Duflo $(1997)$, page $209$. 
\begin{lemma}\label{aL:phi} $ $\\
Let $\left(M_k^{\left(n\right)}\right)_k$ be a martingale such that, for all $k \leq n$, $\left|M_k^{\left(n\right)}-M_{k-1}^{\left(n\right)}\right|\leq c_n$, and set \\$\Phi_c\left(\lambda\right)=c^{-2}\left(e^{\lambda c}-1-\lambda c\right)$. For all $\lambda_n $ such that $\lambda_n c_n \leq 1$ and all $\alpha_n >0$, we have
\begin{eqnarray*}
 \mathbb{P}\left(\lambda_n\left( M_n^{\left(n\right)}-M_0^{\left(n\right)}\right) \geq \Phi_{c_n}\left(\lambda_n\right)<M^{\left(n\right)}>_n+\alpha_n \lambda_n\right) \leq e^{-\alpha_n \lambda_n}.
\end{eqnarray*}
\end{lemma}
In view of~\eqref{abis22}, there exists $C^*>0$ such that, on $\left\{\sup_{l\leq k}\left|Y_l\right|\leq A_n\right\}$,
\begin{eqnarray*}
\left|r_k\left(x\right)-r\left(x\right)\right|&\leq&C^*k\gamma_kh_k^{-1}A_n.
\end{eqnarray*}
Consequently, there exists $C_1>0$ such that 
\begin{eqnarray*}
\left|M_k^{\left(n\right)}\left(x\right)-M_{k-1}^{\left(n\right)}\left(x\right)\right|&\leq &\Pi_k^{-1}\left(f\left(x\right)\right)\gamma_k\left|Z_k\left(x\right)-\mathbb{E}\left(Z_k\left(x\right)\right)\right|\left|\left(r_{k-1}\left(x\right)-r\left(x\right)\right)\mathds{1}_{\sup_{l\leq k-1}\left|Y_l\right|\leq A_n}\right|\\
&\leq & \Pi_k^{-1}\left(f\left(x\right)\right)\gamma_k\left(2h_{k}^{-1}\|K\|_{\infty}\right)\left(C^*\left(k-1\right)\gamma_{k-1}h_{k-1}^{-1}A_n\right)\\
& \leq &C_1\Pi_k^{-1}\left(f\left(x\right)\right)k\gamma_k^2h_k^{-2}A_n.
\end{eqnarray*}
\begin{itemize}
\item In the case $\lim_{n\to\infty}\left(n\gamma_n\right)=\infty$, since $\left(n\gamma_n^2h_n^{-2}\right)\in \mathcal{GS}\left(1-2\alpha+2a\right)$ there exists $\left(u_k\right)\to 0$ such that
\begin{eqnarray*}
\frac{\left(k-1\right)\gamma_{k-1}^2h_{k-1}^{-2}}{k\gamma_{k}^2h_{k}^{-2}}
&=&1-\left[1-2\alpha+2a\right]\frac{1}{k}+o\left(\frac{1}{k}\right)\nonumber\\
&=&1+u_k\gamma_k.
\end{eqnarray*}
It follows that there exists $k_0\geq n_0$ such that, for all $k\geq k_0$ and for all $x\in I$, 
\begin{eqnarray*}
\frac{\Pi_{k-1}^{-1}\left(f\left(x\right)\right)\left(k-1\right)\gamma_{k-1}^2h_{k-1}^{-2}}{\Pi_k^{-1}\left(f\left(x\right)\right)k\gamma_{k}^2h_{k}^{-2}}
&=&\left(1-\gamma_kf\left(x\right)\right)\left(1+u_k\gamma_k\right)\nonumber\\
&=&1-\gamma_kf\left(x\right)+u_k\gamma_k\left(1-\gamma_kf\left(x\right)\right)\nonumber\\
&\leq&1-\gamma_k\varphi+u_k\gamma_k\left(1+\gamma_k\|f\|_{\infty}\right)\nonumber\\
&\leq &1.
\end{eqnarray*}
Consequently, there exists $C>0$ such that, for all $x\in I$ and all $k\leq n$,
\begin{eqnarray}\label{a73bis}
\left|M_k^{\left(n\right)}\left(x\right)-M_{k-1}^{\left(n\right)}\left(x\right)\right|& \leq &C\Pi_n^{-1}\left(f\left(x\right)\right)n\gamma_n^2h_n^{-2}A_n.
\end{eqnarray}
\item In the case $\lim_{n\to\infty}\left(n\gamma_n\right)<\infty$ (in which case $\alpha=1$), we set $\epsilon \in \left]0,\min\left\{\left(1-3a\right)/2;\varphi\xi^{-1}-m^*\right\}\right[$ (where $m^*$ is defined in~\eqref{a10bis}), and write
\begin{eqnarray*}
\left|M_k^{\left(n\right)}\left(x\right)-M_{k-1}^{\left(n\right)}\left(x\right)\right|& \leq &C_1\left[\Pi_k^{-1}\left(f\left(x\right)\right)k^{-\epsilon} m_k\right]A_n\left[m_k^{-1}k^{1+\epsilon}\gamma_k^2h_k^{-2}\right].
\end{eqnarray*}
Since $\left(m_n^{-1}n^{1+\epsilon}\gamma_n^2h_n^{-2}\right)\in \mathcal{GS}\left(m^*+1+\epsilon-2\alpha+2a\right)$ with 
\begin{eqnarray*}
m^*+1+\epsilon-2\alpha+2a &\leq &\frac{1-a}{2}+\epsilon-1+2a\\
&\leq & \epsilon-\frac{1}{2}\left(1-3a\right)<0,
\end{eqnarray*}
the sequence $\left(m_n^{-1}n^{1+\epsilon}\gamma_n^2h_n^{-2}\right)$ is bounded. On the other hand, since $\left(n^{-\epsilon}m_n\right)\in\mathcal{GS}\left(-\epsilon-m^*\right)$, there exists $\left(u_k\right)\to 0$ and $k_0\geq n_0$ such that, for all $k\geq k_0$ and for all $x \in I$,
\begin{eqnarray*}
\frac{\Pi_{k-1}^{-1}\left(f\left(x\right)\right)\left(k-1\right)^{-\epsilon}m_{k-1}}{\Pi_k^{-1}\left(f\left(x\right)\right)k^{-\epsilon}m_{k}}&=&\left(1-\gamma_kf\left(x\right)\right)\left(1+\left(m^*+\epsilon \right)\frac{1}{k}+o\left(\frac{1}{k}\right)\right)\nonumber\\
&=&\left(1-\gamma_kf\left(x\right)\right)\left(1+\left(m^*+\epsilon \right)\xi \gamma_k+u_k\gamma_k\right)\nonumber\\
&\leq&\left(1-\gamma_k\varphi\right)\left(1+\left(m^*+\epsilon \right)\xi\gamma_k+\left|u_k\right|\gamma_k\right)\nonumber\\
&\leq&1-\frac{\left(\varphi-\left(m^*+\epsilon \right)\xi\right)\gamma_k}{2}\nonumber\\
&\leq &1.
\end{eqnarray*}
It follows that there exists $C>0$ such that, for all $x\in I$ and all $k\leq n$, 
\begin{eqnarray}\label{a73ter}
\left|M_k^{\left(n\right)}\left(x\right)-M_{k-1}^{\left(n\right)}\left(x\right)\right|& \leq &C\Pi_n^{-1}\left(f\left(x\right)\right)n^{-\epsilon}m_nA_n.
\end{eqnarray}
\end{itemize}
From now on, we set
\begin{eqnarray*}
c_n\left(x\right)= \left\{\begin{array}{llll}C\Pi_n^{-1}\left(f\left(x\right)\right)n\gamma_n^2h_n^{-2}A_n & \mbox{ if }\lim_{n\to\infty}\left(n\gamma_n\right)=\infty,\\
C\Pi_n^{-1}\left(f\left(x\right)\right)m_nn^{-\epsilon}A_n & \mbox{ if }\lim_{n\to\infty}\left(n\gamma_n\right)<\infty,
\end{array}\right.
\end{eqnarray*}
so that in view of~\eqref{a73bis} and~\eqref{a73ter}, for all $x\in I$ and all $k\leq n$, we have
\begin{eqnarray*}
\left|M_k^{\left(n\right)}\left(x\right)-M_{k-1}^{\left(n\right)}\left(x\right)\right|& \leq & c_n\left(x\right).
\end{eqnarray*}
Now, let $\left(u_n\right)$ be a positive sequence such that, for all $n$,
\begin{eqnarray}\label{aun12}
\left\{\begin{array}{llll}u_n\leq C^{-1}n^{-1}\gamma_n^{-2}h_n^{2}A_n^{-1} & \mbox{ if }\lim_{n\to\infty}\left(n\gamma_n\right)=\infty,\\
u_n \leq C^{-1}m_n^{-1}n^{\epsilon}A_n^{-1} & \mbox{ if }\lim_{n\to\infty}\left(n\gamma_n\right)<\infty,
\end{array}\right.
\end{eqnarray}
and set
\begin{eqnarray*}
\lambda_n\left(x\right)=u_n\Pi_n\left(f\left(x\right)\right).
\end{eqnarray*}
Let us at first assume that the following lemma holds.
\begin{lemma}\label{aL:Mnc} $ $\\
There exist $C_2>0$ and $\rho>0$ such that for all $x,y\in I$ such that $\left|x-y\right| \leq C_2 n^{-\rho}$, we have
\begin{eqnarray*}
\left|\lambda_n\left(x\right)M_n^{\left(n\right)}\left(x\right)-\lambda_n\left(y\right)M_n^{\left(n\right)}\left(y\right)\right|& \leq & 1,\\
\left|\Phi_{c_n\left(x\right)}\left(\lambda_n\left(x\right)\right)<M^{\left(n\right)}>_n\left(x\right)-\Phi_{c_n\left(y\right)}\left(\lambda_n\left(y\right)\right)<M^{\left(n\right)}>_n\left(y\right)\right|& \leq & 1.
\end{eqnarray*}
\end{lemma}
Set
\begin{eqnarray*}
d_n&=&C_2n^{-\rho},\nonumber\\
V_n\left(x\right)&=&\lambda_n\left(x\right)M_n^{\left(n\right)}\left(x\right)-\Phi_{c_n\left(x\right)}\left(\lambda_n\left(x\right)\right)<M^{\left(n\right)}>_n\left(x\right),\nonumber\\
\alpha_n\left(x\right)&=&\frac{\left(\rho+2\right)\ln n}{\lambda_n\left(x\right)}.
\end{eqnarray*}
Let $I_i^{\left(n\right)}$ be $N\left(n\right)$ intervals of length $d_n$ such that $\cup_{i=1}^{N\left(n\right)}=I$, and for all $i\in \left\{1,\ldots,N\left(n\right)\right\}$, set $x_i^{\left(n\right)}\in I_i^{\left(n\right)}$. Applying Lemma~\ref{aL:Mnc}, we get, for $n$ large enough,
\begin{eqnarray*}
\mathbb{P}\left[\sup_{x \in I}V_n\left(x\right) \geq 2\left(\rho+2\right)\ln n\right]
& \leq & \sum_{i=1}^{N\left(n\right)}\mathbb{P}\left[\sup_{x \in I_i^{\left(n\right)}}V_n\left(x\right)\geq 2\left(\rho+2\right)\ln n\right]\nonumber\\
& \leq & \sum_{i=1}^{N\left(n\right)}\mathbb{P}\left[V_n\left(x_i^{\left(n\right)}\right)+2\geq 2\left(\rho+2\right)\ln n\right]\nonumber\\
& \leq & N\left(n\right)\sup_{x \in I}\mathbb{P}\left[V_n\left(x\right)\geq \left(\rho+2\right)\ln n\right].
\end{eqnarray*}
Now, the application of Lemma $\ref{aL:phi}$ ensures that, for all $x\in I$,
\begin{eqnarray*}
\mathbb{P}\left[V_n\left(x\right) \geq \left(\rho+2\right)\ln n\right]
& \leq & \mathbb{P}\left[\lambda_n\left(x\right)M_n^{\left(n\right)}\left(x\right)-\Phi_{c_n\left(x\right)}\left(\lambda_n\left(x\right)\right)<M^{\left(n\right)}>_n\left(x\right) \geq \alpha_n\left(x\right) \lambda_n\left(x\right)\right]\nonumber\\
& \leq & \exp\left[-\alpha_n\left(x\right)\lambda_n\left(x\right)\right]\nonumber\\
& \leq & n^{-\left(\rho+2\right)}.
\end{eqnarray*}
It follows that
\begin{eqnarray*}
\sum_{n\geq 1}\mathbb{P}\left[\sup_{x \in I}V_n\left(x\right) \geq 2\left(\rho+2\right)\ln n\right]=O\left(\sum_{n\geq 1}n^{-2}\right)<+\infty,
\end{eqnarray*}
and, applying Borel-Cantelli Lemma, we obtain 
\begin{eqnarray*}
\sup_{x \in I}\lambda_n\left(x\right)M_n^{\left(n\right)}\left(x\right)& \leq & \sup_{x \in I}\Phi_{c_n\left(x\right)}\left(\lambda_n\left(x\right)\right)<M^{\left(n\right)}>_n\left(x\right)+2\left(\rho +2\right)\ln n \quad a.s.
\end{eqnarray*}
Since $\Phi_{c}\left(\lambda\right)\leq \lambda^2$ as soon as $\lambda c \leq 1$, and since $\lambda_n\left(x\right)=u_n\Pi_n\left(f\left(x\right)\right)$, it follows that
\begin{eqnarray*}
u_n\sup_{x \in I}\Pi_n\left(f\left(x\right)\right)M_n^{\left(n\right)}\left(x\right)& \leq & u_n^2\sup_{x \in I}\Pi_n^2\left(f\left(x\right)\right)<M^{\left(n\right)}>_n\left(x\right)+2\left(\rho +2\right)\ln n \quad a.s.
\end{eqnarray*}
Establishing the same upper bound for the martingale $\left(-M_k^{\left(n\right)}\right)$, we obtain
\begin{eqnarray*}
\sup_{x \in I}\Pi_n\left(f\left(x\right)\right)\left|M_n^{\left(n\right)}\left(x\right)\right|& \leq & u_n \sup_{x \in I}\Pi_n^2\left(f\left(x\right)\right)<M^{\left(n\right)}>_n\left(x\right)+2\frac{\left(\rho +2\right)\ln n}{u_n}\quad a.s.
\end{eqnarray*}
Now, since $\sup_{x\in I}Var\left(Z_k\left(x\right)\right)=O\left(h_k^{-1}\right)$, we have
\begin{eqnarray}\label{a<M*>}
\lefteqn{\sup_{x \in I}\Pi_n^{2}\left(f\left(x\right)\right)<M^{\left(n\right)}>_{n}\left(x\right)}\nonumber\\
&=&O\left(\sum_{k=n_0}^n\sup_{x \in I}U_{k,n}^2\left(f\left(x\right)\right)\gamma_k^2\sup_{x \in I}\left|r_{k-1}\left(x\right)-r\left(x\right)\right|^2\sup_{x \in I}\left(Var\left[Z_k\left(x\right)\right]\right)\right)\nonumber\\
&=&O\left(\sum_{k=n_0}^nU_{k,n}^2\left(\varphi\right)\gamma_k^2h_k^{-1}w_k^2\right)\quad a.s.
\end{eqnarray}
$\bullet$ Let us first consider the case when the sequence $\left(n\gamma_n\right)$ is bounded. In this case,~\eqref{a<M*>} and Lemma~\ref{alap} imply that
\begin{eqnarray*}
\sup_{x \in I}\Pi_n^{2}\left(f\left(x\right)\right)<M^{\left(n\right)}>_{n}\left(x\right)&=&O\left(\sum_{k=n_0}^nU_{k,n}^2\left(\varphi\right)\gamma_km_k^2w_k^2\right)\quad a.s.\nonumber\\
&=&O\left(m_n^2w_n^2\right)\quad a.s.
\end{eqnarray*} 
In this case, we have thus proved that, for all positive sequence $\left(u_n\right)$ satisfying~\eqref{aun12}, we have 
\begin{eqnarray}\label{a3jul}
\sup_{x\in I}\Pi_n\left(f\left(x\right)\right)\left|M_{n}^{\left(n\right)}\left(x\right)\right|=O\left(u_nm_n^2w_n^2+ \frac{\ln n}{u_n}\right)\quad a.s.
\end{eqnarray}
Now, since the sequence $\left(\left[m_n^{-1}w_n^{-1}\sqrt{\ln n}\right]m_nn^{-\epsilon}A_n\right)$ belongs to $\mathcal{GS}\left(-\left(w^*+\epsilon\right)\right)$ with $w^*+\epsilon>0$, there exists $u_0>0$ such that, for all $n$,
\begin{eqnarray*}
u_0m_n^{-1}w_n^{-1}\sqrt{\ln n}\leq C^{-1}m_n^{-1}n^{\epsilon}A_n^{-1}
\end{eqnarray*}
(where $C$ is defined in~\eqref{aun12}). Applying~\eqref{a3jul} with $\left(u_n\right)=\left(u_0m_n^{-1}w_n^{-1}\sqrt{\ln n}\right)$, we obtain
\begin{eqnarray*}
\sup_{x\in I}\Pi_n\left(f\left(x\right)\right)\left|M_{n}^{\left(n\right)}\left(x\right)\right|&=&O\left(m_nw_n\sqrt{\ln n}\right)\quad a.s.,
\end{eqnarray*}
which concludes the proof of Lemma~\ref{aL:rwn2} in this case.\\
$\bullet$ Let us now consider the case $\lim_{n\to\infty}\left(n\gamma_n\right)=\infty $. In this case,~\eqref{a<M*>} implies that
\begin{eqnarray*}
\sup_{x\in I}\Pi_n^2\left(f\left(x\right)\right)<M>_{n}^{\left(n\right)}\left(x\right)&=&O\left(\gamma_nh_n^{-1}w_n^2\right)\quad a.s.
\end{eqnarray*}
In this case, we have thus proved that, for all positive sequence $\left(u_n\right)$ satisfying~\eqref{aun12}, we have
\begin{eqnarray}\label{a21jul}
\sup_{x\in I}\Pi_n\left(f\left(x\right)\right)\left|M_{n}^{\left(n\right)}\left(x\right)\right|&=&O\left(u_n\gamma_nh_n^{-1}w_n^2+\frac{\ln n}{u_n}\right)\quad a.s.
\end{eqnarray}
Now, in view of~\eqref{a30bis}, of~\eqref{aAn*}, and of the assumptions of Lemma \ref{aL:rwn2}, we have
\begin{eqnarray*}
\left[\sqrt{\gamma_n^{-1}h_n\ln n}\,w_n^{-1}\right]\left[n\gamma_n^2h_n^{-2}A_n\right]&=&O\left(\sqrt{\gamma_nh_n^{-1}\ln n}\,w_n^{-1}B_n\right)\\
&=&O\left(1\right).
\end{eqnarray*}
Thus, there exists $u_0>0$ such that, for all $n$,
\begin{eqnarray*}
u_0\sqrt{\gamma_n^{-1}h_n\ln n}\,w_n^{-1}\leq C^{-1}n^{-1}\gamma_n^{-2}h_n^2A_n^{-1}.
\end{eqnarray*}
Applying~\eqref{a21jul} with $\left(u_n\right)=\left(u_0\sqrt{\gamma_n^{-1}h_n\ln n}\,w_n^{-1}\right)$, we obtain
\begin{eqnarray*}
\sup_{x\in I}\Pi_n\left(f\left(x\right)\right)\left|M_{n}^{\left(n\right)}\left(x\right)\right|&=&O\left(\sqrt{\gamma_nh_n^{-1}\ln n}\,w_n\right)\quad a.s.\\
&=&O\left(m_nw_n\sqrt{\ln n}\right)\quad a.s.,
\end{eqnarray*}
which concludes the proof of Lemma~\ref{aL:rwn2}.

\paragraph{Proof of Lemma~\ref{aL:Mnc}} $ $\\
Let $(\delta_n)\in{\cal GS}(-\delta^*)$, 
set $x,y \in I$ such that $\left|x-y\right|\leq\delta_n$, and let $c_i^*$ denote generic 
constants. Let us first note that 
\begin{eqnarray*}
 & & \left|Z_k\left(x\right)-Z_k\left(y\right)\right|
\leq c_1^*\delta_nh_k^{-2}, \\
& & \left|Z_k\left(x\right)\right|\leq c_2^*h_k^{-1}, \\
& & 
\left|Var\left[Z_k\left(x\right)\right]-Var\left[Z_k\left(y\right)\right]\right| \\
&  & \mbox{ } =
\left|\mathbb{E}\left\{\left(Z_k(x)-Z_k(y)-[\mathbb{E}(Z_k(x))-\mathbb{E}(Z_k(y))]\right)
\left(Z_k(x)+Z_k(y)-[\mathbb{E}(Z_k(x))+\mathbb{E}(Z_k(y))]\right)\right\}\right|
\\
& & \mbox{ } \leq 
8c_2^*h_k^{-1}\mathbb{E}\left[\left|Z_{k}\left(x\right)-Z_{k}\left(y\right)\right|
\right]\\
& & \mbox{ } \leq  c_3^*h_k^{-3}\delta_n,
\end{eqnarray*}
and that, in view of (\ref{abis22}), on 
$\left\{\sup_{l\leq k}\left|Y_l\right|\leq A_n\right\}$,
$$
\left|r_k\left(x\right)-r\left(x\right)\right|\leq c_4^*k\g_kh_k^{-1}A_n.
$$
Now, in view of~\eqref{arnMar}, we have 
\begin{eqnarray*}
r_k\left(x\right)-r\left(x\right)=\left[1-\gamma_kZ_k\left(x\right)\right]
\left[r_{k-1}\left(x\right)-r\left(x\right)\right]
+\gamma_k\left[W_k\left(x\right)-r\left(x\right)Z_k\left(x\right)\right],
\end{eqnarray*}
so that 
\begin{eqnarray*}
\lefteqn{\left[r_k\left(x\right)-r\left(x\right)\right]-
\left[r_k\left(y\right)-r\left(y\right)\right]}\nonumber\\
&=&\left[1-\gamma_kZ_k\left(x\right)\right]
\left[r_{k-1}\left(x\right)-r\left(x\right)\right]
-\left[1-\gamma_kZ_k\left(y\right)\right]
\left[r_{k-1}\left(y\right)-r\left(y\right)\right]\nonumber\\
&&+\gamma_k\left(\left[W_k\left(x\right)-W_k\left(y\right)\right]
-\left[r\left(x\right)Z_k\left(x\right)-r\left(y\right)Z_k\left(y\right)\right]
\right)\nonumber\\
&=&\left[1-\gamma_kZ_k\left(x\right)\right]\left(\left[r_{k-1}\left(x\right)
-r\left(x\right)\right]-\left[r_{k-1}\left(y\right)-r\left(y\right)\right]\right)\\
&&-\gamma_k\left[Z_k\left(x\right)-Z_k\left(y\right)\right]
\left[r_{k-1}\left(y\right)-r\left(y\right)\right]\nonumber\\
&&+\gamma_k\left(\left[W_k\left(x\right)-W_k\left(y\right)\right]
-r\left(x\right)\left[Z_k\left(x\right)-Z_k\left(y\right)\right]
-Z_k\left(y\right)\left[r\left(x\right)-r\left(y\right)\right]\right).
\end{eqnarray*}
For $k\geq n_0$, we have $\left|1-\gamma_kZ_k\left(x\right)\right|\leq 1$. It follows that, 
for $k\geq n_0$ and on $\left\{\sup_{l\leq k}\left|Y_l\right|\leq A_n\right\}$,
\begin{eqnarray*}
\lefteqn{\left|\left[r_k\left(x\right)-r\left(x\right)\right]
-\left[r_k\left(y\right)-r\left(y\right)\right]\right|}\nonumber\\
&\leq&
\left|\left[r_{k-1}\left(x\right)-r\left(x\right)\right]
-\left[r_{k-1}\left(y\right)-r\left(y\right)\right]\right|
+\gamma_k\left|Z_k\left(x\right)-Z_k\left(y\right)\right|\left|r_{k-1}\left(y\right)-r\left(y\right)\right|\nonumber\\
&&+\gamma_k\left(\left|Y_k\right|+\left|r\left(x\right)\right|\right)\left|Z_k\left(x\right)-Z_k\left(y\right)\right|+\gamma_k\left|Z_k\left(y\right)\right|\left|r\left(x\right)-r\left(y\right)\right|\nonumber\\
&\leq& \left|\left[r_{k-1}\left(x\right)-r\left(x\right)\right]
-\left[r_{k-1}\left(y\right)-r\left(y\right)\right]\right|
+\left(c^*_5k\gamma_k^2h_k^{-3}\delta_nA_n\right)
+\left(c^*_6\gamma_kA_nh_k^{-2}\delta_n\right)+\left(c^*_7\gamma_kh_k^{-1}\delta_n\right)\nonumber\\
& \leq & \left|\left[r_{k-1}\left(x\right)-r\left(x\right)\right]
-\left[r_{k-1}\left(y\right)-r\left(y\right)\right]\right|
+c^*_8k\gamma_k^2h_k^{-3}\delta_nA_n \\
& \leq & c^*_9\sum_{j=1}^kj\gamma_j^2h_j^{-3}\delta_nA_n\\
& \leq & c^*_{10}k^2\gamma_k^2h_k^{-3}\delta_nA_n.
\end{eqnarray*}
Moreover, we note that, on $\left\{\sup_{l\leq k}\left|Y_l\right|\leq A_n\right\}$,
\begin{eqnarray*}
\lefteqn{\left|\left[r_k\left(x\right)-r\left(x\right)\right]^2
-\left[r_k\left(y\right)-r\left(y\right)\right]^2\right|}\nonumber\\
&\leq &\left|\left[r_k\left(x\right)-r\left(x\right)\right]
-\left[r_k\left(y\right)-r\left(y\right)\right]\right|
\left|\left[r_k\left(x\right)-r\left(x\right)\right]
+\left[r_k\left(y\right)-r\left(y\right)\right]\right|\nonumber\\
& \leq &c^*_{11}k^3\gamma_k^2h_k^{-4}\delta_nA_n^2.
\end{eqnarray*} 
Using all the previous upper bounds, as well as \eqref{a9ju} with $p=1$, we get 
\begin{eqnarray*}
\lefteqn{\left|\lambda_n\left(x\right)M_n^{\left(n\right)}\left(x\right)-\lambda_n\left(y\right)M_n^{\left(n\right)}\left(y\right)\right|}\nonumber\\
&=&u_n\left|\sum_{k=n_0}^nU_{k,n}\left(f\left(x\right)\right)\gamma_k\left(\mathbb{E}\left[Z_k\left(x\right)\right]-Z_k\left(x\right)\right)\left[r_{k-1}\left(x\right)-r\left(x\right)\right]\right.\nonumber\\
&&-\left.\sum_{k=n_0}^nU_{k,n}\left(f\left(y\right)\right)\gamma_k\left(\mathbb{E}\left[Z_k\left(y\right)\right]-Z_k\left(y\right)\right)\left[r_{k-1}\left(y\right)-r\left(y\right)\right]\right|\mathds{1}_{\sup_{l\leq k-1}\left|Y_l\right|\leq A_n}\nonumber\\
&\leq &
u_n\sum_{k=n_0}^nU_{k,n}\left(f\left(x\right)\right)\gamma_k
\left|\mathbb{E}\left[Z_k\left(x\right)\right]-Z_k\left(x\right)\right|
\left|\left(r_{k-1}\left(x\right)-r\left(x\right)\right)-
\left(r_{k-1}\left(y\right)-r\left(y\right)\right)\right|
\mathds{1}_{\sup_{l\leq k-1}\left|Y_l\right|\leq A_n}\nonumber\\
&&+u_n\sum_{k=n_0}^nU_{k,n}\left(f\left(x\right)\right)\gamma_k
\left|r_{k-1}\left(y\right)-r\left(y\right)\right|
\left(\left|Z_{k}\left(x\right)-Z_{k}\left(y\right)\right|
+\mathbb{E}\left[\left|Z_{k}\left(x\right)-Z_{k}\left(y\right)\right|\right]\right)
\mathds{1}_{\sup_{l\leq k-1}\left|Y_l\right|\leq A_n}\nonumber\\
&&+u_n\sum_{k=n_0}^nU_{k,n}\left(f\left(y\right)\right)
\left|\frac{U_{k,n}\left(f\left(x\right)\right)}{U_{k,n}\left(f\left(y\right)\right)}-1\right|
\gamma_k\left|r_{k-1}\left(y\right)-r\left(y\right)\right|
\left(\left|Z_{k}\left(x\right)\right|
+\mathbb{E}\left[\left|Z_{k}\left(x\right)\right|\right]\right)
\mathds{1}_{\sup_{l\leq k-1}\left|Y_l\right|\leq A_n}\nonumber\\
&\leq &
c^*_{12}u_n\sum_{k=n_0}^nU_{k,n}\left(\varphi\right)\gamma_k
\left(h_k^{-1}\right)\left(k^2\gamma_k^2h_k^{-3}\delta_nA_n\right)
+c^*_{13}u_n\sum_{k=n_0}^nU_{k,n}\left(\varphi\right)\gamma_k(k\gamma_kh_k^{-1}A_n)
\left(\delta_nh_k^{-2}\right)\nonumber\\
&&+c^*_{14}u_n\sum_{k=n_0}^nU_{k,n}\left(\varphi\right)\left(\delta_n\gamma_n^{-1}\right)
\gamma_k(k\gamma_kh_k^{-1}A_n)\left(h_k^{-1}\right)\nonumber\\
&\leq &c^*_{15}u_nn^2\gamma_n^2h_n^{-4}\delta_nA_n+c^*_{16}u_nn\gamma_nh_n^{-3}\delta_nA_n
+c^*_{17}u_nnh_n^{-2}\delta_nA_n\nonumber
\end{eqnarray*}
In view of (\ref{aAn*}) and (\ref{aun12}), it follows that there exist $s_1^*>0$ and 
$\tilde S_n^{(1)}\in{\cal GS}(s_1^*)$ such that 
\begin{eqnarray}\label{aMn}
\left|\lambda_n\left(x\right)M_n^{\left(n\right)}\left(x\right)
-\lambda_n\left(y\right)M_n^{\left(n\right)}\left(y\right)\right|
&\leq & \delta_n \tilde S_n^{(1)}.
\end{eqnarray}
Now, we have
\begin{eqnarray*}
\lefteqn{\Phi_{c_n\left(x\right)}\left(\lambda_n\left(x\right)\right)<M^{\left(n\right)}>_n\left(x\right)-\Phi_{c_n\left(y\right)}\left(\lambda_n\left(y\right)\right)<M^{\left(n\right)}>_n\left(y\right)}\nonumber\\
&=&\frac{\Phi_{c_n\left(x\right)}\left(\lambda_n\left(x\right)\right)}{\lambda_n^2\left(x\right)}\lambda_n^2\left(x\right)<M^{\left(n\right)}>_n\left(x\right)-\frac{\Phi_{c_n\left(y\right)}\left(\lambda_n\left(y\right)\right)}{\lambda_n^2\left(y\right)}\lambda_n^2\left(y\right)<M^{\left(n\right)}>_n\left(y\right)\nonumber\\
&= & \frac{\Phi_{c_n\left(x\right)}\left(\lambda_n\left(x\right)\right)}{\lambda_n^2\left(x\right)}\left[\lambda_n^2\left(x\right)<M^{\left(n\right)}>_n\left(x\right)-\lambda_n^2\left(y\right)<M^{\left(n\right)}>_n\left(y\right)\right]\nonumber\\
&&+\left[\frac{\Phi_{c_n\left(x\right)}\left(\lambda_n\left(x\right)\right)}{\lambda_n^2\left(x\right)}-\frac{\Phi_{c_n\left(y\right)}\left(\lambda_n\left(y\right)\right)}{\lambda_n^2\left(y\right)}\right]\lambda_n^2\left(y\right)<M^{\left(n\right)}>_n\left(y\right).
\end{eqnarray*}
Since $c_n\left(x\right)\lambda_n\left(x\right)=c_n\left(y\right)\lambda_n\left(y\right)=\tilde B_n$, 
we have
\begin{eqnarray*}
\frac{\Phi_{c_n\left(x\right)}\left(\lambda_n\left(x\right)\right)}{\lambda_n^2\left(x\right)}
=\tilde B_n^{-2}\left(\exp\left(\tilde B_n\right)-1-\tilde B_n\right)
=\frac{\Phi_{c_n\left(y\right)}\left(\lambda_n\left(y\right)\right)}{\lambda_n^2\left(y\right)}.
\end{eqnarray*}
Using the fact that $\Phi_{c}\left(\lambda\right)\leq \lambda^2$ for $\lambda c\leq 1$, 
and applying \eqref{a9ju} with $p=2$, we deduce that 
\begin{eqnarray*}
\lefteqn{\left|\Phi_{c_n\left(x\right)}\left(\lambda_n\left(x\right)\right)<M^{\left(n\right)}>_n\left(x\right)-\Phi_{c_n\left(y\right)}\left(\lambda_n\left(y\right)\right)<M^{\left(n\right)}>_n\left(y\right)\right|}\nonumber\\
& \leq & \left|\lambda_n^2\left(x\right)<M^{\left(n\right)}>_n\left(x\right)-\lambda_n^2\left(y\right)<M^{\left(n\right)}>_n\left(y\right)\right|\nonumber\\
&\leq&u_n^2\left|\sum_{k=n_0}^nU_{k,n}^2\left(f\left(x\right)\right)\gamma_k^2Var\left[Z_k\left(x\right)\right]\left[r_{k-1}\left(x\right)-r\left(x\right)\right]^2\right.\nonumber\\
&&-\left.\sum_{k=n_0}^nU_{k,n}^2\left(f\left(y\right)\right)\gamma_k^2Var\left[Z_k\left(y\right)\right]\left[r_{k-1}\left(y\right)-r\left(y\right)\right]^2\right|\mathds{1}_{\sup_{l\leq k-1}\left|Y_l\right|\leq A_n}\nonumber\\
&\leq &u_n^2\sum_{k=n_0}^nU_{k,n}^{2}\left(f\left(x\right)\right)\gamma_k^2Var\left[Z_k\left(x\right)\right]\left|\left(r_{k-1}\left(x\right)-r\left(x\right)\right)^2-\left(r_{k-1}\left(y\right)-r\left(y\right)\right)^2\right|\mathds{1}_{\sup_{l\leq k-1}\left|Y_l\right|\leq A_n}\nonumber\\
&&+u_n^2\sum_{k=n_0}^nU_{k,n}^{2}\left(f\left(x\right)\right)\gamma_k^2\left(r_{k-1}\left(y\right)-r\left(y\right)\right)^2\left|Var\left[Z_k\left(x\right)\right]-Var\left[Z_k\left(y\right)\right]\right|\mathds{1}_{\sup_{l\leq k-1}\left|Y_l\right|\leq A_n}\nonumber\\
&&+u_n^2\sum_{k=n_0}^nU_{k,n}^{2}\left(f\left(y\right)\right)\left|\frac{U_{k,n}^{2}\left(f\left(x\right)\right)}{U_{k,n}^{2}\left(f\left(y\right)\right)}-1\right|\gamma_k^2\left(r_{k-1}\left(y\right)-r\left(y\right)\right)^2Var\left[Z_k\left(y\right)\right]\mathds{1}_{\sup_{l\leq k-1}\left|Y_l\right|\leq A_n}\nonumber\\
&\leq &c^*_{18}u_n^2\sum_{k=n_0}^nU_{k,n}^{2}\left(\varphi\right)\gamma_k^2\left(h_k^{-1}\right)
\left(k^3\gamma_k^2h_k^{-4}\delta_nA_n^2\right)
+c^*_{19}u_n^2\sum_{k=n_0}^nU_{k,n}^{2}\left(\varphi\right)\gamma_k^2
\left(k^2\gamma_k^2h_k^{-2}A_n^2\right)\left(\delta_nh_k^{-3}\right)\nonumber\\
&&+c^*_{20}u_n^2\sum_{k=n_0}^nU_{k,n}^{2}\left(\varphi\right)\left(\delta_n\gamma_n^{-1}\right)
\gamma_k^2\left(k^2\gamma_k^2h_k^{-2}A_n^2\right)h_k^{-1}\nonumber
\end{eqnarray*}
In view of (\ref{aAn*}) and (\ref{aun12}), it follows that there exist $s_2^*>0$ and 
$\tilde S_n^{(2)}\in{\cal GS}(s_2^*)$ such that 
\begin{eqnarray}\label{aMncr}
\left|\Phi_{c_n\left(x\right)}\left(\lambda_n\left(x\right)\right)<M^{\left(n\right)}>_n\left(x\right)-\Phi_{c_n\left(y\right)}\left(\lambda_n\left(y\right)\right)<M^{\left(n\right)}>_n\left(y\right)\right|
&\leq &
\delta_n \tilde S_n^{(2)}.
\end{eqnarray}
Lemma~\ref{aL:Mnc} follows from the combination of~\eqref{aMn} and~\eqref{aMncr}.


\begin{thebibliography}{99}

\bibitem{blum}  Blum, J.R. (1954),
Multidimensional stochastic approximation methods,
{\em Ann. Math. Statist.,} {\bf 25},  737-744.


\bibitem{Bo73}
 Bojanic, R. and Seneta, E. (1973), A unified theory of regularly varying sequences, \textit{Math}. Z., {\bf 134}, 91-106.

\bibitem{chen}  Chen, H. (1988), 
Lower rate of convergence for locating a maximum of a function,
{\em Ann. Statist.,} {\bf 16},  1330-1334.

\bibitem{2chen}  Chen, H.F., Duncan, T.E., and Pasik-Duncan, B. (1999), 
A Kiefer-Wolfowitz algorithm with randomized differences,
{\em IEEE Trans. Automat. Control,} {\bf 44},  442-453.

\bibitem{claeskens} 
Claeskens, G. and Van Keilegom, I.  (2003), Bootstrap confidence bands for regression curves and their derivatives,  {\em Ann. Statist.,}  {\bf 31}, 1852-1884. 


\bibitem{delyon}  Delyon, B. and Juditsky, A.B. (1992),
Stochastic optimization with averaging of trajectories,
{\em Stochastics Stochastic Rep.,} {\bf 39},  107-118.

\bibitem{1dippon}  Dippon, J. and Renz J. (1996),
Weighted means of processes in stochastic approximation,
{\em Math. Meth. Statist.,} {\bf 5},  32-60.

\bibitem{2dippon}  Dippon, J. and Renz J. (1997),
Weighted means in stochastic approximation of minima,
{\em SIAM J. Control Optim.,} {\bf 35},  1811-1827.

\bibitem{dippon}  Dippon, J. (2003),
Accelerated randomized stochastic optimization,
{\em Ann. Statist.,} {\bf 31},  1260-1281.

\bibitem{duflo} Duflo, M. (1996), {\em Algorithmes stochastiques,} 
Collection math{\'e}matiques et applications, Springer.


\bibitem{Du97} 
Duflo, M. (1997), \textit{Random Iterative Models}, Collection Applications of mathematics, Springer.

\bibitem{2fabian} Fabian, V. (1967),
Stochastic approximation of minima with improved asymptotic speed,
{\em Ann. Math. Statist.,} {\bf 38},  191-200.

\bibitem{Ga73}
 Galambos, J. and Seneta, E. (1973), Regularly varying sequences, \textit{Proc. Amer. Math. Soc}., $\mathbf{41}$, 110-116.

\bibitem{Gy02}
Gy\"orfi, L., Kohler, M., Krzyzak, A. and Walk, H. (2002), \textit{A distribution-free theory of nonparametric regression,} Springer-Verlag, New York. 


\bibitem{1Hall80}  Hall, P. and Heyde, C.C. (1980), 
{\em Martingale limit theory and its application,} Academic Press, Inc., 
New York-London.

\bibitem{hall1992} Hall, P. (1992), Effect of bias estimation on coverage accuracy of bootstrap confidence intervals for a probability density, {\em Ann. Statist.} {\bf 20}, 675-694.

\bibitem{kiefer} Kiefer, J. and Wolfowitz, J. (1952),
Stochastic estimation of the maximum of a regression function,
{\em Ann. Math. Statist.} {\bf 23}, 462-466.

\bibitem{kushner} K{ushner}  H.J. and C{lark},  D.S. (1978),
{\em Stochastic approximation methods for constrained and unconstrained systems,} 
Springer, New York.

\bibitem{biskushner} K{ushner}  H.J. and Y{ang},  J. (1993),
{Stochastic approximation with averaging: of the iterates: Optimal 
asymptotic rate of convergence for general processes} 
{\em  SIAM J. Control Optim.,} {\bf 31}, 1045-1062.

\bibitem{1lebreton} Le Breton, A.  (1993),
{About the averaging approach schemes for stochastic approximation,}
{\em  Math. Methods Statist.,} {\bf 2}, 295-315.

\bibitem{lebreton} Le Breton, A. and  Novikov, A. (1995),
{Some results about averaging in stochastic approximation,}
{\em  Metrika,} {\bf 42}, 153-171.

\bibitem{ljung} L{jung}, L.  P{flug}, G. and  W{alk}, H. (1992),  
{\em Stochastic approximation and optimization of random systems,} 
Birkh{\"a}user.

\bibitem{mokka} Mokkadem, A. and  Pelletier, M. (2007),
{A companion for the Kiefer-Wolfowitz-Blum stochastic 
approximation algorithm,} {\em Ann. Statist.,} {\bf 35},  1749-1772.

\bibitem{MP06} 
Mokkadem, A. and Pelletier, M. (2008), Compact law of the iterated logarithm for matrix-normalized sums of random vectors, {\em  to Theory Probab. Appl.,}   {\bf 52}, 636-650.

\bibitem{mok07} Mokkadem, A., Pelletier, M. and Slaoui, Y. (2008),
{The stochastic approximation method for the estimation of a multivariate probability density,} 
{\em  J. Statist. Plann. Inference,} doi:10.1016/j.jspi.2008.11.012


\bibitem{nadarayar} Nadaraya, E. A. (1964), On estimating regression. 
\textit{Theory Probab. Appl.} \textbf{10},  186-190.


\bibitem{Neumann}
Neumann, M. H. and Polzehl, J. (1998), Simultaneous bootstrap confidence bands in nonparametric regression, \textit{J. Nonparametr. Statist.}  \textbf{9}  307-333.

\bibitem{nevelson} Nevels'on, M.B. and  Has'minskii , R.Z. (1976), 
{\em Stochastic approximation and recursive estimation,} 
Amer. Math. Soc, Providence, RI.

\bibitem{2pelletier} Pelletier, M. (2000), 
Asymptotic almost sure efficiency of averaged stochastic algorithms,
{\em  SIAM J. Control Optim.,} {\bf 39}, 49-72.

\bibitem{1polyak}  Polyak, B.T. (1990),
New method of stochastic approximation type,
{\em Automat. Remote Control,} {\bf 51},  937-946.

\bibitem{2polyak}  Polyak, B.T. and Juditsky, A.B. (1992),
Acceleration of stochastic approximation by averaging,
{\em SIAM J. Control Optim.,} {\bf 30},  838-855.

\bibitem{polyak}  Polyak, B.T. and Tsybakov, A.B. (1990),
Optimal orders of accuracy for search algorithms of stochastic 
optimization,
{\em Problems Inform. Transmission,} {\bf 26},  126-133.
\bibitem{RP73}

R\'ev\'esz, P. (1973), Robbins-Monro procedure in a Hilbert space and its application in the theory of learning processes I. \textit{Studia Sci. Math. Hung.}, {\bf 8}, 391-398.

\bibitem{RP}
R\'ev\'esz, P. (1977), How to apply the method of stochastic approximation in the non-parametric estimation of a regression function. \textit{Math. Operationsforsch. Statist., Ser. Statistics}, {\bf 8}, 119-126.

\bibitem{ruppert}  R{uppert}, D. (1982),
Almost sure approximations to the Robbins-Monro and Kiefer-Wolfowitz 
processes with dependent noise,
{\em Ann. of Probab.,} {\bf 10},  178-187.

\bibitem{2ruppert}  R{uppert}, D. (1988),
Efficient estimators from a slowly converging Robbins-Monro process,
{\em Tech. Rep. No.781, School of Oper. Res. and Ind. Engrg.,
Cornell University, Ithaca, NY}.
(See also $\S$ 2.8 of R{uppert}, D. (1991), 
{\em Stochastic approximation,}
in Handbook of Sequential Analysis, Ghosh B.K. and Sen P.K., eds., 
Marcel Dekker, New York, 503-529).

\bibitem{spall} Spall, J.C. (1988),
A stochastic approximation algorithm for large-dimensional systems in the 
Kiefer-Wolfowitz setting. In 
{\em Proc. Conference on Decision and Control,}
1544-1548. IEEE, New York.

\bibitem{2spall} Spall, J.C. (1997),
A one-measurement form of simultaneous perturbation stochastic 
approximation, {\em  Automatica J. IFAC,} {\bf 33},  109-112.

\bibitem{walk}
Walk, H. (2001),  Strong universal pointwise consistency of recursive regression estimates, {\em  
Ann. Inst. Statist. Math.,}  {\bf 53}, 691-707. 

\bibitem{watsonr} Watson, G. S. (1964), Smooth regression analysis. \textit{Sankhya Ser. A.} \textbf{26},  359-372.

\bibitem{yin} Yin, G. (1991),
On extensions of Polyak's averaging approach to stochastic approximation,
{\em Stochastics Stochastic Rep.,} {\bf 33},  245-264.



\end{thebibliography}
\end{document}